\documentclass[12pt]{article}
\usepackage{latexsym}
\usepackage{amsfonts}
\usepackage{amssymb}
\parskip=\medskipamount                 
\thicklines                     
\newcommand{\bimn}[7]{\bibitem{#1}#2,
{\em #3},
{ #4}\hspace{0.25em}{\bf
#5}\hspace{0.25em}(#6)\hspace{0.25em}{#7}.}

\def\inbar{\vrule height1.5ex width.4pt depth0pt}
\def\IC{\relax\,\hbox{$\inbar\kern-.3em{\rm C}$}}
\def\IN{\relax{\rm I\kern-.18em N}}
\def\IQ{\relax\,\hbox{$\inbar\kern-.3em{\rm Q}$}}
\def\IR{\relax{\rm I\kern-.18em R}}
\def\ZZ{\relax{\sf Z\kern-.4em Z}}
\def\a{\alpha} \def\b{\beta}   \def\e{\epsilon} \def\g{\gamma}
\def\G{\Gamma} \def\l{\lambda} 
   
 \def\cB{{\cal B}} \def\cC{{\cal C}} \def\cD{{\cal D}}
   
 \def\cK{{\cal K}} \def\cL{{\cal L}} \def\cM{{\cal M}}
   
 \def\cS{{\cal S}} 

\newtheorem{theorem}{Theorem}[section]
\newtheorem{proposition}[theorem]{Proposition}
\newtheorem{corollary}[theorem]{Corollary}
\newtheorem{conjecture}[theorem]{Conjecture}
\newtheorem{lemma}[theorem]{Lemma}
\newtheorem{definition}[theorem]{Definition}
\newtheorem{remark}[theorem]{Remark}
\newtheorem{mt}{Main Theorem}

\catcode`\@=11

\newif\if@fewtab\@fewtabtrue

\catcode`\@=11

\newif\if@fewtab\@fewtabtrue

{\count255=\time\divide\count255 by 60
\xdef\hourmin{\number\count255}
\multiply\count255 by-60\advance\count255 by\time
\xdef\hourmin{\hourmin:\ifnum\count255<10 0\fi\the\count255}}
\def\ps@draft{\let\@mkboth\@gobbletwo
    \def\@oddhead{}
    \def\@oddfoot
      {\hbox to 7 cm{\footnotesize {\em Draft of \jobname:} \draftdate
       \hfil}\hskip -7cm\hfil\rm\thepage \hfil}
    \def\@evenhead{}\let\@evenfoot\@oddfoot}

\def\ceqno{\global\@fewtabfalse
    \ifcase\@eqcnt \def\@tempa{& & &}\or \def\@tempa{& &}
      \or \def\@tempa{&}
      \or\def\@tempa{}\fi\@tempa
{\rm(\theequation)}}

\def\aeqno#1{\global\@fewtabfalse
    \ifcase\@eqcnt \def\@tempa{& & &}\or \def\@tempa{& &}
      \or \def\@tempa{&}
      \or\def\@tempa{}\fi\@tempa
{\rm(\theequation,#1)}}

\def\label#1{\ifnum\draftcontrol=1
 \global\def\draftnote{$\scriptstyle #1$}\fi
 \@bsphack\if@filesw {\let\thepage\relax
   \def\protect{\noexpand\noexpand\noexpand}%
\xdef\@gtempa{\write\@auxout{\string
      \newlabel{#1}{{\@currentlabel}{\thepage}}}}}\@gtempa
   \if@nobreak \ifvmode\nobreak\fi\fi\fi
  \@esphack}

\def\alabel#1#2{\label{#1}\global\@fewtabfalse
    \ifcase\@eqcnt \def\@tempa{& & &}\or \def\@tempa{& &}
      \or \def\@tempa{&}
      \or\def\@tempa{}\fi\@tempa
{\hbox to 3cm{\phantom{\rm(\theequation,#2)}
\draftnote \hfil}\hskip -3cm {\rm(\theequation,#2)}}}

\def\clabel#1{\label{#1}\global\@fewtabfalse
    \ifcase\@eqcnt \def\@tempa{& & &}\or \def\@tempa{& &}
      \or \def\@tempa{&}
      \or\def\@tempa{}\fi\@tempa
{\hbox to 3cm{\phantom{\rm(\theequation)}
\draftnote \hfil}\hskip -3cm{\rm(\theequation)}}}

\def\eqnarray{\def\draftnote{{}}\global\@fewtabtrue
\stepcounter{equation}\let\@currentlabel=\theequation
\global\@eqnswtrue
\global\@eqcnt\z@\tabskip\@centering\let\\=\@eqncr
$$\halign to \displaywidth\bgroup\@eqnsel\hskip\@centering\@eqcnt\z@
  $\displaystyle\tabskip\z@{##}$&\global\@eqcnt\@ne
  \hskip 1\arraycolsep \hfil$\displaystyle{##}$\hfil
  &\global\@eqcnt\tw@ \hskip 1\arraycolsep
$\displaystyle\tabskip\z@{##}$
\hfil  \tabskip\@centering&\global\@eqcnt\thr@@\llap{##}\tabskip\z@
\cr}

\def\endeqnarray{\@@eqncr\egroup
      \global\advance\c@equation\m@ne$$\global\@ignoretrue}

\def\@eqnnum{\hbox to 3cm{\phantom{\rm(\theequation)} \draftnote
                         \hfil}\hskip -3cm {\rm(\theequation)}}

\def\@@eqncr{\let\@tempa\relax
    \ifcase\@eqcnt \def\@tempa{& & &}\or \def\@tempa{& &}
      \or \def\@tempa{&}
      \or\def\@tempa{}
\fi\@tempa
\if@eqnsw
\if@fewtab\@eqnnum\fi
\stepcounter{equation}\fi\global
\@eqnswtrue\global\@eqcnt\z@\global\@fewtabtrue\cr}

\def\draftcite#1{\ifnum\draftcontrol=1#1\else{}\fi}

\def\@lbibitem[#1]#2{\item{}\hskip -3cm \hbox to 2cm
{\hfil$\scriptstyle\draftcite{#2}$}\hskip
1cm[\@biblabel{#1}]\if@filesw
     {\def\protect##1{\string ##1\space}\immediate
      \write\@auxout{\string\bibcite{#2}{#1}}}\fi\ignorespaces}

\def\@bibitem#1{\item\hskip -3cm \hbox to 2cm
{\hfil $\scriptstyle\draftcite{#1}$}\hskip 1cm
\if@filesw \immediate\write\@auxout
       {\string\bibcite{#1}{\the\value{\@listctr}}}\fi\ignorespaces}

\def\nsection#1{\section{#1}\setcounter{equation}{0}}

\def\draftdate{\number\month/\number\day/\number\year\ \ \ \hourmin }

\global\def\draftcontrol{0}
\catcode`\@=12

\def\theequation{{\thesection.\arabic{equation}}}

\def\qq{\begin{eqnarray}}
\def\qqq{\end{eqnarray}}
\def\rx#1{~(\ref{#1})}
\def\ex#1{eq.\hspace*{-3pt}\rx{#1}}
\def\eex#1{eqs.\hspace*{-3pt}\rx{#1}}
\def\cx#1{~\cite{#1}}
\def\rw#1{~\ref{#1}}

\newlength{\shiftwidth}
\def\shift#1{&&\hbox to \shiftwidth{\hfill $\displaystyle#1$}}
\newlength{\sshiftwidth}
\def\sshift#1{\lefteqn{\hbox to
\sshiftwidth{\hfill$\displaystyle#1$}}}

\def\egt{{\it e.g.\ }}

\def\rhs{{\it r.h.s.\ }}
\def\lhs{{\it l.h.s.\ }}

\def\Tr{\mathop{{\rm Tr}}\nolimits}
\def\Vol{\mathop{{\rm Vol}}\nolimits}

\def\deg{ \mathop{{\rm deg}}\nolimits }
\def\p{^{\prime}}

\def\lk{\mathop{{\rm lk}}\nolimits}
\def\max{\mathop{{\rm max}}\nolimits}

\def\diag#1{\mathop{{\rm diag}}(#1)}

\def\mer{ \mathop{{\rm mer}} }

\def\Pexp{\mathop{{\rm Pexp}}\nolimits}

\def\PexpAm#1{ \Pexp \lrbc{ \oint_{#1} A\,dx} }
\def\PexpAmv#1#2{ \Pexp \lrbc{ \oint_{#1} #2\,dx} }

\def\pr#1#2{ \noindent{\em Proof of #1~\ref{#2}.} }
\def\proof{ \noindent{\em Proof.} }

\def\qed{ \hfill $\Box$ }
\def\const{ {\mbox{const}} }

\def\lrbc#1{ \left( #1 \right) }
\def\lrbs#1{ \left[ #1 \right] }

\def\absol#1{ \left| #1 \right| }

\def\qbezier{\bezier{120}}
\def\DottedCircle{
\bezier{4}(0.966,-0.259)(1.04,0)(0.966,0.259)
\bezier{4}(0.966,0.259)(0.897,0.518)(0.707,0.707)
\bezier{4}(0.707,0.707)(0.518,0.897)(0.259,0.966)
\bezier{4}(0.259,0.966)(0,1.04)(-0.259,0.966)
\bezier{4}(-0.259,0.966)(-0.518,0.897)(-0.707,0.707)
\bezier{4}(-0.707,0.707)(-0.897,0.518)(-0.966,0.259)
\bezier{4}(-0.966,0.259)(-1.04,0)(-0.966,-0.259)
\bezier{4}(-0.966,-0.259)(-0.897,-0.518)(-0.707,-0.707)
\bezier{4}(-0.707,-0.707)(-0.518,-0.897)(-0.259,-0.966)
\bezier{4}(-0.259,-0.966)(0,-1.04)(0.259,-0.966)
\bezier{4}(0.259,-0.966)(0.518,-0.897)(0.707,-0.707)
\bezier{4}(0.707,-0.707)(0.897,-0.518)(0.966,-0.259)
}
\def\Endpoint[#1]{
\ifcase#1
\put(1,0){\circle*{0.15}}
\or\put(0.866,0.5){\circle*{0.15}}
\or\put(0.5,0.866){\circle*{0.15}}
\or\put(0,1){\circle*{0.15}}
\or\put(-0.5,0.866){\circle*{0.15}}
\or\put(-0.866,0.5){\circle*{0.15}}
\or\put(-1,0){\circle*{0.15}}
\or\put(-0.866,-0.5){\circle*{0.15}}
\or\put(-0.5,-0.866){\circle*{0.15}}
\or\put(0,-1){\circle*{0.15}}
\or\put(0.5,-0.866){\circle*{0.15}}
\or\put(0.866,-0.5){\circle*{0.15}}
\fi}
\def\Arc[#1]{
\thicklines			
\ifcase#1
\bezier{25}(0.966,-0.259)(1.04,0)(0.966,0.259)
\or
\bezier{25}(0.966,0.259)(0.897,0.518)(0.707,0.707)
\or
\bezier{25}(0.707,0.707)(0.518,0.897)(0.259,0.966)
\or
\bezier{25}(0.259,0.966)(0,1.04)(-0.259,0.966)
\or
\bezier{25}(-0.259,0.966)(-0.518,0.897)(-0.707,0.707)
\or
\bezier{25}(-0.707,0.707)(-0.897,0.518)(-0.966,0.259)
\or
\bezier{25}(-0.966,0.259)(-1.04,0)(-0.966,-0.259)
\or
\bezier{25}(-0.966,-0.259)(-0.897,-0.518)(-0.707,-0.707)
\or
\bezier{25}(-0.707,-0.707)(-0.518,-0.897)(-0.259,-0.966)
\or
\bezier{25}(-0.259,-0.966)(0,-1.04)(0.259,-0.966)
\or
\bezier{25}(0.259,-0.966)(0.518,-0.897)(0.707,-0.707)
\or
\bezier{25}(0.707,-0.707)(0.897,-0.518)(0.966,-0.259)
\fi}
\def\DottedArc[#1]{
\ifcase#1
\bezier{4}(0.966,-0.259)(1.04,0)(0.966,0.259)
\or
\bezier{4}(0.966,0.259)(0.897,0.518)(0.707,0.707)
\or
\bezier{4}(0.707,0.707)(0.518,0.897)(0.259,0.966)
\or
\bezier{4}(0.259,0.966)(0,1.04)(-0.259,0.966)
\or
\bezier{4}(-0.259,0.966)(-0.518,0.897)(-0.707,0.707)
\or
\bezier{4}(-0.707,0.707)(-0.897,0.518)(-0.966,0.259)
\or
\bezier{4}(-0.966,0.259)(-1.04,0)(-0.966,-0.259)
\or
\bezier{4}(-0.966,-0.259)(-0.897,-0.518)(-0.707,-0.707)
\or
\bezier{4}(-0.707,-0.707)(-0.518,-0.897)(-0.259,-0.966)
\or
\bezier{4}(-0.259,-0.966)(0,-1.04)(0.259,-0.966)
\or
\bezier{4}(0.259,-0.966)(0.518,-0.897)(0.707,-0.707)
\or
\bezier{4}(0.707,-0.707)(0.897,-0.518)(0.966,-0.259)
\fi}
\def\Chord[#1,#2]{
\thinlines
\ifnum#1>#2\Chord[#2,#1]
\else\ifnum#1<#2
\ifcase#1
\ifcase#2
\or\qbezier(1,0)(0.516,0.138)(0.866,0.5)
\or\qbezier(1,0)(0.45,0.26)(0.5,0.866)
\or\qbezier(1,0)(0.327,0.327)(0,1)
\or\qbezier(1,0)(0.179,0.311)(-0.5,0.866)
\or\qbezier(1,0)(0.0536,0.2)(-0.866,0.5)
\or\put(1, 0){\line(-2, 0){2}}
\or\qbezier(1,0)(0.0536,-0.2)(-0.866,-0.5)
\or\qbezier(1,0)(0.179,-0.311)(-0.5,-0.866)
\or\qbezier(1,0)(0.327,-0.327)(0,-1)
\or\qbezier(1,0)(0.45,-0.26)(0.5,-0.866)
\or\qbezier(1,0)(0.516,-0.138)(0.866,-0.5)
\fi
\or\ifcase#2\or
\or\qbezier(0.866,0.5)(0.378,0.378)(0.5,0.866)
\or\qbezier(0.866,0.5)(0.26,0.45)(0,1)
\or\qbezier(0.866,0.5)(0.12,0.446)(-0.5,0.866)
\or\qbezier(0.866,0.5)(0,0.359)(-0.866,0.5)
\or\qbezier(0.866,0.5)(-0.0536,0.2)(-1,0)
\or\put(0.866, 0.5){\line(-5, -3){1.73}}
\or\qbezier(0.866,0.5)(0.146,-0.146)(-0.5,-0.866)
\or\qbezier(0.866,0.5)(0.311,-0.179)(0,-1)
\or\qbezier(0.866,0.5)(0.446,-0.12)(0.5,-0.866)
\or\qbezier(0.866,0.5)(0.52,0)(0.866,-0.5)
\fi
\or\ifcase#2\or\or
\or\qbezier(0.5,0.866)(0.138,0.516)(0,1)
\or\qbezier(0.5,0.866)(0,0.52)(-0.5,0.866)
\or\qbezier(0.5,0.866)(-0.12,0.446)(-0.866,0.5)
\or\qbezier(0.5,0.866)(-0.179,0.311)(-1,0)
\or\qbezier(0.5,0.866)(-0.146,0.146)(-0.866,-0.5)
\or\put(0.5, 0.866){\line(-3, -5){1}}
\or\qbezier(0.5,0.866)(0.2,-0.0536)(0,-1)
\or\qbezier(0.5,0.866)(0.359,0)(0.5,-0.866)
\or\qbezier(0.5,0.866)(0.446,0.12)(0.866,-0.5)
\fi
\or\ifcase#2\or\or\or
\or\qbezier(0,1.)(-0.138,0.516)(-0.5,0.866)
\or\qbezier(0,1.)(-0.26,0.45)(-0.866,0.5)
\or\qbezier(0,1.)(-0.327,0.327)(-1,0)
\or\qbezier(0,1.)(-0.311,0.179)(-0.866,-0.5)
\or\qbezier(0,1.)(-0.2,0.0536)(-0.5,-0.866)
\or\put(0, 1){\line(0, -2){2}}
\or\qbezier(0,1.)(0.2,0.0536)(0.5,-0.866)
\or\qbezier(0,1.)(0.311,0.179)(0.866,-0.5)
\fi
\or\ifcase#2\or\or\or\or
\or\qbezier(-0.5,0.866)(-0.378,0.378)(-0.866,0.5)
\or\qbezier(-0.5,0.866)(-0.45,0.26)(-1,0)
\or\qbezier(-0.5,0.866)(-0.446,0.12)(-0.866,-0.5)
\or\qbezier(-0.5,0.866)(-0.359,0)(-0.5,-0.866)
\or\qbezier(-0.5,0.866)(-0.2,-0.0536)(0,-1)
\or\put(-0.5, 0.866){\line(3, -5){1}}
\or\qbezier(-0.5,0.866)(0.146,0.146)(0.866,-0.5)
\fi
\or\ifcase#2\or\or\or\or\or
\or\qbezier(-0.866,0.5)(-0.516,0.138)(-1,0)
\or\qbezier(-0.866,0.5)(-0.52,0)(-0.866,-0.5)
\or\qbezier(-0.866,0.5)(-0.446,-0.12)(-0.5,-0.866)
\or\qbezier(-0.866,0.5)(-0.311,-0.179)(0,-1)
\or\qbezier(-0.866,0.5)(-0.146,-0.146)(0.5,-0.866)
\or\put(-0.866, 0.5){\line(5, -3){1.73}}
\fi
\or\ifcase#2\or\or\or\or\or\or
\or\qbezier(-1,0)(-0.516,-0.138)(-0.866,-0.5)
\or\qbezier(-1,0)(-0.45,-0.26)(-0.5,-0.866)
\or\qbezier(-1,0)(-0.327,-0.327)(0,-1)
\or\qbezier(-1,0)(-0.179,-0.311)(0.5,-0.866)
\or\qbezier(-1,0)(-0.0536,-0.2)(0.866,-0.5)
\fi
\or\ifcase#2\or\or\or\or\or\or\or
\or\qbezier(-0.866,-0.5)(-0.378,-0.378)(-0.5,-0.866)
\or\qbezier(-0.866,-0.5)(-0.26,-0.45)(0,-1)
\or\qbezier(-0.866,-0.5)(-0.12,-0.446)(0.5,-0.866)
\or\qbezier(-0.866,-0.5)(0,-0.359)(0.866,-0.5)
\fi
\or\ifcase#2\or\or\or\or\or\or\or\or
\or\qbezier(-0.5,-0.866)(-0.138,-0.516)(0,-1)
\or\qbezier(-0.5,-0.866)(0,-0.52)(0.5,-0.866)
\or\qbezier(-0.5,-0.866)(0.12,-0.446)(0.866,-0.5)
\fi
\or\ifcase#2\or\or\or\or\or\or\or\or\or
\or\qbezier(0,-1.)(0.138,-0.516)(0.5,-0.866)
\or\qbezier(0,-1.)(0.26,-0.45)(0.866,-0.5)
\fi
\or\ifcase#2\or\or\or\or\or\or\or\or\or\or
\or\qbezier(0.5,-0.866)(0.378,-0.378)(0.866,-0.5)
\fi\fi\fi\fi}
\def\FullChord[#1,#2]{
\Endpoint[#1]
\Endpoint[#2]
\Arc[#1]
\Arc[#2]
\Chord[#1,#2]
}
\def\EndChord[#1,#2]{
\Endpoint[#1]
\Endpoint[#2]
\Chord[#1,#2]
}
\def\Picture#1{
\begin{picture}(2,1)(-1,-0.167)
#1
\end{picture}
}
\def\DottedChordDiagram[#1,#2]{
\Picture{\DottedCircle \FullChord[#1,#2]}
}
\def\ZZ{ \mathbb{Z} }
\def\IQ{ \mathbb{Q} }
\def\IC{ \mathbb{C} }

\def\u#1{ \underline{#1} }
\def\uu#1{ \underline{\underline{#1}} }

\def\ux{ {\u{x}} }
\def\uy{ {\u{y}} }
\def\ua{ {\u{a}} }
\def\ual{ {\u{\a}} }
\def\ube{ {\u{\b}} }
\def\ug{ {\u{g}} }
\def\um{ {\u{m}} }
\def\ut{ {\u{t}} }

\def\umu{ {\u{\mu}} }
\def\ul{ {\u{\lambda}} }
\def\ue{ {\u{\e}} }
\def\ull{ {\u{l}} }
\def\ug{ {\u{g}} }
\def\un{ {\u{n}} }
\def\uM{ {\u{M}} }

\def\uN{ {\u{N}} }
\def\un{ {\u{n}} }

\def\uc{ {\u{c}} }
\def\ull{ {\u{l}} }

\def\ue{ {\u{\e}} }
\def\uue{ {\uu{\e}} }

\def\uul{ {\uu{\l}} }

\def\uun{ {\uu{n}} }

\def\chR{ \check{R} }
\def\chRi{ \chR^{-1} }

\def\ualp{ \ual\p }

\def\suq{ SU_q(2) }

\def\snzi{ \sum_{n=0}^\infty }
\def\slzi{ \sum_{l=0}^\infty }
\def\smzi{ \sum_{m=0}^\infty }
\def\sjoA{ \sum_{j=1}^A }
\def\sinej{ \sum_{i:\; i\neq j} }

\def\sjoL{ \sum_{j=1}^L }
\def\sjoLp{ \sum_{j=1}^{L\p} }
\def\sjoLmo{ \sum_{j=1}^{L-1} }

\def\sjzn{ \sum_{j=0}^n }
\def\sjzi{ \sum_{j=0}^\infty }
\def\sjoLni{ \sum_{1\leq j\leq L \atop j\neq i} }

\def\sjkzi{ \sum_{j,k=0}^\infty }
\def\skzi{ \sum_{k=0}^\infty }
\def\snkzi{ \sum_{n,k=0}^\infty }

\def\soiljL{ \sum_{1\leq i<j\leq L} }
\def\smu{ \sum_{\mu = \pm 1} }
\def\smumu{ \smu \mu }
\def\smuz{ \sum_{\muz=\pm 1} }
\def\smuzmu{ \smuz \muz }
\def\smuL{ \sum_{\muL=\pm 1} }
\def\smuLmu{ \smuL \muL }

\def\skoL{ \sum_{k=1}^L }
\def\sumuL{ \sum_{{\mu_j=\pm 1\atop 2\leq j\leq L}\atop \mu_1=1} }
\def\sumuLb{ \sum_{\mu_j=\pm 1\atop 1\leq j\leq L} }
\def\sumuLum{ \sumuL \prb{\umu} }
\def\sumuLumb{ \sumuLb \prb{\umu} }
\def\smmnp{ \sum_{\um,\um\p,n\geq 0\atop |\um|+|\um\p|\leq n} }
\def\smnzi{ \sum_{m,n=0}^\infty }
\def\smzn{ \sum_{m=0}^n }

\def\sjoimo{ \sum_{j=1}^{i-1} }
\def\sjioN{ \sum_{j=i+1}^N }
\def\soijL{ \sum_{1\leq i,j\leq L} }
\def\soijLp{ \sum_{1\leq i,j\leq L\p} }
\def\soijLLp{ \sum_{1\leq i\leq L\atop 1\leq j\leq L\p} }

\def\snzgi#1#2{ \sum_{0\leq n \leq \max\{m_{#1},\gr_{#2} - m_{#2} -
           1\} } }
\def\snzmi#1{ \sum_{0\leq n \leq m_{#1} } }

\def\sngz{ \sum_{n\geq 0} }

\def\suntM{ \sum_{0\leq \uun \leq 2M\atop
    \absol{\ung} = \absol{\unl} } }
\def\suneqn{ \sum_{{\uun\geq 0 \atop |\ung|=|\unl|} \atop
      |\unl| + |\unpp| \leq n} }
\def\sngn{ \sum_{n\geq |\unl| + |\unpp|} }
\def\suntMlim{ \sum_{0\leq \uun \leq 2M} }

\def\skgz{ \sum_{k\geq n} }
\def\snkzi{ \sum_{k\geq 0\atop 0\leq n\leq k} }

\def\sumzi{ \sum_{\um\geq 0} }
\def\sjoimo{ \sum_{j=1}^{i-1} }

\def\snzk{ \sum_{n=0}^k }
\def\sioL{ \sum_{i=1}^L }

\def\smobj{ \sum_{m_j=0}^{\b_j-1} }

\def\smotzi{ \sum_{m_1,m_2 =0}^\infty }
\def\snmotzi{ \sum_{m_1,m_2,n =0}^\infty }
\def\snmot{ \sum_{m_1,m_2,n} }

\def\sumn{ \sum_n }

\def\snzi{ \sum_{n=0}^\infty }

\def\soknjL{ \sum_{1\leq k\leq L\atop k\neq j} }

\def\empt{ \varnothing }

\def\pioNg{ \sum_{i=1}^{\Ng} }
\def\pjoA{ \prod_{j=1}^A }
\def\pjoL{ \prod_{j=1}^L }
\def\pioL{ \prod_{i=1}^L }
\def\pjoLmo{ \prod_{j=1}^{L-1} }
\def\pjoLi{ \prod_{1\leq j\leq L\atop j\neq i} }
\def\pjoLp{ \prod_{j=1}^{L\p} }
\def\plon{ \prod_{l=1}^n }
\def\plzjmo{ \prod_{l=0}^{j-1} }

\def\pjoLni{ \prod_{1\leq j\leq L\atop j\neq i} }

\def\tppioL#1{ \lrbc{\pioL t_i^{l_{ij}}}^{#1/2} }
\def\ppioL#1{ \lrbc{\pioL t_i^{l_{ij}}}^{#1} }

\def\tpioL#1{ \lrbc{\pioL t_i^{l_{ij}}}^{#1} }

\def\lila{ \uul }

\def\Jbas#1#2{ J_{#1}(#2;q) }
\def\Jbasl#1#2#3{ J_{#1}(#2;q\,|\,#3) }
\def\JaK{ \Jbas{\a}{\cK} }
\def\JuaL{ \Jbas{\ual}{\cL} }
\def\JuapL{ \Jbas{\ual\p}{\cL} }
\def\JubLp{ \Jbas{\ube}{\cLp} }

\def\JuaubLLp{ \Jbas{\ual,\ube}{\cLLp} }

\def\JabLpl{ \Jbasl{\ual,\ube}{\cLLp}{\lila} }

\def\JabLplM{ \Jbasl{\ual,\ube}{\cLLp}{\lila;M} }

\def\Jr{ J^{\rm (r)} }
\def\Jrbas#1#2#3#4{ \Jr_{#1;#2}(#3,#4;q) }
\def\Jrbass#1#2{ \Jr_{#1}(#2;q) }
\def\JruaubLLp{ \Jrbas{\ual}{\ube}{\cL}{\cLp} }

\def\JruapL{ \Jrbass{\ualp}{\cL} }
\def\JruappL{ \Jrbass{\uapp}{\cL} }

\def\Jrmab{ \Jrbas{\umu\ual}{\ube}{\cL}{\cLp} }
\def\Jrbasl#1#2#3{ \Jr_{#1}(#2;q\,|\,#3) }
\def\JrabLplM{ \Jrbasl{\ual;\ube}{\cL,\cLp}{\lila;M} }

\def\JrabunLplM{ \Jrbasl{\ual;\ube\,|\,\uun}{\cL,\cLp}{\lila;M} }
\def\JrmabLplM{ \Jrbasl{\umu\ual;\ube}{\cLLp}{\lila;M} }

\def\JrLabLplMj{ \Jrbasl{\ual\xrem{j};\b_0,\ube}
           {\cL\xrem{j},\cL_j\cup \cLp}{\lila;M} }
\def\JrmLabLplMj{ \lrbar{\Jrbasl{\ual;\ube}
           {\cL,\cLp}{\lila;M}}{\a_j=\mu\b_0} }

\def\Jhr{ \check{J}^{\rm (r)} }
\def\Jhrbas#1#2#3#4{ \Jhr_{#1}(#2,#3;#4;h) }
\def\JhrubLLputh{ \Jhrbas{\ube}{\cL}{\cLp}{\ut} }

\def\Jhremp{ \Jhr_{\ube}(\empt,\cLp;h) }
\def\JhrLut{ \Jhrbas{}{\cL}{\empt}{\ut} }

\def\AP{ \Delta_{\rm A} }
\def\APbas#1#2{ \AP(#1;#2) }
\def\APpbas#1#2#3{ \AP^{#1}(#2;#3) }
\def\APKt{ \APbas{\cK}{t} }

\def\AF{ \nabla_{\rm A} }
\def\AFbas#1#2{ \AF(#1;#2) }

\def\AFLut{ \AFbas{\cL}{\ut} }

\def\AFLuqap{ \AFbas{\cL}{q^{\ua\p}} }
\def\AFLuto{ \AFbas{\cL}{t_1} }

\def\AFL#1{ \AFbas{\cL}{#1} }
\def\AFKt{ \AFbas{\cK}{t} }
\def\AFLLp#1{ \AFbas{\cLLp}{#1} }
\def\AFLLput{ \AFLLp{\ut} }

\def\Pbas#1#2#3#4#5{ P_{#1;#2}(#3,#4;#5) }
\def\pbas#1#2#3#4#5{ p_{#1;#2}(#3,#4;#5) }
\def\Ppbas#1#2#3#4#5{ P\p_{#1;#2}(#3,#4;#5) }
\def\Ppbase#1#2#3{ P\p_{#1}(#2;#3) }
\def\PubenLLp#1{ \Pbas{\ube}{n}{\cL}{\cLp}{#1} }
\def\PubenLLput{ \PubenLLp{\ut} }
\def\PpubenLLp#1{ \Ppbas{\ube}{n}{\cL}{\cLp}{#1} }
\def\PpubenLLput{ \PpubenLLp{\ut} }
\def\Pubumn#1{ \pbas{\um}{n}{\cL}{\cLp}{#1} }

\def\tP{ \tilde{P} }
\def\tp{ \tilde{p} }
\def\tPbas#1#2#3#4#5{ \tP_{#1;#2}(#3,#4;#5) }
\def\tPubenLLp#1{ \tPbas{\ube}{n}{\cL}{\cLp}{#1} }
\def\tPubezLLput{ \tPbas{\ube}{0}{\cL}{\cLp}{\ut} }
\def\tPubenLLput{ \tPubenLLp{\ut} }

\def\tpbassi#1#2#3#4#5{ \tp_{#1;#2}(#3,#4;#5) }
\def\tPubenumLLput{ \tpbassi{\um}{n}{\cL}{\cLp}{\ube} }

\def\Pbass#1#2#3{ P_{#1;#2}(#3) }
\def\PubenLp{ \Pbass{\ube}{n}{\cLp} }
\def\Ppbass#1#2#3{ P\p_{#1;#2}(#3) }
\def\PpubenLp{ \Ppbass{\ube}{n}{\cLp} }
\def\tPbass#1#2#3{ \tP_{#1;#2}(#3) }
\def\tPubenLp{ \tPbass{\ube}{n}{\cLp} }

\def\frPbas#1#2#3#4#5{ {\Pbas{#1}{n}{#2}{#3}{#4} \over
      #5^{2n+1}(#2;#4)} }

\def\frPbasx#1#2#3#4#5{ {\Pbas{#1}{n}{#2}{#3}{#4} \over
      #5^{2n}(#2;#4)} }
\def\frpPbasx#1#2#3#4#5{ {\Ppbas{#1}{n}{#2}{#3}{#4} \over
      #5^{2n}(#2;#4)} }
\def\frpPbasxe#1#2#3{ {\Ppbase{n}{#1}{#2} \over
      #3^{2n}(#1;#2)} }

\def\frPubLLptoAPx{ \frPbasx{\ube}{\cL}{\cLp}{t_1}{\AP} }
\def\frpPubLLptoAPx{ \frpPbasx{\ube}{\cL}{\cLp}{t_1}{\AP} }
\def\frPubLLputAFx{ \frPbasx{\ube}{\cL}{\cLp}{\ut}{\AF} }
\def\frpPubLLputAFx{ \frpPbasx{\ube}{\cL}{\cLp}{\ut}{\AF} }
\def\frpPubLLputAFxeap{ \frpPbasxe{\cL}{q^{\ua\p}}{\AF} }
\def\frpPubLLputAFxe{ \frpPbasxe{\cL}{\ut}{\AF} }

\def\frPubLLpeaoAP{ \frPbas{\ube}{\cL}{\cLp}{e^{2\pi ia_1}}{\AP} }

\def\frPubLLpemuaAF{ \frPbas{\ube}{\cL}{\cLp}{e^{2\pi i\umu\ua}}{\AF}}

\def\dsbas#1#2#3{ d_{#1;#2}(#3) }

\def\duzzL{ d_{\u{0},0}(\cL) }
\def\dmmnL{ \dsbas{\um,\um\p}{n}{\cLLp} }

\def\Zbas#1#2{ Z_{#1}(#2;K) }
\def\Zsbas#1{ Z(#1;K) }
\def\ZuaML{ \Zbas{\ual}{M,\cL} }
\def\ZuaSL{ \Zbas{\ual}{S^3,\cL} }
\def\ZtruaSL{ \Zibas{\triv}{\ual}{S^3,\cL} }
\def\ZS{ \Zsbas{S^3} }
\def\Zibas#1#2#3{ Z^{(#1)}_{#2}(#3;K) }
\def\ZcuaML{ \Zibas{c}{\ual}{M,\cL} }
\def\ZuaubMLLp{ \Zbas{\ual,\ube}{M,\cLLp} }

\def\Zab{ \Zbas{\ual,\ube}{S^3,\cLLp} }
\def\Ztrab{ \Zibas{\triv}{\ual,\ube}{S^3,\cLLp} }
\def\Zmba{ \Zibas{\umu}{\ube}{S^3,\cL,\cLp;\ua} }
\def\ZmbaK{ \Zibas{\umu}{\ube}{S^3,\cL,\cLp;\ua/K} }

\def\ZubLLpua{ \Zbas{\ube}{S^3,\cL,\cLp;\ua} }
\def\ZubMLLpua{ \Zbas{\ube}{M,\cL,\cLp;\ua} }
\def\ZubMLLpuaK{ \Zbas{\ube}{M,\cL,\cLp;\ual/K} }

\def\triv{ {\rm tr} }
\def\intbcbas#1{ \int\limits_{ {\rm b.c.}(#1)} }
\def\intbcLua{ \intbcbas{\cL;\ua} }

\def\SCS{ S_{\rm CS} }
\def\SCSibas#1{ \SCS^{(#1)} }
\def\SCSc{ \SCSibas{c} }
\def\cMc{ \cM_{\ci} }
\def\ci{ c }
\def\sumc{ \sum_{\ci} }

\def\Nc{ N_{\ci} }
\def\NcM{ \Nc(M) }
\def\Ccuan{ C^{(\ci)}_{\ual;n} }
\def\CcuanML{ \Ccuan(M,\cL) }

\def\Ac{ A^{(\ci)} }
\def\Amu{ A^\bumu }
\def\Dc{ D_{\ci} }
\def\Hc{ H_{\ci} }

\def\RRS{ \tau_{\rm R} }
\def\RRSc{ \RRS^{(\ci)} }

\def\bMo{ b_M^1 }

\def\Phibas#1#2{ \Phi_{\rm A}(#1;#2) }
\def\PL#1{ \Phibas{\cL}{#1} }
\def\PLut{ \PL{\ut} }

\def\PLuqap{ \PL{q^{\ua\p}} }
\def\PLuqapp{ \PL{q^{\uapp}} }

\def\uapp{ \ua^{\prime\prime} }
\def\app{ \a^{\prime\prime} }

\def\Phibasi#1#2#3{ \Phi_{#1}(#2;#3) }
\def\PLo#1{ \Phibasi{1}{\cBL}{#1} }
\def\PLout{ \PLo{\ut} }
\def\Phibasii#1#2#3#4{ \Phi^{#1}_{#2}(#3;#4) }
\def\PLouti#1{ \Phibasii{#1}{1}{\cBL}{\ut} }

\def\sqi#1{ #1,#1^{-1} }

\def\ZZtti{ \ZZ[t,t^{-1}] }
\def\ZZqqi{ \ZZ[\sqi{q}] }
\def\ZZttoi{ \ZZ[t_1,t_1^{-1}] }
\def\ZZutthi{ \ZZ[\ut^{1/2},\ut^{-1/2}] }
\def\ZZutti{ \ZZ[\ut,\ut^{-1}] }
\def\ZZutthih{ \ZZ[\ut^{1/2},\ut^{-1/2},1/2] }

\def\ZZtotq{ \ZZ[\sqi{t_1},\sqi{t_2},\sqi{q}] }

\def\ZZutqpm{ \ZZ[\lila,\sqi{\ut},\sqi{q}] }
\def\ZZulutpm{ \ZZ[\lila,\sqi{\ut}] }
\def\ZZutulz{ \ZZ[\sqi{\ut},\lzl] }

\def\IQhhbas#1{ \IQ[#1][[h]] }
\def\IQmomtn{ \IQ[m_1,m_2,n] }
\def\IQulutpm{ \IQ[\lila,\sqi{\ut}] }

\def\IQbj{ \IQ[\b_j] }
\def\IQutuul{ \IQ[\sqi{\ut},\uul] }
\def\IQube{ \IQ[\ube] }

\def\IQbua{ \IQ[[\ua]] }
\def\IQbuab{ \IQ[\ube][[\ua]] }

\def\ICN{ \IC^N }
\def\ICNo{ \IC^{|\uN|+1} }
\def\ICz{ \IC[z] }
\def\ICzozN{ \IC[z_1,\ldots,z_N] }
\def\ICzzzN{ \IC[z_0,\ldots,z_{|\uN|}] }

\def\mtibas#1#2{ (#1_1,\ldots,#1_{#2}) }
\def\prb#1{ \left\{ #1 \right\} }

\def\dmnK{ d_{m;n}(\cK) }
\def\qmon{ (q-1)^n }
\def\fbas#1#2{ f_{#1}(\cK;#2) }
\def\PnKbas#1{ P_n(\cK;#1) }
\def\PnKt{ \PnKbas{t} }
\def\frKAPbas#1{ {\PnKbas{#1}\over \APpbas{2n+1}{\cK}{#1}} }

\def\cLp{ \cL\p }
\def\cLb{ \bar{\cL} }
\def\cLj{ \cL_j }
\def\cLpj{ \cL\p_j}

\def\cLLp{ \cL\cup\cLp }

\def\lij{ l_{ij} }
\def\lpij{ l\p_{ij} }
\def\lpp{ l^{\prime\prime} }
\def\lppij{ \lpp_{ij} }
\def\cLbas#1{ \cL_{[#1]} }

\def\phbas#1#2{ #1^{#2/2} - #1^{-#2/2} }
\def\phq#1{ \phbas{q}{#1} }

\def\fm{ f_m }
\def\fmmo{ f_{m-1} }
\def\fmpo{ f_{m+1} }
\def\fmom{ \ftm{m_1}{m_2} }
\def\fmomz{ \ftm{m_0}{m_1} }
\def\fmomn#1#2{ \ftm{m_2#1 n}{m_1#2 n} }
\def\fmnmo#1#2{ \ftm{m_1#1 n}{m_2#2 n} }
\def\fmnmoz#1#2{ \ftm{m_0#1 n}{m_1#2 n} }

\def\Vgr{ V_\gr }
\def\gri{ {\gr,\infty} }
\def\grip{ {\gr,\pm\infty} }
\def\Vgri{ V_\gri }
\def\Vgrip{ V_\grip }
\def\Vmgri{ V_{-\gr,\infty} }

\def\Vaj{ V_{\a_j} }
\def\Vaji{ V_{\a_j,\infty} }
\def\Vai{ V_{\a,\infty} }
\def\Vaii#1{ V_{\a_{#1},\infty} }
\def\Vbiip#1{ V_{\b_{#1},\pm\infty} }
\def\Vtaa#1#2{ V_{\gr_{#1}}\otimes V_{\gr_{#2}} }
\def\Vtaai#1#2{ V_{\gr_{#1,\infty}}\otimes V_{\gr_{#2,\infty}} }
\def\Vtaaip#1#2{ V_{\gr_{#1,\pm\infty}}\otimes
           V_{\gr_{#2,\pm\infty}} }
\def\Vtgg#1#2{ V_{\g_{#1}}\otimes V_{\g_{#2}} }
\def\VBua{ V_{\cB;\ual} } 
\def\VBLua{ V_{\cBL;\ual} }
\def\VBLab{ V_{\cBLLp;\ual,\ube} } 

\def\Vtipi#1{ (\Vtip)^{[#1]} }

\def\VBLpb{ V_{\cBLp;\ube} }
\def\VBLai{ V_{\cBL;\ual,\infty} }
\def\VBLaiz{ \Vzi\otimes\VBLai }
\def\VBLpbip{ V_{\cBLp;\ube,\pm\infty} }
\def\Vti{ \VBLai\otimes\VBLpb }
\def\Vtip{ \VBLai\otimes\VBLpbip }

\def\Vgriia#1{ V_{\a_{#1},\infty} }

\def\VBLpjubip{ V_{\cBLpj;\b_j,\pm\infty} }

\def\cBLpj{ \cB_{\cLp_j} }

\def\Vzi{ V_{0,\infty} }
\def\Vzi{ V_{1,\infty} }
\def\lst{ \l_* }

\def\Vot#1{ \Aot{V}{#1} }
\def\brot#1#2{ \lrbc{#1}^{\otimes #2} }
\def\Aot#1#2{ #1^{\otimes #2} }

\def\Vbj{ V_{\b_j} }

\def\Vaj{ V_{\a_j} }

\def\Vas{ V_{\as} }

\def\VBLabo{ \VBLab^{[1]} }

\def\Vbj{ V_{\b_j} }
\def\Vbjip{ V_{\b_j,\pm\infty} }

\def\lo{ \l_0 }
\def\lopa{ \l_0^{\as-1} }
\def\lomH{ \lo^{-H} }
\def\loz{ \l_0=0 }

\def\otidot{ \otimes\cdots\otimes }

\def\ftm#1#2{ f_{#1}\otimes f_{#2} }

\def\qa{ q^\a }
\def\qua{ q^{\ual} }

\def\qao{ q^{\alo} }

\def\muz{ \mu_0 }

\def\aj{ \a_j }
\def\bj{ \b_j }
\def\qaj{ q^{\aj} }
\def\qai#1{ q^{\a_{#1}} }

\def\qaoat#1{ q^{#1(\gr_1-1)(\gr_2-1)/4} }
\def\qaoati#1#2#3{ q^{#1(#2_1-1)(#3_2-1)/4} }

\def\qhaH#1#2{ q^{#1{3\over 8} (\gh-1)\otimes H #2 {1\over 8}H\otimes
          (\gh-1)} }

\def\gh{ \hat{\gr} }

\def\muz{ \mu_0 }

\def\atmto{\gr_2 - m_2 - 1}

\def\Qq{ Q_q }

\def\Qqo{ Q_q^{[1]} }
\def\Qqom{ Q_{q,*}^{[1]} }

\def\qHt{ q^{H/2} }

\def\lmH{ \l^{-H} }

\def\fsl{ \phi_{\rm sl} }
\def\fslLLp{ \fsl(\cLLp) }
\def\qmfsl{ q^{-\fsl} }
\def\qmfslLLp{ q^{-\fsl(\cLLp)} }

\def\fo{ \phi_0 }
\def\foL{ \fo(\cL) }
\def\fz{ \phi_0 }
\def\fzL{ \fz(\cL) }
\def\qfzL{ q^{\fzL} }
\def\ovhqfzL{ {1\over h}\,\qfzL }

\def\ftwo{ \phi_2 }
\def\fth{ \phi_3 }

\def\Phub{ \Phi_{\ube} }

\def\Labul{ \Lii{\ual,\ube}{\lzl} }

\def\qual{ q^{\ual} }
\def\qmual{ q^{-\ual} }

\def\lkbas#1#2{ \lk(#1;#2) }
\def\lkLua{ \lkbas{\cL}{\ual} }
\def\lkLuap{ \lkbas{\cL}{\ualp} }
\def\lkLuapp{ \lkbas{\cL}{\uapp} }

\def\qrtr{ {1\over 4} }
\def\hlf{ {1\over 2} }

\def\lrbar#1#2{ \left. #1 \right|_{#2} }

\def\bumu{ {(\umu)} }
\def\prmu{ \prb{\umu} }

\def\Hmu{ H_\bumu }
\def\Nmu{ N^\bumu }

\def\Ki{ K^{-1} }

\def\serc#1#2{ \lrbs{#1}_{(#2)} }

\def\spb#1#2#3#4{ \plon\lrbc{ {(q^{#1}-q^{#2})(q^{#3}-1)\over
        q^{#4}-1} } }

\def\hnn{ \hlf n(n+1) }

\def\gr{ \gamma }

\def\Qpulj{ Q\p_{\lzl;j} }

\def\hQ{ \hat{Q} }

\def\Qul{ Q_{\lila} }
\def\Quli#1{ Q^{(#1)}_{\lzl} }
\def\QulN{ \Quli{N} }
\def\lzl{ \ul }
\def\lila{ \lzl }

\def\Qulm{ Q_{\lila,*} }

\def\Qulomi#1{ Q^{(#1)}_{\lzl,*} }
\def\QulomN{ \Qulomi{N} }
\def\Quloms{ Q_{\lstar;\lzl} }

\def\Qulmstar{ (\lstar^{2\hm}\otimes\QulomN) }
\def\hQulm{ \hQ_{\lstar,\lzl} }

\def\PM{ P_M }
\def\PMt{ \PM\otimes\PM }
\def\PMtb{ (\PMt) }
\def\IPMtb{ (I\otimes\PM) }
\def\PMItb{ (\PM\otimes I) }

\def\BM{ B_M }
\def\BMLLp{ B_{\cLLp;M} }

\def\BMLLpm{ B_{\cLLp;M,*} }

\def\BMLLpun{ B_{\cLLp;M;\uun} }

\def\BMLLpunps{ \BMLLpun[\ut;\uue] }
\def\BC{ B^{(C)} }
\def\BCMLLpunps{ \BMLLpun^{(C)} [\ut;\uue] }

\def\hB{ \hat{B} }

\def\BCMLunv{ \BC_{\cL;\uun}[\ut;\uue] }
\def\BMLunv{ B_{\cL}[\ut;\ue] }
\def\hBMLunv{ \hB_{\cL}[\ut;\ue] }
\def\BMLjpunv{ B_{\cLpj;M;\uun}[\ut;\uep] }

\def\Dsue{ D^*_{\uun}[\ut;\ube;\del_\uue;h] }
\def\Dsue{ D_{\uun}(\ut;\ube;\del_\uue;h) }
\def\Dsnue{ D_{\uun;n}(\ut;\ube;\uue) }
\def\Dsnued{ D_{\uun;n}(\ut;\ube;\del_\uue) }

\def\Dsnmue{ d_{\uun;n;\um} }

\def\Duulh{ D(\l;\del_{\l};h) }
\def\Dtot{ D_{\uun}(\ut;\l;\ube;\del_{\uue};\del_{\l};h) }
\def\Dntot{ D_{\uun;n}(\ut;\l;\ube;\del_{\uue};\del_{\l}) }

\def\Duunn { D_{\uun;n} }

\def\uep{ \ue\p }

\def\Tjki#1{ T_{j,k}^{(#1)} }
\def\Tjkp{ \Tjki{+} }
\def\Tjkm{ \Tjki{-} }
\def\Tjkpm{ \Tjki{\pm} }
\def\Tzzpm{ T^{(\pm)}_{0,0} }
\def\momtn{ (m_1,m_2,n) }

\def\Tnki#1{ T_{n,k}^{(#1)} }

\def\Tnkop{ \Tnki{1,+} }
\def\Tnkom{ \Tnki{1,-} }
\def\Tnkopm{ \Tnki{1,\pm} }
\def\momtgz{ m_1,m_2,\gr }

\def\Tzzpm#1{ T_{0,0}^{(#1,\pm)} }
\def\Tnkl#1{ T_{n,k,l}^{(2,#1)} }

\def\Tnktp{ \Tnki{2,+} }
\def\Tnktm{ \Tnki{2,-} }
\def\Tnktpm{ \Tnki{2,\pm} }

\def\momgt{ m_1,m_2,\gr,t }
\def\momgti{ m_1,m_2,\gr,\sqi{t} }

\def\Tnkthp{ \Tnki{3,+} }
\def\Tnkthm{ \Tnki{3,-} }
\def\Tnkthpm{ \Tnki{3,\pm} }
\def\momtgzth{ \gr_1,\gr_2,m_1,m_2 }

\def\Smni#1#2{ S^{(#1)}(#2,q) }

\def\qao{ q^{\a_1} }
\def\qat{ q^{\a_2} }

\def\bqmo#1{ (1-q^{#1}) }
\def\bomq#1{ (q^{#1}-1) }
\def\omq#1{ q^{#1} - 1 }

\def\inm{ ^n_m }
\def\imnum{ ^{n_1, n_2}_{m_1,m_2} }
\def\momt{ _{m_1,m_2} }
\def\momz{ _{m_0,m_1} }
\def\mspmn{ \momt^{m_1-n,m_2+n} }

\def\mspmnz{ \momz^{m_0+n,m_1-n} }
\def\msmpnz{ \momz^{m_0-n,m_1+n} }

\def\mrpmn{ \momt^{m_2+n,m_1-n} }
\def\mrmpn{ \momt^{m_2-n,m_1+n} }
\def\motp{ _{m_1,m_2}^{m_2\p,m_1\p} }
\def\motpe{ _{m_0,m_1,m_2}^{m_0\p,m_2\p,m_1\p} }

\def\dlti#1{ \delta^{m_0\p}_{m_0#1n} }

\def\rem#1{ _{[#1]} }
\def\xrem#1{ _{(#1)} }

\def\ccLLp{\cL,\cLp}
\def\cLLz{\cLz\cup\cL}
\def\cLz{ \cL_0 }
\def\cLzLp{ \cLz\cup\cLp }

\def\AMqi#1#2#3{ A_{#1}(#2;\lila;#3;q|M) }
\def\AMq#1{ A_{\ube}(\ccLLp;\lila;#1;q|M) }

\def\AMb#1{ A_{\ube;n}(\ccLLp;\lila;#1|M) }
\def\AMbi#1#2#3{ A_{#1;n}(#2;\lila;#3|M) }

\def\AMqun#1{ A_{\ube\,|\,\uun}(\ccLLp\lila;#1;q|M) }
\def\AMbun#1{ A_{\ube;n\,|\,\uun}(\ccLLp\lila;#1|M) }

\def\cLLt{ (\cL,\cL) }
\def\cLLpt{ (\cL,\cLp) }
\def\cLpLt{ (\cLp,\cL) }
\def\cLpLpt{ (\cLp,\cLp) }
\def\PLM{ P_{\cL;M} }
\def\PLpM{ P_{\cLp;M} }
\def\PLMs{ P^*_{\cL;M} }
\def\PLLp{ P_{\cLLp;M} }

\def\Mi{ M\rightarrow\infty }
\def\lMi{ \lim_{\Mi} }
\def\hz{ h\rightarrow 0 }

\def\hR{ \hat{R} }

\def\Rpm{ R^{\pm} }
\def\Rp{ R^{+} }
\def\Rm{ R^{-} }

\def\hRpm{ \hR^{\pm} }
\def\hRp{ \hR^{+} }
\def\hRm{ \hR^{-} }

\def\tteee{ [t_1,t_2;\uesp] }

\def\Rl{ R^{(-)} }
\def\Rg{ R^{(+)} }
\def\Rgl{ R^{(\pm)} }
\def\hRl{ \hR^{(-)} }
\def\hRg{ \hR^{(+)} }
\def\hRgl{ \hR^{(\pm)} }

\def\el{ \e }
\def\eel{ [\el_1,\el_2] }
\def\eele#1{t_{#1}; \el_1,\el_2,\e\p }
\def\eg{ \e }
\def\eeg{ [\eg_1,\eg_2] }
\def\eege#1{t_{#1}; \eg_1,\eg_2,\e\p }
\def\egl{ \e^{\pm} }
\def\eegl{ [\egl_1,\egl_2] }

\def\Nl{ N^- }
\def\Ng{ N^+ }
\def\Npp{ N^{\prime\prime} }

\def\nl{ n^- }
\def\ng{ n^+ }
\def\npp{ n^{\prime\prime} }
\def\unl{ \un^- }
\def\ung{ \un^+ }
\def\unpp{ \un^{\prime\prime} }

\def\Rni#1{ R^{#1,n} }
\def\Rnp{ \Rni{+} }
\def\Rnm{ \Rni{-} }
\def\Rnpm{ \Rni{\pm} }

\def\Ri#1{ R_n^{(#1)} }

\def\Rppn#1{ R^{\prime\prime}_{n,#1} }

\def\Oni#1{ O^{#1}_n }
\def\Onp{ \Oni{+} }
\def\Onm{ \Oni{-} }
\def\Onpm{ \Oni{\pm} }
\def\Onmp{ \Oni{\mp} }

\def\Chsn{ C_n }

\def\del{ \partial }
\def\deee{ (\del_{\e_1},\del_{\e_2},\del_{\e_3}) }

\def\clsin{ t_{1,2} = q^{\a_{1,2}}\atop \e=0 }

\def\HLp{ H_{\cLp} }

\def\hm{ \hat{m} }

\def\lstar{ \l_{*} }

\def\Ckn{ C_{k,n} }

\def\Dumn{ d_{\um,\uun,n} }

\def\cBL{ \cB_{\cL} }

\def\cBLLp{ \cB_{\cLLp} }
\def\cBLp{ \cB_{\cLp} }

\def\BL{ B_{\cL} }
\def\BLLp{ B_{\cLLp} }

\def\Fung{ F(\ung) }

\def\Oh{ \hat{O} }
\def\liMi{ M\rightarrow \infty }
\def\limMi{ \lim_{\liMi} }

\def\Qiii#1#2#3#4{ Q^{#1}_{#2}(#3;#4;\lzl) }
\def\Qzutule{ \Qii{0}{\cL}{\ut} }
\def\Qizutule#1{ \Qiii{#1}{0}{\cL}{\ut} }
\def\Qzutule{ Q(\cL;\ut;\lzl) }
\def\Qizutule#1{ \Qiii{#1}{}{\cL}{\ut} }
\def\Qnzutule{ \Qizutule{2n+1} }
\def\Qili#1#2#3{ Q_{#1}(#2;#3;\lzl) }
\def\Qnubut{ \Qili{\ube;n}{\ccLLp}{\ut} }
\def\Qzubut{ \Qili{\ube;0}{\ccLLp}{\ut} }

\def\Qmnubut{ q_{\um;n}(\ccLLp;\ube) }

\def\ulzcl{ {\lo=0\atop\l=1} }
\def\lslo{ {\lstar=0\atop \l=1} }

\def\uer{ {\u{e}} }

\def\rhop{ \rho_+ }
\def\rhom{ \rho_- }
\def\rhopm{ \rho_\pm }
\def\tontw{ [t_1,t_2] }
\def\thtmh{ t_1^{1/2} - t_1^{-1/2} }

\def\fii#1#2{ \phi_{#1}(#2) }
\def\fii#1#2{ \phi(#2\,|\,#1) }
\def\qfii#1#2{ q^{\fii{#1}{#2}} }

\def\Lii#1#2{ \Lambda(#1;#2) }
\def\as{ \a_* }
\def\bas{ [\as] }

\def\Foi#1{ \Phi_2(#1;\ut\rem{1};\lzl;\ube) }
\def\Foxi#1{ \Phi_2(#1;\ut\xrem{i};\lzl;\ube) }
\def\muaL{ \mu_L\a_L }

\def\lepcl{ {\l=1 \atop \lo,\uep=0} }
\def\elt{ \ut,\ul,\uep }

\def\muL{ \mu_L }

\def\bLp{ \b_{L\p} }

\def\bcLs{ \mathfrak{L}_* }

\def\cLro{ \cL\rem{1} }

\def\Ci#1#2{ \cC_{#1;#2} }
\def\CLl{ \Ci{\ccLLp}{l} }

\def\Gbas#1#2#3{ G^{(#1)}(#2;#3) }
\def\Gbasi#1{\Gbas{c_{#1}}{\cL}{\ua} }
\def\GcLua{ \Gbas{c}{\cL}{\ua} }

\def\Fbas#1#2#3#4{ F^{(#1)}_{#2}(#3;#4) }
\def\FcnLab{ \Fbas{c}{n}{\ccLLp}{\ua;\ube} }
\def\Fbasi#1{ \Fbas{#1}{n}{\ccLLp}{\ua;\ube} }

\def\Hbas#1#2#3#4{ F^{(#1)}_{#2}(#3;#4) }
\def\Hcmni#1{ \Hbas{c_{#1}}{m,n}{\ccLLp}{\uaro;\ube} }

\def\Hbas#1#2#3#4{ H^{(#1)}_{#2}(#3;#4) }
\def\Hbasi#1{ \Hbas{c_{#1}}{m_1,m_2,n}{\ccLLp}{\ube} }
\def\Hbasii#1#2{ \Hbas{c_{#1}}{#2}{\ccLLp}{\ube} }

\def\Ibas#1#2#3{ I^{(#1)}(#2;#3) }
\def\Ibasi#1{\Ibas{c_{#1}}{\ccLLp}{\ua;\ube;\e} }
\def\IcLLabe{ \Ibas{c}{\ccLLp}{\ua;\ube;\e} }
\def\Ibasii#1{ I^{(c_{#1})}_{m,n}(\ccLLp;\uaro;\ube) }
\def\Iot{ I^{(c_1,c_2)}_{m,n}(\ccLLp;\uaro;\ube) }
\def\Iotm{ I^{(c_1,c_2)}_{m,n-m}(\ccLLp;\uaro;\ube) }
\def\Ibasw#1#2{ I(#1;#2;\e) }
\def\Ibasws{ \Ibasw{\ccLLp}{\ua;\ube} }

\def\IGbas#1#2{ I_S(#1;#2;\e) }

\def\sg{ \sum_{g\in \GcLLp} }
\def\sgp{ \sg (-1)^{|g|} }

\def\XXbas#1{ X_{#1}^{(c_1,c_2)} }

\def\cLpro{ \cLp\rem{1} }

\def\uaro{ \ua\rem{1} }

\def\mc#1{ \mathcal{#1} }

\def\cB{ \mc{B} }
\def\cC{ \mc{C} }
\def\cD{ \mc{D} }

\def\cK{ \mc{K} }
\def\cL{ \mc{L} }
\def\cM{ \mc{M} }

\def\cS{ \mc{S} }

\def\fPan{ P_{\a;n} }
\def\qnal{ {q^{\a/2} - q^{-\a/2} \over q^{1/2} - q^{-1/2} } }
\def\baas{ [\a] }

\def\URC{$U(1)$-RC\ }

\def\Lmipo{ i }

\def\GcL{ G_{\cL} }
\def\GcLLp{ G_{\ccLLp} }

\def\cLpk{ \cLp_k }

\def\cLkLLz{ \cLpk\cup\cLLz }

\def\bLp{ \b_{L\p} }

\def\Xvar#1{ X^{(#1)} }
\def\xargv{ (t_1,t_2;m_1,m_2,n;q) }

\def\spN{ {N^{\prime}} }

\def\pkoNp{ \prod_{l=1}^{\spN} }
\def\es{ \sigma }
\def\esv#1#2{ \es_{#1}^{(#2)} }

\def\uesp{ \e_1,\e_2,\e_3,\e_4 }
\def\uespi#1{ \e_{1,#1},\e_{2,#1},\e_{3,#1},\e_{4,#1} }

\def\uer{ \u{e} }

\def\phlk{ \phi_{\rm lk} }
\def\phlkbas#1#2{ \phlk(#1;#2) }
\def\phlkLpb{ \phlkbas{\cL,\cLp}{\ube} }
\def\phlkL{ \phlk(L) }

\def\PhubLpt{ \Phub(\cL,\cLp;\ut) }
\def\PhubLpto{ \Phub(\cL,\cLp;t_1) }

\def\cLzL{ \cLz\cup\cL }

\begin{document}

\begin{titlepage}
\centerline{\hfill                 math.QA/9806004}
\vfill
\begin{center}

{\large \bf
A contribution of a $U(1)$-reducible connection to quantum
invariants of links I: $R$-matrix and Burau representation} \\

\bigskip
\centerline{L. Rozansky\footnote{
This work was supported by NSF Grant DMS-9704893}
}


\centerline{\em Department of Mathematics, Yale University}
\centerline{\em 10 Hillhouse Ave., P.O. Box 208283}
\centerline{\em New Haven, CT 06520-8283}
\centerline{\em E-mail address: rozansky@math.yale.edu}

\vfill
{\bf Abstract}

\end{center}
\begin{quotation}

We introduce a new invariant of links in $S^3$ which we call a
contribution of a $U(1)$-reducible connection to the colored Jones
polynomial. This invariant is a formal power series in $h$, whose
coefficients are rational functions of parameters
$t_j$
assigned to the link components. The denominators of the
rational functions are the powers of the colored
Alexander-Conway function. Similarly to the case of a knot, the Jones
polynomial of an algebraically connected link is determined by the
contribution of the $U(1)$-reducible connection. Namely, if we
substitute $t_j= (1+h)^{\a_j}$
in the formal power
series and then expand the resulting expressions in powers of $h$ and
$\a_j$ through a special step-by-step procedure, then we obtain the
expansion of the colored Jones polynomial of the link in powers of
$h=q-1$ and colors $\a_j$. The proof of these results is based on a
similarity between $R$-matrix and Burau representations of the
braid group.




\end{quotation} \vfill \end{titlepage}

\pagebreak
\tableofcontents

\nsection{Introduction}
\label{s1}
\hyphenation{Re-she-ti-khin}
\hyphenation{Tu-ra-ev}
\subsection{Motivation}

The famous topological invariant of knots and links, the
Alexander-Conway polynomial, can be calculated by various means. In
particular, it can be calculated
by Conway's skein relations or by expressing it as a
determinant of a matrix related to the Burau representation of the
braid group. The polynomial that comes out of these calculations is
known to be determined by the group of the knot. Therefore the
Alexander-Conway polynomial has a clear topological interpretation.

The colored Jones polynomial of knots and links can be also
calculated by either skein relations and cabling or by taking a trace
of the product of $R$-matrices coming from the quantum algebra
$\suq$. However, in contrast to the Alexander polynomial, the
emerging expressions do not appear to have a simple structure that
would point towards their topological interpretation. Therefore some
tricks are needed in order to see topology inside the colored Jones
polynomial.

P. Melvin and H. Morton\cx{MeMo} suggested to study the expansion of
the colored Jones polynomial $\JaK$ of a knot $\cK\subset S^3$ in
powers of $q-1$ and color $\a$
\qq
\JaK = \snzi \fPan(\cK)\,(q-1)^n.
\label{1.1*}
\qqq
They proved that the coefficients $\fPan(\cK)$ can be presented as
polynomials of $\a$ and they conjectured the bound on the degree of
these polynomials
\qq
\fPan(\cK) \in \IQ[\a],\qquad \deg \fPan(\cK) \leq n/2.
\label{1.1*1}
\qqq
%
%
Since
\qq
q^\a - 1 = \a\,(q-1) + O((q-1)^2)\qquad\mbox{as $q\rightarrow 1$},
\label{1.2}
\qqq
then
we can rewrite the series\rx{1.1*} as
\qq
\JaK = \baas\snzi \fbas{n}{q^\a}\,\qmon,
\label{1.3}
\qqq
where
\qq
\baas = \qnal
\qqq
and the coefficients
$\fbas{n}{t}$ are formal power series in $t^{1/2}-t^{-1/2}$
\qq
\fbas{n}{t} = \smzi \dmnK\,(t^{1/2}-t^{-1/2})^{2n}.
\label{1.4}
\qqq
(only even powers of $t^{1/2}-t^{-1/2}$ participate in this
expansion, because the polynomials\rx{1.1*1} are known to contain
only odd powers of $\a$ and $[\a]$ is an odd function of $\a$).

Based on their calculations of the expansion\rx{1.3} for torus knots,
Melvin and Morton conjectured that
\qq
\fbas{0}{t} = {1\over\APKt},
\label{1.5}
\qqq
where $\APKt$ is the Alexander-Conway polynomial of $\cK$.
The conjectures of Melvin and Morton were
proved by D.~Bar-Natan and S.~Garoufalidis\cx{BG}.
We showed in\cx{Ro9} that the other coefficients $\fbas{n}{t}$ of
\ex{1.3} are also related to the Alexander-Conway polynomial. We
proved that
\qq
\fbas{n}{t} = \frKAPbas{t},\qquad
\PnKt\in\ZZtti.
\label{1.6}
\qqq
The manner in which the polynomials $\PnKt$ appear in the
expansion\rx{1.3},\rx{1.6} suggests that they might present a natural
extension of the Alexander-Conway polynomial and therefore may have a
direct topological interpretation. Another reason to explore the
polynomials $\PnKt$ is that the expansion
\qq
\JaK = \baas\snzi \frKAPbas{\qa}\,\qmon
\label{1.7}
\qqq
is well-suited for studying $p$-adic properties of the
Reshetikhin-Turaev invariant of rational homology
spheres\cx{Lw},\cx{Ro8}.

The goal of this paper is to extend the expansion\rx{1.7} to the
links in $S^3$. As we will see, the Melvin-Morton expansion of the
colored Jones polynomial of a link can not be
directly presented in the
form\rx{1.7}. However, we will be able to define the analogs of
polynomials $\PnKt$ for links whose Alexander-Conway polynomial is not
identically equal to zero. We will do it in such a way that the link
analog of the relation\rx{1.7} will hold (in a certain way) if the
link has sufficiently many non-zero linking numbers.

\subsection{Notations, definitions and basic properties}
\subsubsection{Multi-index notations}
For a non-negative integer $A$ we denote $
%
\ux = \mtibas{x}{A}
%
$
and for a positive $k$, $\ux\rem{k} = \mtibas{x}{A-k}$ and
$\ux\xrem{k} = (x_1,\ldots,x_{k-1},x_{k+1},\ldots x_A)$.
We also
use the following notations
\qq
\begin{array}{c@{\qquad} c@{\qquad} c}
y\ux  =  \mtibas{yx}{A}, &
\ux^y  =  \mtibas{x^y}{A}, &
\ux\uy  =  (x_1 y_1,\ldots,x_A y_A),\\[2ex]
\prb{f(\ux)}  =  \pjoA f(x_j), &
|\ux|  =  \sjoA x_j, &
\ux^{\uy} = \pjoA x_j^{y_j}.
\end{array}
\label{1.9}
\qqq
For a constant $C$ we denote
\qq
\u{C} = (C,\ldots,C).
\label{1.9*}
\qqq
We assume concatenation of multi-indices
\qq
x,\uy = (x,y_1,\ldots, y_A),\qquad
\ux,\uy = (x_1,\ldots, x_A, y_1,\ldots, y_B)\quad
\mbox{if $\uy = \mtibas{y}{B}$}.
\label{1.9*1}
\qqq
Also
%
$\ux = x$ means $x_j=x$ for all $1\leq j\leq A$.

Let $f(\ux)$ be a meromorphic function of the variables $\ux$ which
may have singularities in the vicinity of $\ux=0$.
$\serc{f(\ux)}{\ux}$ denotes a result of sequential expansion of
$f(\ux)$ in powers of $\ux$ when we start with expanding in $x_L$ and
finish with expanding in $x_1$.

\subsubsection{Linear algebra notations}

Let $V$ be a linear space with a basis $\fm$, $m\geq 0$. If an
operator $O$ acts on $V$, then we denote its matrix elements as
$O\inm$ so that
\qq
O(\fm) = \sum_n O\inm\,f_n.
\label{1.9*1*1}
\qqq
Suppose that $V_1$ and $V_2$ are two linear spaces with bases $\fm$,
$m\geq 0$. Let $O$ be an operator which maps the tensor product
$V_1\otimes V_2$ either into itself or into $V_2\otimes V_1$. We
denote the matrix elements of $O$ as $O\imnum$ so that
\qq
O(\fmom) = \sum_{n_1,n_2} O\imnum \, f_{n_1}\otimes f_{n_2}.
\label{1.9*1*2}
\qqq
If $V$ is a linear space, then for a positive
integer $N$ we denote
\qq
\Vot{N} = \underbrace{V\otidot V}_{N\ \rm times}.
\label{1.9*2}
\qqq
Similarly, if $O$ is an operator acting on $V$, then $\Aot{O}{N}$
denotes a corresponding operator acting on $\Vot{N}$.

Let $O$ be an operator acting on a tensor product of linear
spaces $V\otimes W$. If we take a `partial' trace of $O$ over
the space $W$, then we obtain an operator acting on $V$. We denote it
as $Tr_W O$. Suppose that there is a subspace $U\subset V$
such that $(\Tr_W O) (U)\subset U$. This means that the action of
$\Tr_W O$ on $V$ can be restricted to $U$. We denote the trace of
$\Tr_W O$ over $U$ simply as $\Tr_{U\otimes W} O$:
\qq
\Tr_{U\otimes W} O = \Tr_U \lrbc{ \Tr_W O}
\label{1.9*3}
\qqq
\begin{remark}
\rm
Note that this notation is a bit misleading because, generally
speaking, the action of $O$ on $U\otimes W$ is not well-defined.
However, we hope that there will be no confusion in the context of
our calculations.
\end{remark}

\subsubsection{Topological notations}
In this paper $\cL$ (or $\cLp$) usually denotes an
$L$ (or $L\p$)-component link in $S^3$, while $\cLj$,
$1\leq j\leq L$ (or $\cLpj$, $1\leq j\leq L\p$)
denote its components.
The numbers $\lij$, $1\leq i,j\leq L$
(or $\lpij$, $1\leq i,j\leq L\p$) are the linking numbers of $\cL$
(or $\cLp$). The numbers $\lpp_{ij}$, denote the linking numbers
between $\cL_i\in\cL$ and $\cLpj\in\cLp$. We denote
\qq
\lkbas{\cL}{\ux} = \hlf\soiljL \lij\,x_i x_j
\label{1.9*4}
\qqq

For $1\leq k\leq L-1$, $\cLbas{k}$ denotes a sublink of
$\cL$ which consists of components $\cLj$, $1\leq j\leq L-k$ and for
$1\leq k \leq L$, $\cL\xrem{k}$ denotes the link $\cL$ with the
component $\cL_k$ removed.

\begin{definition}
\label{t1.1}
\rm
A link $\cL$ is \emph{algebraically connected} to another link $\cLp$
if either one of the links is empty, or at least one of the linking
numbers $\lpp_{ij}$ is non-zero.
\end{definition}
\begin{definition}
A link $\cL$ is called \emph{algebraically connected} if it can not
be presented as a union of two links which are not algebraically
connected to each other.
\end{definition}
%
\begin{definition}
\label{t1.2}
\rm
The indexing of the components of an algebraically connected link
$\cL$ is called \emph{admissible} if $\cLbas{k}$ is algebraically
connected for any $1\leq k\leq L-1$.
\end{definition}
\begin{lemma}
\label{t1.3}
If $\cL$ is an algebraically connected $L$-component link, then for
any selected component of $\cL$ there exists an admissible indexing
such that the selected component becomes $\cL_1$.
\end{lemma}
\proof
We leave the proof to the reader.\qed

Let $\cB$ be an $N$-component braid. Its \emph{closure} is
constructed by connecting the upper ends of its strands with the
lower ends in the same order. For any link $\cL$ there exists a braid
$\cBL$ whose closure is $\cL$. The numbers $\uN=\mtibas{N}{L}$ denote
the numbers of strands of $\cB$ which belong to the link components
$\cLj$ of $\cL$.

In this paper we will always assume that the strands of $\cL$ appear
at the bottom and top of $\cB$ in a particular order: if $i<j$, then
the strands of $\cL_i$ precede the strands of $\cLj$. If we present a
union of two links $\cLLp$ as a closure of a braid $\cBLLp$, then we
also assume that the strands of $\cL$ precede the strands of $\cLp$.
$\uN\p = \mtibas{N}{L\p}$ denote the numbers of strands in the link
components of $\cLp$.

For any link
$\cL$ one can always find a braid $\cBL$ which satisfies the ordering
property. Note, that this property of $\cBL$ is completely unnecessary
for the calculations that we will perform in this paper, but it will
simplify our notations.

We say that an elementary crossings of
$\cBLLp$ is of $\cLLt$, $\cLLpt$,
$\cLpLt$ or $\cLpLpt$ type depending on whether the incoming strands
belong to $\cL$ or $\cLp$.

\subsubsection{Alexander-Conway polynomial}
Let $\cK$ be a knot in $S^3$. $\APKt$ denotes its Alexander-Conway
polynomial. This polynomial satisfies the skein relation of
Fig.\rw{fig:skein}
and it is normalized by the condition
\qq
\APbas{\cK}{1}=1.
\label{1.10}
\qqq
%
\begin{figure}[hbt]
\def\KZero{
\begin{picture}(2,2)(-1,-1)
\put(0,0){\circle{2}}
\put(-0.707,-0.707){\vector(1,1){1.414}}
\put(0.707,-0.707){\vector(-1,1){1.414}}
\put(0,0){\circle*{0.15}}
\end{picture}}

\def\KPlus{ \Picture {
\put(0,0){\circle{2}}
\put(-0.707,-0.707){\vector(1,1){1.414}}
\put(0.707,-0.707){\line(-1,1){0.6}}
\put(-0.107,0.107){\vector(-1,1){0.6}}
} }


\def\KMinus{ \Picture {
\put(0,0){\circle{2}}
\put(-0.707,-0.707){\line(1,1){0.6}}
\put(0.107,0.107){\vector(1,1){0.6}}
\put(0.707,-0.707){\vector(-1,1){1.414}}
} }


\def\KII{ \Picture {
\put(0,0){\circle{2}}
\qbezier(-0.707,-0.707)(0,0)(-0.707,0.707)
\qbezier(0.707,-0.707)(0,0)(0.707,0.707)
\put(-0.607,0.607){\vector(-1,1){0.1414}}
\put(0.607,0.607){\vector(1,1){0.1414}}
} }

\begin{displaymath}
{\Delta_A\left(\KPlus\; ;t\right)}~\qquad - \qquad
\Delta_A\left(\KMinus\; ;t\right)~\qquad = \qquad
(t^{1/2} - t^{-1/2})\;
\Delta_A\left(\KII\; ;t\right)~.
\end{displaymath}
\caption{The skein relation for the Alexander-Conway polynomial}
\label{fig:skein}
\end{figure}

For an \emph{oriented} $L$-component link $\cL\subset S^3$ we
define the Alexander-Conway function $\AFLut$, $\ut=\mtibas{t}{L}$
according to the definition
of\cx{Tu} except that our $\ut$ is equal to $\ut^2$ of Turaev.
The function $\AFLut$ is almost a polynomial:
\qq
\AFLut \in
\left\{
\begin{array}{cl}
{\ZZttoi\over \phbas{t_1}{1} } &\mbox{if $L=1$}\\
\ZZutthi &\mbox{if $L\geq 2$}.
\end{array}
\right.
\label{1.11}
\qqq
It satisfies the following
properties:
\qq
&\AFKt = {\APKt \over \phbas{t}{1} }\qquad\mbox{for $\cL$ being
a knot $\cK$},
\label{1.12}
\\
&\AF(\mbox{Hopf link};t_1,t_2) = 1,
\label{1.13}
\\
&\AFL{\ut^{-1}} = (-1)^L\,\AFLut.
\label{1.14}
\qqq
If we change the orientation of a component $\cL_j$
and denote the new oriented link as $\cL\p$, then
\qq
\AFbas{\cLp}{\ut} = -\AFL{t_1,\ldots,t_j^{-1},\ldots,t_L}.
\label{1.15}
\qqq

The Alexander-Conway function satisfies the Torres formula: if
$L\geq 2$, then for $i$, $1\leq i\leq L$
\qq
\lrbar{\AFLut}{t_i=1} = \lrbc{ \pjoLi t_j^{l_{ij}/2} -
\pjoLi t_j^{-l_{ij}/2} }
\AFbas{\cL\xrem{i}}{\ut\xrem{i}}.
\label{1.15*}
\qqq

A correction factor
\qq
\PLut = \pjoL t_j^{-{1\over 2} \lrbc{ \sinej l_{ij} + 1} }
\label{1.16}
\qqq
allows us to remove the fractional powers of $\ut$ from $\AFLut$: if
$L\geq 2$ then
\qq
\PLut \AFLut \in \ZZutti.
\label{1.17}
\qqq
\subsubsection{Burau representation}
The Alexander-Conway function of a link can be calculated with the
help of Burau representation of the braid group. Let $\cBL$ be an
$N$-strand braid whose closure is a link $\cL$. Consider a linear
space $\ICN$ whose basis vectors $\uer=\mtibas{e}{N}$ are associated
with $N$ positions of the braid strands. To a positive (negative)
elementary
braid in positions $j$, $j+1$
we associate an operator $\rhop\tontw$ ($\rhom\tontw$)
acting on the vectors $e_j$, $e_{j+1}$
\qq
\rhop\tontw = \pmatrix{ 1 - t_2^{-1} & t_1^{-1} \cr 1 & 0},\qquad
\rhom\tontw = \pmatrix{ 0 & 1 \cr t_2 & t_1^{-1} t_2(1-t_1) },
\label{1.17*1}
\qqq
where $t_1$ and $t_2$ denote the variables related to the incoming
link components. Note that $\rhop^{-1}\tontw = \rhom[t_2,t_1]$. If we
present the braid $\cBL$ as a product of elementary braids and then
replace those braids with operators $\rhopm$, then we obtain an
operator $\BL[\ut]$ which acts on $\ICN$. Let $Q$ be an operator
which acts on $\ICN$ as
\qq
Q(e_1) = 0,\qquad Q(e_j)=e_j \quad\mbox{for $j\neq 1$}.
\label{1.17*1*1}
\qqq
Then
\qq
\AFLut = \PLout\,
{\det\nolimits_{\ICN}(1-Q\BL[\ut]) \over t_1^{1/2} - t_1^{-1/2} },
\label{1.17*2}
\qqq
where
\qq
\PLout = t_1^{1/2}\,\pjoL t_j^{\hlf \lrbc{-N_j + \sioL \lij} }.
\label{1.17*3}
\qqq

\subsubsection{$\suq$-modules and $R$-matrix}

The quantum invariants in 3-dimensional topology depend on
a `quantum deformation parameter' $q$. We do not assume that $q$ is a
root of unity unless stated explicitly. We also use a parameter
\qq
h = q-1.
\label{1.18}
\qqq

Let $X$, $Y$ and $H$ be the generators of the quantum
group $\suq$ which satisfy the commutation relations
\qq
[H,X] = 2X,\quad [H,Y] = -2Y, \quad [X,Y] = [H].
\label{1.19}
\qqq
In the last equation of\rx{1.19} we used a standard notation
\qq
[x] = {\phq{x} \over \phq{1}}.
\label{1.20}
\qqq
We choose a basis $\fm$, $0\leq m\leq \gr-1$ of a $\gr$-dimensional
$\suq$ module $\Vgr$ in such a way that
\qq
H\,\fm & = & (\gr - 2m - 1)\, \fm,
\label{1.21}\\
X\,\fm & = & [m]\, \fmmo,
\label{1.22}\\
Y\,\fm & = & [ \gr - m - 1]\,\fmpo,
\label{1.23}
\qqq
We will also use an auxiliary operator $\hm$
\qq
\hm(\fm) = m\fm.
\label{1.23*}
\qqq
%


$R$-matrix acts on a tensor product of two $\suq$ modules
$\Vtaa{1}{2}$
\qq
R(\fmom) & = & \snzgi{1}{2} R\mspmn\fmnmo{-}{+},
\label{1.24**1}
\qqq
%
where
\qq
R\mspmn & = &
(\phq{1})^n \,
{ [\atmto]! \over [\atmto - n]! }\,
{ [m_1]! \over [m_1-n]!\,[n]! }
\label{1.24}\\
&&
\qquad\times
q^{ {1\over 4} \big( (\gr_1-2m_1-1)(\gr_2-2m_2-1) -
n(\gr_1 - \gr_2 -2m_1+2m_2+n+1)\big)}.
\nonumber
\qqq
%

\subsubsection{Colored Jones polynomial}
\label{Jpol}
Let $\cL\subset S^3$ be an $L$-component link which is colored by
positive integers $\ual = \mtibas{\a}{L}$. We present
$\cL$ as a closure of an $N$-strand braid $\cBL$. We index the
strands of $\cBL$ according to the order in which they appear at the
bottom of the braid.


Consider a tensor product of $\suq$ modules `attached' to the braid
strands
\qq
\VBLua = \bigotimes_{j=1}^L \Vot{N_j}_{\a_j}.
\label{1.25}
\qqq
We can construct an action of the braid $\cB$ on $\VBua$ if we
present $\cB$ as a product of positive and negative elementary braids
and replace each positive braid with the operator
\qq
\chR = P\,R, \qquad P\colon
\; \Vtaa{1}{2}  \rightarrow \Vtaa{2}{1},\;\;
P(\ftm{m_1}{m_2}) = \ftm{m_2}{m_1}
\label{1.26}
\qqq
acting on a product of $\suq$ modules assigned to the braided pair of
strands, while replacing each negative elementary braid with the
inverse operator $\chRi$. We denote the resulting operator as
\qq
\BL\colon\VBLua \rightarrow \VBLua.
\label{1.27}
\qqq
%
%
We introduce the `quantum trace' operator
\qq
\Qq = \brot{\qHt}{N}.
\label{1.28}
\qqq
The colored Jones polynomial of $\cL$ is defined as a trace over
$\VBLua$
\qq
\JuaL = \qmfsl\Tr_{\VBLua} \Qq \BL.
\label{1.29}
\qqq
The factor $\qmfsl$ is necessary to correct self-linking of the
components of $\cL$
\qq
\fsl(\cL) = \qrtr \sjoL l_{jj}\,(\a_j^2-1),
\label{1.30}
\qqq
where $l_{jj}$ is the number of positive self-crossings minus the
number of negative self-crossings of the strands of $\cB$ which
constitute $\cL_j$.

Being a polynomial of
$q^{1/2},q^{-1/2}$, the colored Jones polynomial of
$\cL$ can be expanded in powers of $h=q-1$. Melvin and
Morton\cx{MeMo} proved that the coefficients in this expansion are
polynomials of the color variables
\qq
& \JuaL = \snzi \Pbass{\ual}{n}{\cL}\, h^n,
\label{mm1}\\
&
\Pbass{\ual}{n}{\cL} \in \IQ [\ual], \qquad
\deg_{\ual} \Pbass{\ual}{n}{\cL} \leq 2n+1.
\nonumber
\qqq
and the polynomials $\Pbass{\ual}{n}{\cL}$ are odd functions of all
variables $\ual$. We call \ex{mm1} \emph{the Melvin-Morton
expansion}.

\subsection{Results}
The purpose of this paper is to prove the following
\def\themt{}
\begin{mt}
\nonumber
Let $\cL$ be an oriented $L$-component link in $S^3$, whose
Alexander-Conway function is not identically equal to zero
\qq
\AFLut \not\equiv 0,\quad \ut = \mtibas{t}{L}.
\label{1.32}
\qqq
Let $\cLp$ be an
$L\p$-component link in $S^3$ with positive
integers $\ube=\mtibas{\b}{L\p}$
assigned to its components.
Then there exists a unique sequence of polynomial invariants
\qq
\PubenLLput \in \ZZutthih, \quad n\geq 0,
\label{1.33}
\qqq
such that
an expansion of $\PubenLLput$ in powers of $\ut-1$ has the
form
\qq
& \PubenLLput =
\sumzi\Pubumn{\ube}\,(\ut-1)^{\um},
\qquad \Pubumn{\ube} \in \IQube,
\label{1.34}
\qqq
and the formal power series in $h$
%
%
%
\qq
&\JhrubLLputh =
\left\{
\begin{array}{cl}   \displaystyle
 \snzi \PubenLp\,h^n
&\mbox{if $L=0$,}
\vspace{5pt}
\\
\displaystyle
{1\over h}\,{1 \over \AFbas{\cL}{t_1}}
\snzi \frPubLLptoAPx\,h^n
&\mbox{if $L=1$,}
\vspace{5pt}
\\
\displaystyle
{1\over h}\, {1 \over \AFLut}
\snzi \frPubLLputAFx\,h^n
&\mbox{if $L\geq 2$}
\end{array}
\right.
\label{1.35}
\qqq
%
satisfies the following properties:
\begin{itemize}
\item[(1)]
If the link $\cL$ is empty (that is, if $L=0$), then
$\Jhremp$ is the
Melvin-Morton expansion\rx{mm1} of $\JubLp$
\qq
\Jhremp
= \lrbar{\JubLp}{q=1+h}\quad\mbox{in $\IQhhbas{\ube}$.}
\label{1.38}
\qqq
%
\item[(2)]
For any $j$, $1\leq j\leq L\p$,
\qq
\lrbar{\JhrubLLputh}{\b_{j}=1} =
\Jhrbas{\ube\xrem{j}}{\cL}{\cLp\xrem{j}}{\ut}.
\label{1.39}
\qqq
%
\item[(3)]
Let $\cL$ be a non-empty link. Suppose that for a number $i$,
$1\leq i\leq L$, the link component $\cL_i$ is algebraically
connected to $\cL\xrem{i}$ and
$\AFbas{\cL\xrem{i}}{\ut\xrem{i}}\not\equiv 0$. Then
%
\qq
\hspace{-50pt}
\Jhrbas{\b_0,\ube}{\cL\xrem{i}}{\cL_i\cup\cLp}{\ut\xrem{i}} =
\left\{
\begin{array}{cl}
\displaystyle
\smumu
\lrbc{ \pjoLi t_j^{l_{ij}} }^{\mu\b_0/2}
\lrbar{
\Jhrbas{\ube}{\cL}{\cLp}{\ut}
}{t_i = q^{\mu\b_0}}
& \mbox{if $L\geq 2$,}
\vspace{5pt}
\\
\displaystyle
\Jhrbas{\ube}{\cL}{\cLp}{q^{\b_0}}
& \mbox{if $L=1$.}
\end{array}
\right.
\label{1.42}
\qqq
\end{itemize}

Thus defined, the polynomials $\PubenLLput$ have additional
properties:
%
%
%
\qq
&\Pbas{\ube}{0}{\cL}{\cLp}{\ut} = 
\left\{
\begin{array}{cl}
\displaystyle
\{\ube\} & \mbox{if $L=0$,}
\vspace{5pt}
\\
\displaystyle
\pjoLp{ \tppioL{\b_j} - \tppioL{-\b_j} \over \tppioL{1} -
\tppioL{-1}} &
\mbox{if $L\geq 1$,}
\end{array}
\right.
\label{1.43}
\qqq
\qq
&\PubenLLp{\ut^{-1}} = \PubenLLput,
\label{1.44}
\qqq
and if we reverse the orientation of a link component $\cL_j$, then
for the new oriented link $\cLb$
\qq
\Pbas{\ube}{n}{\cLb}{\cLp}{\ut} =
\PubenLLp{t_1,\ldots,t_j^{-1},\ldots,t_L}.
\label{1.45}
\qqq
The polynomials $\Pubumn{\ube}$ contain only odd powers of $\b_j$
\qq
\Pubumn{\b_1,\ldots,-\b_j,\ldots,\b_L} =
- \Pubumn{\ube},
\label{1.45*2}
\qqq
and their degree is bounded
\qq
\deg p_{\um,n} \leq
\left\{
\begin{array}{cl}
|\um| + 2n + 1 & \mbox{if $L=1$,}\\
|\um| + 4n + 1 & \mbox{if $L\geq 1$.}
\end{array}
\right.
\label{1.45*0}
\qqq
Finally, if we select an orientation for
the components of $\cLp$, then
the formal power series\rx{1.35} can be rewritten as
\qq
&\JhrubLLputh =
\left\{
\begin{array}{cl}   \displaystyle
q^{\phlkLpb}\snzi \PpubenLp\,h^n
&\mbox{if $L=0$,}
\vspace{5pt}
\\
\displaystyle
{1\over h}\,q^{
\phlkLpb}\,
{\PhubLpto\over \AFbas{\cL}{t_1}}
\snzi \frpPubLLptoAPx\,h^n
&\mbox{if $L=1$,}
\vspace{5pt}
\\
\displaystyle
{1\over h}\, q^{
\phlkLpb}\,
{\PhubLpt\over \AFLut}
\snzi \frpPubLLputAFx\,h^n
&\mbox{if $L\geq 2$,}
\end{array}
\right.
\label{1.45**0}
\qqq
where
\qq
\lefteqn{
\phlkLpb
}
\label{1.45**1}
\\
&& =
\hlf
\lrbc{
L + \soiljL \lij +
\sum_{1\leq i<j\leq L\p} \lpij
(\b_i-1)(\b_j-1)
-\soijLLp \lpp_{ij}(\b_j-1)
-\sjoLp (\b_j-1)
}
\nonumber
\qqq
\qq
&\PhubLpt  =
\pioL t_i^{\hlf \sjoLp\lpp_{ij}(\b_j-1) },
\label{2.103}
\qqq
and
\qq
\PpubenLLput \in \ZZutti
\label{1.45**2}
\qqq
for any integer $\ube$.

\end{mt}

\begin{remark}
\rm
A combination of \eex{1.44} and\rx{1.45} implies that if we
take a link $\cL$ and reverse the orientation of all its
components, then for the resulting link $\cLb$
\qq
\Pbas{\ube}{n}{\cLb}{\cLp}{\ut} =
\PubenLLput.
\label{1.45*0*1}
\qqq
\end{remark}

The uniqueness of the polynomials\rx{1.33} means that the formal
power series\rx{1.35} is an invariant of a pair of
(oriented and unoriented) links $\cL,\cLp$.
\begin{definition}
\rm
We call the expression
\qq
\JruaubLLp = q^{\lkLua}\,\Jhrbas{\ube}{\cL}{\cLp}{\qual}
\label{1.45**}
\qqq
\emph{the $U(1)$ reducible connection contribution} to
the Jones polynomial, or simply the \emph{\URC invariant}.
\end{definition}

\begin{corollary}
\label{cor1}
Suppose that the link $\cL$ is algebraically connected and that it
has admissible indexing. Let us expand the sum
\qq
\sumuLum \Jrmab
\label{1.45*}
\qqq
in powers of $h$ and $\ual$ going in sequence of
$L$ steps, at
$(L-i+1)$-st step
($1\leq i\leq L$) expanding in powers of $\a_{\Lmipo}$ and $h$
by using the formulas
\qq
q^{a_{\Lmipo}} & = & (1+h)^{\a_{\Lmipo}} = \snzi
{\a_{\Lmipo} \choose n} h^n,
\label{1.40}\\
q^{\a_{\Lmipo} \aj} & = & (1+ (\qaj - 1) )^{\a_{\Lmipo}} =
\snzi
{\a_{\Lmipo} \choose n} (\qaj -1)^n \quad\mbox{for $1\leq j<\Lmipo$.}
\label{1.41}
\qqq
%
The resulting power series
coincides with the Melvin-Morton
expansion of the colored Jones polynomial $\JuaubLLp$
\qq
\serc{\sumuLum\Jrmab}{\ual} =
\JuaubLLp\quad\mbox{in $\IQhhbas{\ual,\ube}$}.
\label{1.45*1}
\qqq
\end{corollary}

\begin{remark}
\label{remord}
\rm
Since the Alexander-Conway function $\AFLLput$ may have zeroes in the
vicinity of $\ut=1$, then the order of expanding $\Jrmab$ in
powers of $\ual$ is important. In fact, the order suggested in
Corollary\rw{cor1} guarantees that no negative powers of $\ual$
appear at any stage in the expansion of the sum over $\umu$.
\end{remark}
\pr{Corollary}{cor1}
Since $\cL$ is algebraically connected, then for any $j$,
$2\leq j\leq L$, the link component $\cLj$ has non-zero linking
number with at least one of the components of $\cLbas{j}$.
Therefore we can apply Claim~3 of Main Theorem at each step of
the expansion of the sum\rx{1.45*}. After we complete the expansion
in powers of $\ual$, we apply Claim~2 which leads us to the claim
of Corollary\rw{cor1}.\qed

\begin{remark}
\rm
Suppose that $q=e^{2\pi i/K}$, where $K$ is a positive integer and
the colors $\ual$ are positive integers between 1 and $K-1$. Then
Equations\rx{1.45*0} and\rx{1.45*1} lead to a simple proof of the
so-called Symmetry Principle established by R.~Kirby and
P.~Melvin\cx{KiMe}, which describes how the colored Jones
polynomial of a link changes when one of the colors is changed
from $\a_j$ to $K-\a_j$. We give a precise formulation and a
detailed proof of the Symmetry Principle in Appendix\rw{c2}.

\end{remark}

\subsection{A quick recipe for calculating \URC invariant}
Our proof of Main Theorem is constructive. In Section\rw{s2} we give
a detailed recipe for calculating the \URC invariant $\JhrubLLputh$.
The presence of the link $\cLp$ makes this calculation more involved,
so here we will explain briefly how to obtain the \URC invariant
$\JhrLut$ in case when $\cLp$ is empty.

First of all, we have to get the polynomials $\Tjkpm \momtn$,
described in Theorem\rw{expand}. These polynomials come from
expanding the following expressions
\qq
\lefteqn{
\Xvar{+}\xargv
}
\label{1.45*1*1}\\
&
\hspace{1in}= &
q^{(m_1+1)m_2 + \hlf n(n+1)}
\plon \lrbc{
(q^{-m_2-l} - t_2^{-1})(q^{m_1-n+l} - 1)\over q^l-1},
\nonumber\\
\lefteqn{
\Xvar{-}\xargv
}
\label{1.45*1*2}\\
&
\hspace{1in} = &
q^{-m_1(m_2+1) - \hlf n(n+1)}
\plon \lrbc{
(q^{m_1+l} - t_1)(q^{-m_2+n-l} - 1) \over q^{-l} - 1}
\nonumber
\qqq
in powers of $h=q-1$ (while keeping all other variables fixed):
\qq
\lefteqn{
\Xvar{+}\xargv
}
\label{1.45*1*3}\\
&
\hspace{1in}=
&
{m_1 \choose n}\, (1-t_2^{-1})^n
\sjzn \lrbc{
{\plzjmo (n-l)\over (1-t_2^{-1})^j} \skzi
\Tjkp\momtn\,h^{j+k}
},
\nonumber\\
\lefteqn{
\Xvar{-}\xargv
}
\label{1.45*1*4}\\
&
\hspace{1in}=
&
{m_2 \choose n}\, (1-t_1)^n
\sjzn \lrbc{
{\plzjmo (n-l)\over (1-t_1)^j} \skzi
\Tjkm\momtn\,h^{j+k}
}.
\nonumber
\qqq
We will also need the coefficients $\Ckn$ which
appear in the expansion
\qq
{1\over n!}\,\lrbc{(1+h)^{-1/2} -1}^n =
\skgz \Ckn\,h^k.
\label{2.69x}
\qqq

Suppose that an $L$-component link $\cL$ can be constructed as a
closure of an $N$-strand braid $\cBL$. Let us present $\cBL$ as a
product of $\spN$ elementary braids
\qq
\cBL = \pkoNp \esv{p(l)}{s(l)},
\label{1.45*1*5}
\qqq
where $\esv{p}{\pm}$ is a positive (negative) elementary braid which
twists the strands in the $p$-th and $(p+1)$-th positions, while
$p(l)$ and $s(l)$ are the functions describing the positions and
signs of the elementary braids which comprise a presentation for
$\cBL$. We are going to map this presentation into an operator acting
on an $N$-dimensional space $\ICN$ with a chosen basis
$\uer = \mtibas{e}{N}$. Let $\hRpm_p\tteee$, $1\leq p\leq N-1$ be the
matrices\rx{2.80},\rx{2.81} acting on the basis elements
$e_p, e_{p+1}$. Then we define the operator
\qq
\hBMLunv = \pkoNp \hR^{s(l)}_{p(l)}
[t_{n_1(l)},t_{n_2(l)},\uespi{l}],
\label{1.45*1*6}
\qqq
where $\ue = (\e_{j,l}\,|\, 1\leq j\leq 4,\, 1\leq l\leq \spN )$
and the functions $n_1(l),n_2(l)$ indicate to which link components
belong the strands which come into the $l$-th elementary braid
$\esv{p(l)}{s(l)}$. We also need an operator
\qq
\hQ = \diag{0,1,\ldots,1},
\label{1.45*1*7}
\qqq
which projects out $e_1$.

Next, we construct a formal power series in $h$ whose coefficients
are polynomials of the variables $\l$ and $\ue$
\qq
D(\l,\ue;h) =
\lrbc{
\snkzi \Ckn\,\l^n\,h^k
}
\pkoNp
\lrbc{
\sjkzi T_{j,k}^{(s(l))}
({\e_{1,l}},{\e_{2,l}},{\e_{3,l}})
\,{\e_{4,l}^j}\,h^{j+k}
}.
\label{1.45*1*8}
\qqq
Now
\qq
\JhrLut & = &
{1- t_1^{-1} \over 1 - q^{-1}}\,
q^{\hlf\sjoL (l_{jj} - N_j)}
\pjoL t_j^{\hlf \lrbc{N_j-\sioL \lij} }
\label{1.45*1*9}\\
&&\hspace{1.5in}\times
\lrbar{
D(\l\del_\l,\del_{\ue};h)
{1 \over \det_{\ICN}(1-\l \hQ \hBMLunv)}
}
{\ue = 0\atop\l=1},
\nonumber
\qqq
where $N_j$ is the number of the strands of the braid $\cBL$ which
are part of the link component $\cL_j$ and $l_{jj}$ is the
self-linking number of $\cL_j$ defined by the closure of $\cBL$ (in
other words, $l_{jj}$ equals the number of positive self-crossings of
the strands of $\cLj$ minus the number of negative self-crossings).

\begin{remark}
\rm
The formula\rx{1.45*1*9} and the detailed calculations of
subsections\rw{s2.2.3} -- \ref{s2.2.5} may seem rather
`unnatural'. The natural way of calculating the \URC invariant
through $R$-matrix formula for the colored Jones polynomial\rx{1.29}
is to use the machinery of free field representation for $\suq$
generators and the perturbative expansion of quantum field theory. In
this paper we tried to avoid any explicit reference to these methods
and stick to simple algebra. However, we hope to explain the
application of QFT methods in a future publication dealing with the
colored Jones polynomials associated to simple Lie algebras of
higher rank.

\end{remark}

\subsection{Path integral interpretation}

E.~Witten showed\cx{Wi} that the colored Jones polynomial of a link
$S^3$ can be expressed in terms of a path integral over $SU(2)$
connections. For a positive integer $K$ and an $L$-component link
$\cL$ in a 3-manifold $M$ consider the following path integral which
goes over connections $A$ in the trivial $SU(2)$ bundle over $M$
\qq
\ZuaML = \int \cD A\;
e^{(iK/2\pi)\SCS[A]}
\pjoL \Tr_{\aj} \PexpAm{\cLj}.
\label{1.46}
\qqq
Here $\ual = \mtibas{\a}{L}$
are positive integers, $\SCS[A]$ is the Chern-Simons invariant
of $A$, $\PexpAm{\cLj}$ is a holonomy of $A$ along $\cLj$ and
$\Tr_{\aj}$ is a trace in the $\aj$-dimensional $SU(2)$ module. If
$M=S^3$ then, according to\cx{Wi},
\qq
\JuaL = {\ZuaSL \over \ZS},\qquad
\ZS = \sqrt{2\over K}\sin\lrbc{\pi\over K},
\label{1.47}
\qqq
assuming that
\qq
q= e^{2\pi i/K}.
\label{1.48}
\qqq

At the `physical level of rigor' the path integral\rx{1.46} can be
evaluated in the stationary phase approximation as
\qq
K\rightarrow \infty,\qquad\ual = \const.
\label{1.49}
\qqq
The stationary points of the Chern-Simons action are flat $SU(2)$
connections. Therefore the integral\rx{1.46} should split into a sum
of contributions of connected components $\cMc$ of the moduli space
$\cM$ of flat $SU(2)$ connections on $M$, each contribution being
proportional to the `classical exponential'
\qq
&\ZuaML = \sumc\ZcuaML,
\label{1.50}\\
&\ZcuaML = e^{(iK/2\pi)\SCSc}
\lrbc{4\pi^2\over K}^{\NcM}\snzi \CcuanML\,K^{-n},
\label{1.51}\\
&\NcM\in\IQ,\qquad \CcuanML\in\IC,
\nonumber
\qqq
Here $\ci$ indexes the connected components of $\cM$ and $\SCSc$ is
the Chern-Simons invariant of any connection belonging to $\cMc$.

If a component $\cMc$ represents a single connection $\Ac$ (that is,
if $\cMc$ is a separate point in $\cM$), then there are specific
`physical' predictions about the numbers $\NcM$ and $\CcuanML$. Let
$\Dc$ be a covariant derivative in the trivial $su(2)$ bundle over
$M$ associated with the connection $\Ac$. $\RRSc$ denotes the
Reidemeister-Ray-Singer torsion of $\Dc$,
$\Hc$ denotes the isotropy group of $\Ac$
($\Hc\subset SU(2)$ commutes with the holonomies of $\Ac$),
$\Vol(\Hc)$ is the volume of $\Hc$ (if $\Hc$ is discrete, then
$\Vol(\Hc)$ is the number of elements, otherwise the volume is
calculated with the metric in which the roots of $su(2)$ have length
$\sqrt{2}$) and $\bMo$ is the first Betti
number of $M$. Then, according to\cx{Wi},
\qq
&\NcM = \hlf\,\dim\Hc,
\label{1.52}\\
&
\absol{\CcuanML} =
{
\pjoL \Tr_{\aj} \PexpAmv{\cLj}{\Ac} \over
\Vol(\Hc)\; (\RRSc)^{1/2} }.
\label{1.53}
\qqq

If $M=S^3$, then the only flat $SU(2)$ connection on $M$ is the
trivial connection, so the sum of \ex{1.50} contains only one term
\qq
\ZuaSL = \ZtruaSL.
\label{1.53*}
\qqq
We use $\triv$ as the index $c$
of the trivial connection. It is easy to
see that
\qq
H_\triv = SU(2), \qquad \Vol(H_\triv) = 2^{5/2} \pi^2,\qquad
\Nc(S^3) = 3/2,\qquad
\RRS^{(\triv)} = 1,
\label{1.54}
\qqq
so that the absolute value of the first term in the sum of \ex{1.51}
is
\qq
\absol{C^{(\triv)}_{\ual;0}(S^3,\cL)}
=  \prb{\ual} 2^{-5/2}\pi^{-2}.
\label{1.55}
\qqq
This is consistent with
the first term of the Melvin-Morton expansion\rx{mm1}. Indeed, since
\qq
h = 2\pi i\Ki + O(K^{-3}),\quad
\ZS = 2^{1/2}\pi K^{-3/2} + O(K^{-7/2})
\quad\mbox{as $K\rightarrow \infty$}
\label{1.55*}
\qqq
then \eex{1.47},\rx{1.51} predict that
\qq
\absol{\duzzL} = {K^{3/2}\over 2^{1/2}\pi}
\lrbc{ 4\pi^2\over K}^{\Nc(S^3)} \absol{C^{(\triv)}_{\ual;0}(S^3,\cL)}
=1,
\label{1.55*1}
\qqq
while according to\cx{MeMo}, $\duzzL=1$.

The paper\cx{EMSS} describes a path integral representation of the
colored Jones polynomial which is slightly different from \ex{1.46}.
Let $\cL,\cLp$ be a pair of links in $M$ endowed with colors
$\ual,\ube$.
For rational numbers
$\ua=\mtibas{a}{L}$ such that $\ua/K\in \ZZ$, consider
the following expression
\qq
\ZubMLLpua = \intbcLua
\cD A\;e^{(iK/2\pi)\SCS[A]}\; \pjoLp \Tr_{\bj} \PexpAm{\cLpj}.
\label{1.56}
\qqq
Here $\intbcLua$ denotes a path integral over $SU(2)$ connections in
the complement of the tubular neighborhood of $\cL$ (we denote this
neighborhood simply as $M\setminus \cL$). The connections should
satisfy the following boundary condition: for each connection there
should exist the elements $\ug \in SU(2)$, $\ug=\mtibas{g}{L}$ such
that
\qq
\PexpAm{\mer(\cLj)} = g_j e^{i\pi a_j \sigma_3} g_j^{-1}.
\label{1.57}
\qqq
%
In \ex{1.57} we used the following notations: $\mer(\cLj)$ is the
meridian of $\cLj$ (that is, a simple cycle on the boundary of the
tubular neighborhood of $\cLj$ which is contractible through that
neighborhood) and $\sigma_3$ is one of the Pauli matrices:
\qq
\sigma_3 = \pmatrix{1 & 0\cr 0 & -1}.
\label{1.59}
\qqq
According to\cx{EMSS}, the integral\rx{1.56} is equal to $\ZuaubMLLp$
at $\ua=\ual/K$
\qq
\ZuaubMLLp =
\ZubMLLpuaK.
\label{1.59*}
\qqq

Representation\rx{1.51} coupled with \ex{1.59}
allows us to apply the stationary phase
approximation to the calculation of $\ZuaubMLLp$ in the limit
\qq
K\rightarrow \infty, \qquad \ua = \ual / K = \const,\qquad
\ube=\const.
\label{1.60}
\qqq
The integral\rx{1.51} should split into a sum of contributions of
flat $SU(2)$ connections on $M\setminus \cL$ which satisfy the
condition\rx{1.52}.

If $M=S^3$, then among all flat connections satisfying\rx{1.57} there
are exactly $2^{L-1}$ $U(1)$-reducible connections. Their holonomies
belong to the same subgroup $U(1)\subset SU(2)$. These
$U(1)$-reducible connections can be indexed by a set of numbers
$\umu = \mtibas{\mu}{L}$ such that $\mu_j=\pm 1$, $2\leq j\leq L$,
$\mu_L=2$. Using $\umu$ as an index of the corresponding
$U(1)$-reducible connection
and identifying $U(1)$ with $\{z\in\IC\,|\;|z|=1\}$ we find that
\qq
\PexpAmv{\mer(\cLj)}{A^\bumu} = e^{i\pi \mu_j a_j}.
\label{1.61}
\qqq
\begin{conjecture}
\label{conj1}
The formal sum\rx{1.35} is related to the contributions of
$U(1)$-reducible connections. Namely,
\qq
{\ZmbaK \over \ZS} = \prmu\Jrmab,
\label{1.62}
\qqq
or, equivalently,
\qq
{\Zmba\over \ZS} =
\left\{
\begin{array}{cl}   \displaystyle
 \snzi \PubenLp\,h^n
&\mbox{if $L=0$,}
\vspace{5pt}
\\
\displaystyle
{\sin(\pi a_1/K) \over \sin(\pi/K)}
\snzi \frPubLLpeaoAP\,h^n
&\mbox{if $L=1$,}
\vspace{5pt}
\\
\displaystyle
{\prb{\umu}\over h}\, q^{\foL}\,
e^{2\pi iK\lkbas{\cL}{\umu\ua}} \snzi \frPubLLpemuaAF\,h^n
&\mbox{if $L\geq 2$.}
\end{array}
\right.
\label{1.62*}
\qqq
\end{conjecture}

We have two reasons to make this conjecture. The first reason is that
the first term in the sum of \ex{1.62*} matches the first term of the
sum of \ex{1.51} written for the $U(1)$-reducible connection $\Amu$.
Indeed, in view of \ex{1.43},
the first term in the expansion of the \rhs of \ex{1.62*}
for $L\geq 1$ in
powers of $\Ki$ becomes
\qq
{\prmu\over 2\pi i K}\,
{e^{2\pi iK\lkbas{\cL}{\umu\ua}}\over \AFL{e^{2\pi i \umu\ua}}}\,
\pjoLp {\sin\lrbc{ \pi\b_j \skoL \mu_k a_k} \over
\sin\lrbc{ \pi\skoL \mu_k a_k} }.
\label{1.63}
\qqq
On the other hand, the Chern-Simons action of $\Amu$ (with
appropriate boundary terms included) is
\qq
\SCS^\bumu = 2\pi^2 \lkbas{\cL}{\umu\ua},
\label{1.64}
\qqq
while the isotropy group is
\qq
\Hmu = U(1),\quad \Vol(\Hmu) = 2^{3/2}\pi.
\label{1.65}
\qqq
The trivial $su(2)$ bundle with connection $\Amu$ splits into a sum
of a 1-dimensional real bundle with zero connection and a
1-dimensional complex bundle with the corresponding $U(1)$
connection. The Reidemeister-Ray-Singer torsion of the former is
equal to 1, while the torsion of the latter is equal to the absolute
value of the Alexander-Conway function $\AFL{e^{2\pi i \umu\ua}}$.
Therefore
\qq
(\RRS^\bumu)^{1/2} = \absol{ \AFL{e^{2\pi i \umu\ua}} }.
\label{1.66}
\qqq
Finally,
\qq
\Tr_{\b_j}\PexpAmv{\cLp_j}{\Amu} =
{\sin\lrbc{ \pi\b_j \skoL \mu_k a_k} \over
\sin\lrbc{ \pi\skoL \mu_k a_k} },
\label{1.67}
\qqq
so according to \ex{1.53},
\qq
\lrbc{ 4\pi^2\over K}^{\Nmu(S^3)}
\absol{C^\bumu_{\ual,\ube;0}(S^3,\cLLp)} =
{(2K)^{-1/2} \over
\absol{ \AFL{e^{2\pi i \umu\ua}} } }
\pjoLp \absol{\sin\lrbc{ \pi\b_j \skoL \mu_k a_k} \over
\sin\lrbc{ \pi\skoL \mu_k a_k} }.
\label{1.68}
\qqq
If we substitute \eex{1.64} and\rx{1.68} into the
expression\rx{1.51}, then we can see easily that the dominant terms
of the $\Ki$ expansion of both sides of \ex{1.62} coincide at least
up to a $K$-independent phase factor which may come from
$C^\bumu_{\ual,\ube;0}(S^3,\cLLp)$.

The second reason for making Conjecture\rw{conj1} is the relation
between $\JruaubLLp$ and $\JuaubLLp$ established in
Corollary\rw{cor1}. Its path integral interpretation relies on the
following
theorem.

\def\Asp#1#2{ {\cal A}_{#1}}
\def\Aspucna{ \Asp{\uc,\un}{a} }

For positive numbers $\uc=\mtibas{c}{L}$ and positive integers
$\un = \mtibas{n}{L-1}$ let $\Aspucna$ denote a set of all
points $\ua=\mtibas{a}{L}$ which satisfy the inequalities
\qq
a_1<c_1, \qquad a_{j+1}<
c_{j+1}\, a_{j}^{n_j},\quad 1\leq j\leq L-1,
\qquad a_L>0.
\label{1.69}
\qqq

\begin{theorem}
\label{red.conn}
Let $\cL\subset S^3$ be an algebraically connected link with
admissible indexing. Then there exist positive numbers
$\ux=\mtibas{x}{L}$ and positive integers $\un=\mtibas{n}{L-1}$ such
that for any $\ua\in\Aspucna$
the only flat connections in the link complement which satisfy
conditions\rx{1.57} are the $U(1)$-reducible ones.
\end{theorem}
We will prove this theorem in Appendix\rw{a1}.

Theorem\rw{red.conn} implies that $U(1)$-reducible connections $\Amu$
are the only contributors to $\ZubLLpua$ within the domain\rx{1.69}.
In other words, in that domain
\qq
\ZubLLpua = \sumuL \Zmba
\label{1.70}
\qqq
in the limit\rx{1.60}.
Now we want to pass from the limit\rx{1.60} to the
limit
\qq
K\rightarrow \infty,\qquad\ual,\ube = \const.
\label{1.71}
\qqq
This limit is similar to \ex{1.49}, so $\Zab$ receives the
contribution only from the trivial connection on $S^3$, and this
contribution is of the form
\qq
\Zab & = & \Ztrab
\label{1.72}\\
& = & \ZS\smmnp\dmmnL\, \ual^{2\um+1}\,\ube^{2\um\p+1}\,h^n
\nonumber
\qqq
(recall that in view of the first of \eex{1.55*} expansion in powers
of $h$ is equivalent to expansion in powers of $\Ki$).
A natural way to pass from \ex{1.70} to \ex{1.72} is to expand
$\Zmba$ in powers of $\ua$ while substituting $\ua = \ual/K$.
However, we have to stay within the domain\rx{1.69} while performing
this expansion. This can be done by expanding in powers of $\ua$ in
sequence: starting with $a_L$ and ending with $a_1$. Therefore
 we conclude
that
\qq
\Zab = \serc{ \sumuL \ZmbaK}{\ual}.
\label{1.73}
\qqq
Thus the second reason for making Conjecture\rw{conj1} is that
relation\rx{1.62} makes \ex{1.73} equivalent to \ex{1.45*1}.

\nsection{Calculating the RC invariant}
\label{s2}

The best way to derive $\JruaubLLp$ would be to study the asymptotic
behaviour of the colored Jones polynomial
$\JuaubLLp$ in the limit\rx{1.60} or, in other words, in
the limit
\qq
q\rightarrow 1,\qquad \qua,\ube = \const.
\label{2.1}
\qqq
Since we can not do this directly yet, we have to take a detour. We
will study the limit
\qq
q\rightarrow 1,\qquad \ual,\ube=\const.
\label{2.1*}
\qqq
This means that we will work with the Melvin-Morton expansion
\qq
\JuaubLLp = \snzi \Pbass{\ual,\ube}{n}{\cLLp}\,h^n,
\qquad \Pbass{\ual,\ube}{n}{\cLLp} \in \IQ[\ual,\ube]
\label{2.2}
\qqq
as it comes from the $\chR$-matrix representation described in
subsection\rw{Jpol}. However we will derive the expansion\rx{2.2} in
two stages. At the first stage we will expand all instances of $q$ in
$\chR$ and $\Qq$ except for the expressions of the form
\qq
q^{\a_i\a_j},\; \qaj,\quad 1\leq i,j\leq L,
\label{2.3}
\qqq
which we will leave intact. This way we will present $\JuaubLLp$ in
the form which is very close to \ex{1.45*}. At the second stage we
will expand the expressions\rx{2.3} thus coming to the Melvin-Morton
expansion\rx{2.2} in the form\rx{1.45*1}.

\subsection{Preliminary calculations}
\label{s2.1}
\subsubsection{Breaking a strand}
\label{s2.1.2}
We start with \ex{1.29}, according to which
\qq
\JuaubLLp = \qmfslLLp \Tr_{\VBLab} \Qq \BLLp,
\label{2.*1}
\qqq
where $\cBLLp$ is an $N$-strand
braid whose closure produces the link $\cLLp$.

Let us introduce a notation
\qq
\as = \left\{
\begin{array}{cl}
\a_1 & \mbox{if $\cL$ is non-empty}\\
\b_1 & \mbox{if $\cL$ is empty}
\end{array}
\right.
\qqq
The space $\VBLab$ is a tensor product
\qq
\VBLab = \Vas \otimes \VBLabo,
\qqq
where $\VBLabo$ is what remains of the tensor product $\VBLab$ after
we remove the first factor $\Vas$.
We introduce a modified version of the operator\rx{1.28}
\qq
\Qqo = I\otimes \brot{\qHt}{(N-1)}.
\label{2.*2}
\qqq
Obviously,
\qq
\Tr_{\VBLab} \Qq \BLLp =
\Tr_{\Vas} \qHt \Tr_{\VBLabo} \Qqo \BLLp.
\label{2.*3}
\qqq
On the other hand, since the operator $\BLLp$ commutes with the
generators of $\suq$ and since the module $\Vas$ is irreducible, then
the operator $\Tr_{\VBLabo}\Qqo \BLLp$ is proportional to $I$
\qq
\Tr_{\VBLabo}\Qq \BLLp = CI, \qquad C\in\IC.
\label{2.*4}
\qqq
Therefore since $\Tr_{\Vas} \qHt = \bas$, then
\qq
\Tr_{\VBLab} \Qq \BLLp = \bas\, C.
\label{2.*5}
\qqq

There are various ways to determine the constant $C$. For example,
since
\qq
\lrbar{\lopa \Tr_{\Vas} \lomH}{\loz}=1,
\label{2.*5*1}
\qqq
then
\qq
C = \lrbar{ \lopa \Tr_{\Vas} \lomH \Tr_{\VBLabo} \Qqo \BLLp}{\loz},
\label{2.*7}
\qqq
and combining \eex{1.29} and\rx{2.*5} we conclude that
\qq
\JuaubLLp = \qmfslLLp \,\bas
\lrbar{ \lopa
\Tr_{\VBLab}
(\lomH \otimes \Aot{I}{(N-1)})
\, \Qqo \BLLp}{\loz}.
\label{2.*8}
\label{2.10}
\qqq
%
%
%


\subsubsection{Twisting $\chR$-matrix}
\label{s2.1.1}
In subsection\rw{Jpol} we described a representation of the braid
group in which a positive elementary braid is mapped into the matrix
$\chR=PR$. It is better for our purposes to use a slightly
modified representation in which the image of a positive elementary
braid is defined as
\qq
\chR = 
\qhaH{-}{-} PR\, \qhaH{}{+}.
\label{2.4}
\qqq
Here $\gh$ is an operator which multiplies the elements of the module
$\Vgr$ by $\gr$
\qq
\gh\colon \Vgr\rightarrow \Vgr,\qquad \gh(\fm) = \gr \fm.
\label{2.5}
\qqq
As a result, the matrices $\chR$ and $\chRi$ map $\Vtgg{1}{2}$ into
$\Vtgg{2}{1}$ according to the formulas
\qq
\chR(\fmom) & = & \snzgi{1}{2} \chR\mrpmn\fmomn{+}{-},
\label{2.5*1}\\
\chRi(\fmom) & = & \snzgi{2}{1} (\chRi)\mrmpn \fmomn{-}{+},
\label{2.5*2}
\qqq
where
%
\qq
\chR\mrpmn & = & \qaoat{}
\spb{-m_2-l}{-\gr_2}{m_1-n+l}{l}\,
\label{2.6}\\
&&\qquad\times
q^{-m_2(\gr_1-1)+m_1 m_2 + \hnn},
\nonumber\\
(\chRi)\mrmpn & = &
\qaoat{-} \spb{m_1+l}{\gr_1}{-m_2+n-l}{-l}\,
\label{2.7}\\
&&\qquad\times
q^{m_1(\gr_2-1) + n(\gr_2-\gr_1) - m_1 m_2 - \hnn}.
\nonumber
\qqq
%

The reason for twisting $\chR$-matrix is that the matrix
elements\rx{2.6}
and\rx{2.7} have certain integrality properties.
\begin{lemma}
\label{lem.int}
There exist the polynomials $\Smni{\pm}{t_1,t_2}\in\ZZtotq$ such that
\qq
\chR\mrpmn & = & \qaoat{} \Smni{+}{\qao,\qat},
\\
(\chRi)\mrmpn & = & \qaoat{-} \Smni{-}{\qao,\qat}.
\qqq
\end{lemma}
\proof
Since
\qq
\plon\lrbc{ \omq{m_1-n+l}\over \omq{l}},\quad
\plon\lrbc{ \omq{-m_2+n-l} \over \omq{-l}},\quad q^{\pm\hlf n(n+1)}
\in\ZZqqi,
\qqq
then we can simply choose
\qq
\Smni{+}{t_1,t_2} & = &
\plon\lrbc{ (q^{-m_2-l} - t_2^{-1})\bomq{m_1-n+l}\over \omq{l}}
t_1^{-m_2}\,q^{m_2+m_1 m_2 + \hlf n(n+1)},
\\
\Smni{-}{t_1,t_2} & = &
\plon\lrbc{ (q^{m_1+l} - t_1) \bomq{-m_2+n-l} \over \omq{-l}}
t_1^{-n}\,t_2^{m_1+n}\,q^{-m_1 - m_1 m_2 - \hlf n(n+1)}.
\qqq
\qed

From now on we will use the matrices\rx{2.6},\rx{2.7} instead of $PR$
and its inverse in order to construct the operator $\BLLp$
which represents the braid $\cBLLp$.
However, this substitution does not affect \ex{2.*8}.
Indeed, if $\Vtaa{1}{2}$ is a part of the tensor
product $\VBLab$,
then a conjugation of the matrix $PR$ by the operator
$\qhaH{}{+}$ in \ex{2.4} is equivalent to the conjugation of the
whole space $\VBLab$ by
\qq
\bigotimes_{i=1}^N q^{ {H\over 8} \lrbc{
3\sjoimo (\gr(j)-1) + \sjioN (\gr(j)-1)}},
\label{2.8}
\qqq
where $\gr$ maps the number of the strand to the color of the
corresponding link component.
Thus we see that
the conjugation in \ex{2.4} has no effect on the trace of
\ex{2.*8}.
\subsubsection{Parametrized and regularized trace}

We want to introduce two more operators acting on $\VBLab$.
For a positive integer $N$ and for
two complex variables $\lzl=(\lo,\l)$ we introduce an operator
\qq
\QulN = \lo^{-H}\otimes \brot{ \l^{-H}}{(N-1)}
\label{2.12}
\qqq
and a parametrized trace
\qq
\JabLpl = \qmfslLLp \,\bas
\Tr_{\VBLab}\QulN\,\Qqo \BLLp.
\label{2.13}
\qqq
In contrast to the Jones polynomial, the parametrized
trace\rx{2.13} is not
an invariant of the link
$\cLLp$, because it depends on the choice of
the braid $\cBLLp$ whose closure is $\cLLp$.
Equation\rx{2.10} implies a relation between the parametrized trace
and the colored Jones polynomial
\qq
\lrbar{\lopa\JabLpl}{\l=1\atop\loz} = \JuaubLLp.
\label{2.14}
\qqq

Next we introduce a projection operator acting on $\Vgr$
\qq
\PM(\fm) =
\left\{
\begin{array}{cl}   \displaystyle
\fm
&\mbox{if $|\gr-1 -2m|\leq M$,}
\\
\displaystyle
0
&\mbox{if $|\gr-1-2m|> M$.}
\end{array}
\right.
\label{2.15}
\qqq
%
Let us replace the operators $\chR$ in the
expression for $\BLLp$
with `regularized' operators $\PMtb\,\chR\,\PMtb$
at $\cLpLpt$ crossings, with operators $\PMItb\,\chR\,\IPMtb$ at
$\cLLpt$ crossings and with operators $\IPMtb\,\chR\,\PMItb$ at
$\cLpLt$ crossings. We perform the same substitution for the
operators $\chRi$. We denote the resulting
regularized braid operator as $\BMLLp$.
We also introduce a projection operator
\qq
\PLLp = \Aot{\PM}{N}
\qqq
and we define a parametrized regularized trace as
\qq
\JabLplM =
\qmfslLLp \,\bas
\Tr_{\VBLab} \PLLp\,\QulN\,\Qqo \BMLLp.
\label{2.18}
\qqq
Since $\PM = I$ on $\Vgr$ if $M\geq \gr-1$, then
\qq
\JabLplM = \JabLpl\quad\mbox{if $\uM\geq\max\{\ual,\ube\}-1$.}
\label{2.19}
\qqq
%
\subsubsection{Verma module resolution}
Our next step is to resolve the traces over $\suq$ modules attached
to the strands of $\cL$ in terms of Verma modules. For $\gr\in\ZZ$, an
$\suq$ Verma module $\Vgri$ is an infinite-dimensional space with the
basis $\fm$, $m\geq 0$. The action of the $\suq$ generators on
$\Vgri$ is given by \eex{1.21}--(\ref{1.23}),
while the operators $\chR$ and $\chRi$ map
the tensor product $\Vtaai{1}{2}$ into $\Vtaai{2}{1}$ according to
the formulas which are similar to\rx{2.5*1},\rx{2.5*2}
\qq
\chR(\fmom) & = & \snzmi{1} \chR\mrpmn \fmomn{+}{-},
\label{2.19*1}\\
\chRi(\fmom) & = & \snzmi{2} (\chRi)\mrmpn \fmomn{-}{+},
\label{2.19*2}
\qqq
where the matrix elements are given by \eex{2.6} and\rx{2.7}.

For a positive integer $\gr$ there is a natural immersion
\qq
\Vgr \rightarrow \Vgri,\qquad \fm\mapsto \fm.
\label{2.27}
\qqq
Moreover, there is an isomorphism of $\suq$ modules
\qq
\Vgri/\Vgr = \Vmgri.
\label{2.28}
\qqq
If we choose the vectors $\fm$, $m\geq \gr$ as the basis of
$\Vgri/\Vgr$, then this isomorphism identifies
\qq
f_{m+n} \leftrightarrow { [\gr+m]!\over [m]!}\,\fm.
\label{2.29}
\qqq

We need another family of infinite-dimensional modules. For an
integer $\gr$,
$\Vgrip$ is an infinite-dimensional space with the basis $\fm$,
$m\in\ZZ$.
The action of the $\suq$ generators on
$\Vgri$ is given again by \eex{1.21}--(\ref{1.23}),
while the operators $\chR$ and $\chRi$ map
the tensor product $\Vtaaip{1}{2}$ into $\Vtaaip{2}{1}$ according to
the formulas which are similar to\rx{2.5*1},\rx{2.5*2}
\qq
\chR(\fmom) & = & \sngz \chR\mrpmn \fmomn{+}{-},
\label{2.19*3}\\
\chRi(\fmom) & = & \sngz (\chRi)\mrmpn \fmomn{-}{+},
\label{2.19*4}
\qqq
where the matrix elements are given again by \eex{2.6} and\rx{2.7}.
Also, for a positive $\gr$ there is an immersion
\qq
\Vgr \rightarrow \Vgrip,\qquad \fm\mapsto \fm.
\label{2.27*1}
\qqq

Consider two tensor products of $\suq$ modules
\qq
&\VBLai = \bigotimes_{j=1}^L \Aot{\Vaji}{N_j},
\qquad
\VBLpbip =
\bigotimes_{j=1}^{L\p} \VBLpjubip,
\label{2.21**1}
\qqq
where in the second tensor product
\qq
\VBLpjubip =
\Vbj\otimes \Aot{\Vbjip}{(N\p_j-1)}.
\label{2.21**2}
\qqq
%
Since the action of $\chR$-matrix was extended to spaces $\Vgri$ and
$\Vgrip$, then the action of the operator $\Qul\,\Qqo \BM$ can be
extended to $\Vti$ and $\Vtip$. For a set of integer (not
necessarily positive) colors $\ual$ and for a set of positive
integer colors $\ube$ we define the following trace
\qq
\JrabLplM = \qmfslLLp\,\bas \Tr_{\Vti}\PLLp\, \QulN\,\Qqo \BMLLp.
\label{2.22*}
\qqq
%

\begin{proposition}
\label{resol}
The trace $\JrabLplM$ has the following properties:
\begin{itemize}
\item[(1)]
If $\cL$ is empty, then
\qq
\Jrbasl{\ube}{\empt,\cLp}{\lzl;M} =
\Jbasl{\ube}{\cLp}{\lzl;M}.
\label{2.22**}
\qqq
\item[(2)]
For any $j$, $1\leq j\leq L\p$,
\qq
\lrbar{\JrabLplM}{\b_{j}=1}
=
\Jrbasl{\ual;\ube\xrem{j}}{\ccLLp\xrem{j}}{\lzl;M}.
\label{2.22**1**1}
\qqq

\item[(3)]
The space $\VBLpb$ can be replaced by $\VBLpbip$ in \ex{2.22*}
\qq
\hspace{-20pt}\JrabLplM =
\qmfslLLp\,\bas \Tr_{\Vtip}\PLLp \QulN\,\Qqo \BMLLp
\label{2.22*1}
\qqq

\item[(4)]
If $\cL$ is non-empty, then for any $j$, $1\leq j \leq L$
\qq
\JrLabLplMj =
\smumu \JrmLabLplMj.
\label{2.22*2}
\qqq
\end{itemize}
\end{proposition}

\pr{claims 1 and 2 of Proposition}{resol}
Equation\rx{2.22**1**1} is obvious if we compare \eex{2.22*}
and\rx{2.18}. To prove \ex{2.22**1*1} we observe that
$\dim V_{\bLp}=1$ if $\bLp=1$. As a result, the matrix elements of
$\chR$-matrices\rx{2.5*1},\rx{2.5*2} which involve at least one space
$V_{\bLp}$ are equal to 1. The same is true for the matrix elements
of the operators $\PLLp$, $\QulN$ and $\Qqo$ which act on
$V_{\bLp}$.
Finally, it
follows from \ex{1.30} that
\qq
\lrbar{\qmfslLLp}{\bLp=1} = q^{-\fsl(\cL\cup\cLp\rem{1})}
\label{2.22**1**2}
\qqq
and thus we come to \ex{2.22**1**1}.\qed


The proof of claims~3 and 4 of
Proposition\rw{resol} requires a simple lemma
\begin{lemma}
\label{lemresol}
Let $V$ be a finite-dimensional linear space with a subspace
$W\subset V$. For a fixed positive integer $N$ we denote
\qq
U_j = \Vot{(j-1)} \otimes W\otimes \Vot{(N-j)},\qquad 1\leq j\leq N.
\label{2.23}
\qqq
Suppose that a linear operator $O$ acting on the space $\Vot{N}$ has
the following property: there exists a cyclic permutation of $N$
elements $\sigma$ such that
\qq
O(U_j) \subset U_{\sigma(j)}\quad\mbox{for any $j$, $1\leq j\leq N$}.
\label{2.24}
\qqq
Then
\begin{itemize}
\item[(1)]
The action of $O$ is well-defined on the spaces $\Aot{W}{N}$,
$\Aot{(V/W)}{N}$ and
there is a simple trace formula
\qq
\Tr_{\Aot{W}{N}} O = \Tr_{\Vot{N}} O - \Tr_{\Aot{(V/W)}{N}} O.
\label{2.25}
\qqq
\item[(2)]
Let $j$ be a fixed number, $1\leq j\leq N$.
A partial trace of $O$ over $N-1$ spaces $V$ of the tensor product
$\Vot{N}$ except the first space $V$ defines an operator
$\Tr_{\Vot{(N-1)}} O$ acting on $V$. Then
$\Tr_{\Vot{(N-1)}}(W)\subset W$ and
\qq
\Tr_{\Aot{W}{N}} O = \Tr_{U_1} O .
\label{2.26}
\qqq
\end{itemize}
\end{lemma}
\proof
We leave the proof to the reader. \qed

\pr{claims~3 and~4 of Proposition}{resol}
Since the $R$-matrix, as well as the operators $\PM$, $\QulN$ and
$\Qqo$ can be expressed in terms of the generators of $\suq$, then
their action commutes with immersions\rx{2.27} and\rx{2.27*1}. Upon
being acted on by the $\suq$ generators of $R$-matrices, the elements
of $\suq$ modules are directed along the braid strands by
the permutation operators $P$ which are present at every
elementary crossing. The permutation generated by $B$ among
the strands of any particular link component is cyclic. Therefore it
is easy to see that the operator
$\QulN\,\Qqo \BMLLp$
satisfies the
conditions of Lemma\rw{lemresol} for any space $\Aot{\Vaji}{N_j}$
($V$ being $\Vaji$ and $W$ being $\Vaj$) and for any space
$\Aot{\Vbjip}{N\p_j}$ ($V$ being $\Vbjip$ and $W$ being $\Vbj$).
The infinite dimension of the spaces $\Vaji$ and
$\Vbjip$ is not important, since the presence of the projection
operators $\PM$ in the expression for $\BMLLp$
makes the calculation of the traces effectively
finite-dimensional.
Therefore \ex{2.22*1} follows from \ex{2.26}, and
\ex{2.22*2} follows from \ex{2.25}. \qed

Proposition\rw{resol} has the following
\begin{corollary}
\label{cor.resol}
The regularized trace\rx{2.18} can be presented as a linear
combination of the traces\rx{2.22*}
\qq
&\JabLplM = \sumuLumb \JrmabLplM,
\label{2.22}
\qqq
\end{corollary}
\proof
Equation\rx{2.22} is a result of applying \ex{2.22*2} consequently
to all the components of the link $\cL$. \qed

\subsection{Expansion of parametrized trace in powers of $h$}
Our next step is to perform a partial expansion of the trace
$\JrabLplM$ in powers of $h$. We want to expand all instances of
$q$ except for the expressions\rx{2.3}.

\subsubsection{General properties}

Let us first establish some properties of the trace $\JrabLplM$.
We introduce the modified operators
\qq
&\Qqom = I\otimes \brot{q^{-\hm}}{(N-1)},
\qquad
\QulomN = \lo^{2\hm}\otimes\brot{ \l^{2\hm} }{(N-1)}.
\label{2.12*1}
\qqq
%
%
We also introduce a modified operator $\BMLLpm$ which is constructed
in the same way as $\BMLLp$, except that we drop the factors
$\qaoat{\pm}$ from the matrix elements of $\chR$ and $\chRi$.

\begin{proposition}
\label{prop.int}
For a given braid $\cBLLp$ and
for fixed $M$ and $\ube$ there exists a unique polynomial
$$\AMq{\ut}\in\ZZutqpm$$ such that
\qq
&\JrabLplM = \bas\,\qfii{\ual,\ube}{\cLLp}\,
\Labul
\,\AMq{\qual},
\label{2.22**1}\\
&\AMq{\qual} = \Tr_{\Vtip}\PLLp\, \QulomN\,\Qqom \BMLLpm,
\label{2.22**1*1}
\qqq
where
\qq
\lefteqn{
\fii{\ual,\ube}{\cLLp}  =  -\fslLLp
}&&
\label{2.22**2}\\
&&
+\;
\qrtr\Bigg(
\soijL \lij(\a_i-1)(\a_j-1) + \!\!\soijLp \lpij (\b_i-1)(\b_j-1)
+2\soijLLp \lppij (\a_i-1)(\b_j-1)\Bigg)
\nonumber\\
&&\hspace{2.5in}
+ \hlf\Bigg(\sjoL N_j (\a_j-1) + \sjoLp N\p_j (\b_j-1)
+ (1-\as)
\Bigg)
\nonumber
\qqq
and
\qq
\Labul = \lo^{1-\as}\,
\l^{\sjoL N_j(1-\a_j) + \sjoL N\p_j(1-\b_j) + (\as-1)}.
\label{2.22**3}
\qqq
\end{proposition}
\proof
Let us review the structure of the matrix elements of the
operators $\QulN$,
$\Qqo$ and $\BMLLp$ in the basis $\fm$. For $\fm\in\Vgr$
\qq
\lmH (\fm) = \l^{1-\gr}\,\l^{2m}\,\fm,
\qquad
\qHt(\fm) = q^{(\gr-1)/2}\,q^{-m}\,\fm.
\label{2.22**4}
\qqq
The expression for $\BMLLp$ contains the matrices
$\chR$ and $\chRi$, whose matrix elements are described by
Lemma\rw{lem.int}. Thus, if we assemble the factors $\l^{1-\gr}$ of
the first of \eex{2.22**4} into $\Labul$ and if we assemble the
factors $q^{(\gr-1)/2}$ of the second of \eex{2.22**4} together with
the factors $\qaoat{\pm}$ of \eex{2.6},\rx{2.7} and with the
factor $\qmfslLLp$ of \ex{2.22*1}
into the factor
$\qfii{\ual,\ube}{\cLLp}$, then the rest of the
product of matrix elements of the operators $\chR$, $\chRi$, $\PM$,
$\QulN$ and $\Qqo$ is a polynomial of $\qual,\qmual$,
and $\sqi{q}$ with integer coefficients. Note that if we use either
\ex{2.22} or \ex{2.22*1} for $\JrabLplM$, then due to the presence
of projection operators $\PM$, the range of summation of the indices
of the basis vectors $f_{m_1}\otidot f_{m_N}$ used to calculate the
matrix elements is finite, and the boundaries of this range do
not depend on $\ual$. Thus, taking a sum of
the product of matrix elements over the indices of basis vectors we
come to \ex{2.22**1}. The uniqueness of the polynomial $\AMq{\ut}$
which satisfies \ex{2.22*1} for any $\ual$ is obvious.\qed

%
%

If we expand the polynomials $\AMq{\ut}$ in powers of $h=q-1$
then we come to the following

\begin{proposition}
\label{cor.expint}
For a braid $\cBLLp$
whose closure is $\cLLp$, there is a unique set of
polynomials
\qq
\AMb{\ut}\in\ZZulutpm
\label{2.22**5*a1}
\qqq
such that
%
for fixed $\ual,\ube$ and $M$ the following series
in powers of $h$ converges in the vicinity of $h=0$
\qq
\JrabLplM =
\bas\,\qfii{\ual,\ube}{\cLLp}\,\Labul \snzi\AMb{\qual}\,h^n.
\label{2.22**5}
\qqq
The polynomials $\AMb{\ut}$ satisfy the
following property: for any $k\geq 0$
%
%
\qq
\lefteqn{
\Foxi{\a_i}\snzk
\AMbi{\a_i,\ube}{\cL\xrem{i},\cL_i\cup\cLp}{\ut\xrem{i}}
}&&
\label{2.22**5**1}\\
&& = \smumu\, \Foxi{\mu\a_i}\snzk
\lrbar{
\AMbi{\ube}{\ccLLp}{\ut}
}{t_i = q^{\mu\a_i} }
+ O(h^{k+1})\quad\mbox{as $h\rightarrow 0$},
\nonumber
\qqq
where
%
\qq
\Foxi{\a_i} =
\lo^{-\delta_{1i}\a_i}\,\l^{(\delta_{1i}-N_i)\a_i}
q^{ {\a_i\over 2}
\lrbc{
\sjoLp \lpp_{ij} (\b_i-1) - \sjoL l_{ij} + N_i - \delta_{1i}
} }
\pjoLni t_j^{l_{ij}\a_i/2}
\label{2.22**5**4}
\qqq
and $\delta_{1i}$ is the Kronecker symbol.
\end{proposition}
\proof
The uniqueness of $\AMb{\ut}$ is obvious.
A set of polynomials satisfying \ex{2.22**5} can be obtained by
expanding the polynomials $\AMq{\ut}$ in powers of $h$
\qq
\AMq{\ut} = \snzi \AMb{\ut}\,h^n,
\quad
\AMb{\ut}\in\ZZulutpm,
\label{2.22**5*1}
\qqq

Now we prove \ex{2.22**5**1}.
It follows from
\eex{2.22*2},\rx{2.22**1} and from the uniqueness of the polynomials
$\AMq{\ut}$ that
%
%
\qq
\lefteqn{
\Foxi{\a_i}\,\AMqi{\a_i,\ube}{\cL\xrem{i},\cL_i\cup\cLp}{\ut\xrem{i}}
} &&
\label{2.22**5**2}\\
&&\hspace{2in} =
\lrbar{
\smumu\,\Foxi{\mu\a_i}\AMqi{\ube}{\ccLLp}{\ut}}{t_i=q^{\mu\a_i}}.
\nonumber
\qqq
If we expand the polynomials $A$ in
both sides of this equation in powers of $h$, then we
come to \ex{2.22**5**1}.\qed

It turns out that the polynomials $\AMb{\ut}$ have a limit as
$M\rightarrow\infty$.

\begin{proposition}
\label{prop.rat}
For a braid $\cB$ whose closure is $\cLLp$ there
exist the polynomials
\qq
&\Qzutule,\Qnubut \in \ZZutulz
\label{2.87}
\qqq
such that there is a vicinity of $\ut=1$, $\lzl=0$ in which for
fixed $\ube$ there exists a limit
\qq
\lMi \AMb{\ut} = {\Qnubut\over \Qnzutule}.
\label{2.22*3*2}
\qqq
The polynomials\rx{2.87} can be selected in such a way that
\qq
&\lrbar{\Qzutule}{\ulzcl}   =
\left\{
\begin{array}{cl}
\PLouti{-1}\,(\thtmh)\,\AFLut & \mbox{if $\cL$ is non-empty,} \\
1 & \mbox{if $\cL$ is empty,}
\end{array}
\right.
\label{2.88}
\qqq
where $\PLout$ is defined by \ex{1.17*3}.
Then
\qq
\lrbar{\Qzubut}{\ulzcl}
=
\pjoLp{1-\tpioL{-\b_j} \over 1 - \tpioL{-1}},
\label{2.88*}
\qqq
and the expansion of
$\Qnubut$ in powers of $\ut-1$ at $\lo=0$, $\l=1$ has a form
\qq
&\lrbar{\Qnubut}{\ulzcl} = \sumzi \Qmnubut\,(\ut-1)^\um,
\label{2.22*3*3}\\
&\Qmnubut \in \IQube,\quad \deg_m q_{\um;n}\leq 2(|\um|+n) + 1.
\label{2.22*3*4}
\qqq
\end{proposition}

We will prove this proposition in the remainder of this
section. However, before we proceed with the detailed calculations of
$\AMb{\ut}$ we want to use Proposition\rw{prop.rat} in order to
derive some properties of the polynomials\rx{2.87}.
%

\begin{proposition}
The polynomials\rx{2.87} satisfy the following properties
\begin{itemize}
\item[(1)]
If the link $\cL$ is empty, then
\qq
[\b_1]\,\qfii{\ube}{\cLp}\,\Lii{\ube}{\lzl}
\snzi
Q_{\ube;n}(\cLp;\lzl)\,h^n
=
\Jbasl{\ube}{\cLp}{\ul}\quad\mbox{in $\IQhhbas{\ube}$.}
\label{2.22*3*6}
\qqq
\item[(2)]
For any $j$, $1\leq j\leq L\p$
%
\qq
\lrbar{\Qnubut}{\b_{j}=1} = \Qili{\ube\xrem{j};n}{\cLLp\xrem{j}}{\ut}.
\label{2.22*3*7}
\qqq
\item[(3)]
For any $i$, $1\leq i\leq L$
and for fixed $\ut\xrem{i}$ and $\uul$ such that
$\lrbar{\Qzutule}{t_i=1} \neq 0$
%
%
\qq
\lefteqn{
\smumu\,\Foxi{\mu\a_i}
\snzk
\lrbar{\Qnubut\over \Qnzutule}{t_i=q^{\mu\a_i}}
} &&
\label{2.22*3*5}\\
&&
\hspace{1.5in}
=
\Foxi{\a_i}
\snzk{\Qili{\a_i,\ube;n}{\cL\xrem{i},\cL_i\cup\cLp}{\ut\xrem{i}} \over
\Qiii{2n+1}{}{\cL\xrem{i}}{\ut\xrem{i}} }
+O(h^{k+1})
\nonumber
\qqq
as $\hz$.
\end{itemize}
\end{proposition}
\proof
Equation\rx{2.22*3*6} follows from \eex{2.22**} and\rx{2.19}.
Equation\rx{2.22*3*7} follows from \ex{2.22**1**1}. In order to prove
\ex{2.22*3*5} we consider
the relation \ex{2.22**5**1} in the vicinity of $\ut=1$,
$\lzl= 0$ in which a uniform convergence\rx{2.22*3*2} holds. Since
the substitution $t_i=q^{\mu\a_i}$ leaves us within the same
vicinity as we set $\hz$, then the uniformity of
convergence allows us to apply the limit $\Mi$ to both
sides of \ex{2.22*3*5}. Thus we obtain\rx{2.22*3*5} in the vicinity
of $\ut\xrem{i}=1$, $\lzl=0$. Since both sides of \ex{2.22*3*5} are
rational (and hence analytic) functions, then the
relation\rx{2.22*3*5} should also hold at any point
$\ut\xrem{i}, \lzl$ such that
$\lrbar{\Qzutule}{t_i=1}\neq 0$.\qed

\subsubsection{Fixing transitions}
We will prove Proposition\rw{prop.rat} by calculating the limit
$\lMi\AMb{\ut}$ explicitly. However, before we do this,
we have to introduce some notations. Let $V_1$
and $V_2$ be two $\suq$ modules (either finite, or
infinite-dimensional) with colors $\gr_1$ and $\gr_2$ and with
bases $\fm$, $m\in\cS_1$ and $\fm$, $m\in\cS_2$,
$\cS_{1,2}\subset \ZZ$. For $n\geq 0$ we define two operators
$\Rnpm\colon  V_1\otimes V_2\rightarrow V_2\otimes V_1$ such that
\qq
\Rnp(\fmom) & = &
\left\{
\begin{array}{cl}
\chR\mrpmn \fmomn{+}{-} &\mbox{if $m_2+n\in\cS_1$, $m_1-n\in\cS_2$,}\\
0 &\mbox{otherwise},
\end{array}
\right.
\label{2.41*1}\\
\Rnm(\fmom) & = &
\left\{
\begin{array}{cl}
(\chRi)\mrpmn \fmomn{-}{+}
&\mbox{if $m_2-n\in\cS_1$, $m_1+n\in\cS_2$,}\\
0 &\mbox{otherwise}.
\end{array}
\right.
\label{2.41*2}
\qqq
Obviously,
\qq
\chR = \sngz \Rnp,\qquad \chRi = \sngz \Rnm.
\label{2.41*3}
\qqq

Let us call the positive
$\cLLpt$ crossings and the negative $\cLpLt$ crossings \emph{the
losing crossings} and let us call the negative $\cLLpt$ and positive
$\cLpLt$ crossings \emph{the gaining crossings}. Suppose that the
braid $\cBLLp$
has $\Nl$ losing crossings, $\Ng$ gaining crossings and
$\Npp$ crossings of $\cLpLpt$ type and that the crossings within each
group are indexed in some order. Then three lists of non-negative
numbers $\unl=\mtibas{\nl}{\Nl}$, $\ung=\mtibas{\ng}{\Ng}$ and
$\unpp=\mtibas{\npp}{\Npp}$ assign a non-negative number to any
crossing other than those of $\cLLt$ type. After this assignment we
construct an operator
$\BMLLpun$, $\uun = \unl;\ung;\unpp$ by performing
the following transformation or the expression for $\BMLLp$:
we replace the operators $\chR$ and $\chRi$ at all the crossings
except of $\cLLt$ type by the operators
$\Rnp$ and $\Rnm$ whose numbers $n$ are equal to the numbers assigned
by $\uun$ to the corresponding crossings.


If we replace the braiding
operator $\BMLLp$ by $\BMLLpun$ in \ex{2.22*1}, then we get the trace
\qq
\lefteqn{
\JrabunLplM}
\label{2.43*4}\\
&& =
\bas\,\qfii{\ual,\ube}{\cLLp}\,\Labul
 \Tr_{\Vtip}\PLLp \QulomN\,\Qqom \BMLLpun
\nonumber\\
\qqq
Obviously, the trace\rx{2.43*4} shares the properties of
$\JrabLplM$ listed in Proposition\rw{prop.int} and
Corollary\rw{cor.expint}:
\begin{proposition}
For fixed $\ube$ and $M$ there exists a unique set of polynomials
\qq
&\AMqun{\ut}\in\ZZutqpm,\qquad
\AMbun{\ut}\in\ZZulutpm,
\\
& \AMqun{\ut} = \snzi \AMbun{\ut}\,h^n,
\qqq
such that
\qq
&\JrabunLplM  =
\bas\,\qfii{\ual,\ube}{\cLLp}\,\Labul
\,\AMqun{\qual},
\label{2.43*5}\\
&\AMqun{\qual} = \Tr_{\Vtip}\PLLp\, \QulomN\,\Qqom\,\BMLLpun,
\label{2.43*5*1}
\qqq
and for fixed $\ual$, $\ube$
\qq
\JrabunLplM  = \bas\,\qfii{\ual,\ube}{\cLLp}
\,\Labul \snzi\AMbun{\qual}\,h^n.
\label{2.43*6}
\qqq
\end{proposition}
\proof
The proof is absolutely similar to that of Proposition\rw{prop.int}
and Corollary\rw{cor.expint}.\qed

There is a simple relation between $\JrabunLplM$ and
$\JrabLplM$.
\begin{proposition}
\label{prop.easysum}
\qq
&\JrabunLplM = \suntM \JrabunLplM,
\label{2.43*7}\\
&\AMq{\ut} = \suntM\AMqun{\ut},
\label{2.43*8*1}\\
&\AMb{\ut} = \suntM\AMbun{\ut}.
\label{2.43*8}
\qqq
\end{proposition}
\proof
It is easy to see that due to the presence of projection
operators $\PM$ in the expression for $\BMLLpun$, its matrix elements
are zero unless $\uun\leq 2M$. Therefore in view of \eex{2.41*3}
\qq
\JrabunLplM = \suntMlim \JrabunLplM.
\label{2.43*7*1}
\qqq
Let $\HLp$ be a sum of $\suq$ generators $H$ acting on
the components of $\VBLpbip$ within the
tensor product $\Vtip$
The space $\Vtip$ is a sum
of eigenspaces $\Vtipi{n}$
of $\cLp$ corresponding to its eigenvalues
$n$
\qq
\Vtip = \bigoplus_{n\in\ZZ} \Vtipi{n}
\qqq
It is easy to see that
\qq
\lefteqn{
\PLLp\, \QulomN\,\Qqom \BMLLpun\,\colon
}\\
&&
\hspace{1.7in}
 \Vtipi{n} \rightarrow
\Vtipi{n+2(|\ung|-|\unl|)}.
\nonumber
\qqq
Therefore
\qq
\Tr_{\Vtip} \PLLp\,\QulomN\,\Qqom\BMLLpun = 0,\qquad
\mbox{unless $\absol{\ung}=\absol{\unl}$.}
\qqq
Then \ex{2.43*7} follows
from \ex{2.43*7*1}.
Relations\rx{2.43*8} follow from \ex{2.43*7} because of the
uniqueness of polynomials $A$.\qed

\subsubsection{Expansion of matrix elements}
\label{s2.2.3}
Now it's time to expand the matrix elements of
the operators $\chR$ and $\chRi$ participating in the expression for
$\BMLLpun$, as well as the matrix elements of the operator $\Qqom$.
Since we
have to study the expansion of the trace $\JrabLplM$ in powers of $h$
while keeping the expressions\rx{2.3} unexpanded then we will
need different formulas for expanding the matrix elements\rx{2.6}
and\rx{2.7} depending on the type of the crossing that generated them.
\begin{theorem}
\label{expand}
There exist the polynomials
\qq
&\Tjkpm
\momtn \in \IQmomtn,\quad j,k\geq 0,
\nonumber\\
&\deg_{m_1}\Tjkp, \deg_{m_2}\Tjkm\leq k,\qquad
\deg_{m_2}\Tjkp, \deg_{m_1}\Tjkm\leq j+k,
\nonumber\\
&\deg_{m_1}\Tjkpm + \deg_{m_2}\Tjkpm + \deg_{n} \Tjkpm \leq j + 2k,
\nonumber\\
&\Tzzpm \momtn =1,
\label{2.30}
\qqq
such that for fixed $\gr_1$, $\gr_2$, $m_1$, $m_2$,and $n$
\qq
\chR\mrpmn  & = &
\qaoat{}
{m_1\choose n}\, \bqmo{-\gr_2}^n\, q^{-m_2\gr_1}
\label{2.31}\\
&&
\hspace{1.5in}
\times
\sjzn \lrbc{ {\plzjmo (n-l)\over
\bqmo{-\gr_2}^j } \skzi \Tjkp\momtn\,h^{j+k} },
\nonumber\\
(\chRi)\mrmpn & = &
\qaoat{-}
{m_2\choose n}\bqmo{\gr_1}^n q^{m_1\gr_2 + n(\gr_2-\gr_1)}
\label{2.32}\\
&&
\hspace{1.5in}
\times
\sjzn \lrbc{ {\plzjmo (n-l)\over
\bqmo{\gr_1}^j }\skzi\Tjkm\momtn\,h^{j+k} }.
\nonumber
\qqq
The polynomials $\Tjkp$ and $\Tjkm$ are related
\qq
\skzi \Tjkm\momtn\,h^{j+k}  =
\skzi \Tjkp (m_2,m_1,n)\lrbc{ -{h\over 1+h}}^{j+k}.
\label{2.33}
\qqq
\end{theorem}
This theorem essentially says that there are at most two powers of
$m_1$, $m_2$ or $n$ for each power of $h$ in the expansion of the
matrix elements. This is easy to see by inspecting \eex{2.6}
and\rx{2.7}. The only non-trivial part of the theorem is that $j$
powers of $\bqmo{\pm\gr_{1,2}}$ in denominators of
\eex{2.31},\rx{2.32} are accompanied by the factor $\plzjmo(n-l)$.

\pr{Theorem}{expand}
%
This theorem is a simple corollary of Lemma~3.1
of\cx{Ro9}.\qed

We present explicit formulas for some polynomials $\Tjkpm$ in
Appendix\rw{a2}.

\begin{corollary}
\label{cor.expand}
There are three other ways to expand the matrix elements
$\chR\mrpmn$ and $(\chRi)\mrmpn$.
\begin{itemize}
\item[(1)]
There exist the polynomials
\qq
&\Tnkopm (\momtgz) \in \IQ[\momtgz],
\nonumber\\
&\deg_{m_1}
\Tnkopm \leq k,\quad \deg_{\gr}\Tnkopm + \deg_{m_2}\Tnkopm
\leq n+k,
\label{2.33*}\\
&\Tzzpm{1}(\momtgz) = 1,
\nonumber
\qqq
such that
\qq
\chR\mrpmn & = & \qaoat{}\, {m_1\choose n}\, q^{-m_2\gr_1}
\skzi \Tnkop (m_1,m_2,\gr_2)\, h^{n+k},
\label{2.34}\\
(\chRi)\mrmpn & = & \qaoat{-}\,{m_2\choose n}\,
q^{(m_1+n)\gr_2}
\skzi \Tnkom(m_1,m_2,\gr_1)\, h^{n+k}.
\label{2.35}
\qqq
\item[(2)]
There exist the polynomials
\qq
&\Tnktpm(\momgt) \in \IQ[\momgti],
\nonumber\\
&\deg_{m_2}\Tnktp,\deg_{m_1}\Tnktm\leq k,
\label{2.36}\\
&\Tzzpm{2}(\momgt) = 1,
\nonumber
\qqq
%
such that there is a bound on the powers of $m$ and $\gr$ in the
coefficients of the expansion of the polynomials in powers of $t-1$
\qq
&\Tnktpm(\momgt) = \slzi \Tnkl{\pm}(m_1,m_2,\gr)\,(t-1)^l,
\nonumber\\
&\Tnkl{\pm}(m_1,m_2,\gr)\in\IQ[m_1,m_2,\gr],
\nonumber\\
&\deg_{m_1} \Tnkl{+} + \deg_{\gr} \Tnkl{+} \leq k + l,
\label{2.36*}\\
&\deg_{m_2} \Tnkl{-} + \deg_{\gr} \Tnkl{-} \leq k + l,
\label{2.36*1}
\qqq
%
%
%
and for fixed $\gr_1$, $\gr_2$, $m_1$, $m_2$,and $n$
\qq
\chR\mrpmn & = & \qaoat{}\, \,
\skzi \Tnktp (m_1,m_2,\gr_1, q^{\gr_2} )\, h^{n+k},
\label{2.37}\\
(\chRi)\mrmpn & = & \qaoat{-}\,
\skzi \Tnktm(m_1,m_2,\gr_2,q^{\gr_1})\, h^{n+k}.
\label{2.38}
\qqq
\item[(3)]
There exist polynomials
\qq
&\Tnkthpm(\momtgzth) \in \IQ[\momtgzth],
\nonumber\\
&\deg_{m_1}\Tnktpm + \deg_{\gr_1}\Tnktpm \leq k,\qquad
\deg_{m_2}\Tnktpm + \deg_{\gr_2}\Tnktpm \leq k,
\label{2.38*}\\
&\Tzzpm{3}(\momtgzth) = 1,
\nonumber
\qqq
such that
for fixed $\gr_1$, $\gr_2$, $m_1$, $m_2$,and $n$
\qq
\chR\mrpmn & = & \qaoat{}\, \,
\skzi \Tnkthp(\momtgzth) \, h^{n+k},
\label{2.39}\\
(\chRi)\mrmpn & = & \qaoat{-}\,
\skzi \Tnkthm(\momtgzth)\, h^{n+k}.
\label{2.40}
\qqq
\end{itemize}
\end{corollary}
\proof
One simply has to expand \eex{2.31} and\rx{2.32} in the corresponding
parameters. We leave the details to the reader.\qed

Finally, we expand the matrix elements of the operators $q^{\hm}$
which
comprise the operator $\Qqom$. The matrix elements of
$q^{\hm}$ acting on
$\Vgr$ are all diagonal and
\qq
(q^{\hm})_m^m = q^m =
\skzi {m\choose k} h^k.
\label{2.41}
\qqq
(note that ${m \choose k}=0$ for $k>m$ if we treat ${m\choose k}$ as
a polynomial in $m$ of degree $k$).

\begin{proposition}
\label{prop.subst}
Suppose that we make the following substitution in the expression for
the regularized braiding matrix $\BMLLpun$: we use
expressions\rx{2.31},\rx{2.32} for the matrix elements of $\chR$ and
$\chRi$ at $\cLLt$ crossings, expressions\rx{2.34},\rx{2.35} at
losing crossings, expressions\rx{2.37},\rx{2.38} at gaining
crossings and expressions\rx{2.39},\rx{2.40} at $\cLpLpt$ crossings.
Suppose that we also use the expression\rx{2.41} for the matrix
elements of the operator $\Qqo$. Then after this substitution the
trace\rx{2.43*4} transforms into the series\rx{2.43*6}.
\end{proposition}
\proof
The presence of the projection operator $\PLLp$ in the
trace\rx{2.43*4} guarantees that only a finite number of matrix
elements participating in that trace is non-zero. Therefore the
non-zero products of expanded matrix elements yield a series of the
form\rx{2.43*6} with $\AMbun{\ut}\in\IQulutpm$, and this series
converges uniformly in a vicinity of $h=0$. Now the claim of the
proposition follows from the uniqueness of the expansion\rx{2.43*6}.
\qed
\begin{corollary}
The bounds in the sum of \ex{2.43*8} can be made $M$-independent
\qq
\AMb{\ut} = \suneqn\AMbun{\ut}.
\label{2.41**1}
\qqq
\end{corollary}
\proof
Expansions\rx{2.34},\rx{2.35} of the matrix elements which
appear at losing crossings and expansions \rx{2.39},\rx{2.40} of the
matrix elements which appear at the $\cLpLpt$ crossings
all start at $h^n$. Since, according to
Proposition\rx{prop.subst}, these expansions lead to \ex{2.43*6},
then we conclude that
\qq
\AMqun{\ut} = \sngn \AMbun{\ut}\,h^n.
\label{2.41**2}
\qqq
If we substitute this expansion into \ex{2.43*8*1} and compare the
result with \ex{2.22**5*1}, then we come to \ex{2.41**1}.\qed

\subsubsection{Parametrized matrices}
\label{s2.2.4}
Now we are going to use a special trick in order to calculate the
trace of the expanded matrices. This trick is a usual method of
calculating perturbative expansions in quantum field theory. Its
application is mathematically rigorous in the finite-dimensional
calculations of this paper. The idea is that in order to calculate a
trace of a complicated matrix one may introduce a simpler matrix
which depends on some artificial parameters and which satisfies two
properties: its trace is easy to calculate and the complicated matrix
can be expressed as a combination of derivatives of the parametrized
matrix over its parameters.

The first step in implementing our strategy is to construct a
matrix $\BMLLpunps$ depending on a list of parameters
$\ut;\uue$, together with a differential operator $\Dsue$ which is a
power series in $h$ whose coefficients are polynomials in the
derivatives $\del_\uue$. We will construct $\BMLLpunps$ in such a way
that the matrix elements of
$\lrbar{\Dsue\BMLLpunps}{\ut=\qua\atop\uue=0}$
are equal to those of the
matrix
$\BMLLpun$ expanded as described in Proposition\rw{prop.subst}.
In fact, we will have to consider an extra $\suq$ module $\Vzi$, so
that our matrix $\BMLLpunps$ will act on a tensor product
$\Vzi\otimes\Vtip$
and the expansion of $\BMLLpun$ will be reproduced
by an operator which is a partial trace
$\lrbar{\Dsue\Tr_{\Vzi}\lst^{-H}\BMLLpunps}{\ut=\qua\atop\lst,\uue=0}$.

The basis of our construction is the matrix $\BMLLpun$ itself. We will
present every $\chR$, $\chRi$, $\Rnpm$ participating in the
expression for $\BMLLpun$
as $D(\del_\ue) R[\ue]$, where $D(\del_\ue)$
is a power series in $h$ with coefficients which are polynomials in
the derivatives over some parameters $\ue$, and $R[\ue]$ are
operators depending on parameters $\ue$. These expressions will match
the expansion formulas of Theorem\rw{expand} and
Corollary\rw{cor.expand}. Then we will replace the original
matrices of $\BMLLpun$ by $R[\ue]$ thus producing $\BMLLpunps$, while
the operator $\Dsue$ will be the product of individual operators
$D(\del_\ue)$.

We begin by defining a parametrized matrix $\Rpm\tteee$,
$\ue = (e_1,e_2,e_3,e_4)$ which acts on
a tensor product  of two Verma modules $\Vgriia{1}\otimes\Vgriia{2}$:
\qq
\lefteqn{
\Rp\tteee(\fmom)
}
\label{2.42}\\
&
\hspace{2in}
= & \snzmi{1} (\Rp\tteee)\mrpmn \fmomn{+}{-},
\nonumber\\
\lefteqn{
\Rm\tteee(\fmom)
}
\label{2.43}\\
&
\hspace{2in}
= & \snzmi{2} (\Rm\tteee)\mrmpn \fmomn{-}{+},
\nonumber
\qqq
where
\qq
\hspace{-30pt}
(\Rp\tteee)\mrpmn & = & {m_1\choose n}(1-t_2^{-1}+\e_4)^n\,
t_1^{-m_2}\, e^{m_1\e_1 + m_2\e_2 + n \e_3},
\label{2.44}\\
\hspace{-30pt}
(\Rm\tteee)\mrmpn & = & {m_2\choose n}(1-t_1+\e_4)^n\, t_1^{-n}\,
t_2^{m_1+n}\, e^{m_1\e_1 + m_2\e_2 + n\e_3}.
\label{2.45}
\qqq
The matrix elements\rx{2.44},\rx{2.45} are related to those
of\rx{2.6},\rx{2.7}, because according to \eex{2.31},\rx{2.32}
\qq
\lefteqn{
\chR\mrpmn   =
\qaoati{}{\a}{\a}
}
\label{2.46}\\
&&\times
\lrbar{
\lrbc{  \sjkzi \Tjkp\deee\,\del^j_{\e_4}\,h^{j+k} }
(\Rp\tteee)\mrpmn}{\clsin},
\nonumber\\
\lefteqn{
(\chRi)\mrpmn  =
\qaoati{-}{\a}{\a}
}
\label{2.47}\\
&&\times
\lrbar{
\lrbc{ \sjkzi \Tjkm\deee\,\del^j_{\e_4} \,h^{j+k} }
(\Rm\tteee)\mrmpn}{\clsin}.
\nonumber
\qqq
Note that we could replace the sums $\sjzn$ of \eex{2.31},\rx{2.32}
with the sums $\sjzi$, because the terms with $j>n$ in
\eex{2.46},\rx{2.47} are zero due to the fact that if
$j>n$, then $\del^j_{\e_4} (1-t_{1,2}^{\mp 1} +\e_4)^n=0$. Following
our strategy, we replace the operators $\chR$, $\chRi$ which appear
in the expression for $\BMLLpun$
at the $\cLLt$ crossings, with operators
$\Rpm\tteee$, while the product of operators
$\sjzi \lrbc{ \del^j_{\e_4} \skzi \Tjkp\deee\,h^{j+k} }$
coming from all $\cLLt$ crossings
is the first
factor of $\Dsue$.

Matrix elements\rx{2.34},\rx{2.35} and\rx{2.39},\rx{2.40}
require a more elaborate
treatment. We define the operators
$\Rgl\eegl$ which act on $\Vzi\otimes\Vai$ according to the formula
\qq
\Rg\eeg(\fmomz) & = & \snzmi{0} (\Rg\eeg)
\msmpnz \fmnmoz{-}{+},
\label{2.48}\\
\Rl\eel(\fmomz) & = & \snzmi{1} (\Rl\eel)
\mspmnz \fmnmoz{+}{-},
\label{2.49}
\qqq
where
\qq
(\Rg\eeg)\msmpnz & = & {m_0\choose n} e^{\eg_1 m_1}\,(\eg_2)^n,
\label{2.50}\\
(\Rl\eel)\mspmnz & = & {m_1\choose n} e^{\el_1 m_1}\,(\el_2)^n.
\label{2.51}
\qqq
We also define another pair of operators depending on parameters
$t_1, t_2$
which act on $\Vgrip$
\qq
\Onpm[t_1,t_2](\fm) = t_1^m t_2^n\,f_{m\pm n},
\label{2.51*}
\qqq
and an operator $\Chsn$ which acts on the space $\Vzi$ as
\qq
\Chsn( \fm ) = {\fm\over {m\choose n}}.
\label{2.51*1}
\qqq

Now we have to consider four types of mixed crossings.
\begin{enumerate}
\item
Losing crossings.
\begin{enumerate}
\item
Positive $\cLLpt$ crossing. Consider an operator
\qq
&\Ri{-,+}
[\eele{1}]\,\colon\,
\Vzi\otimes\Vaii{1}\otimes\Vbiip{2} \rightarrow
\Vzi\otimes\Vbiip{2}\otimes\Vaii{1},
\label{2.52}\\
&\Ri{-,+}[\eele{1}] =
(I\otimes P)\,(\Rl\eel\otimes\Onp[e^{\e\p}t_1^{-1},1]).
\label{2.53}
\qqq
It is easy to see that according to \ex{2.34},
\qq
\lefteqn{
\dlti{+}\,
(\Rnm)\motp = \qaoati{}{\a}{\b}
}\hspace{0.4in}&&
\label{2.54}
\\
&&\hspace{-10.5pt}\times\lrbc{ {1\over n!}\,\del_{\el_2}^n\skzi
\Tnkop (\del_{\el_1},\del_{\e\p},\b_2)\,h^{n+k}
}
\lrbar{
(\Ri{-,+}[\eele{1}])\motpe
}{t_1=\qai{1}\atop\el_1,\el_2,\e\p=0},
\nonumber
\qqq
where $\dlti{+}$ is the Kronecker symbol.
%
\item
Negative $\cLpLt$ crossing.
Consider an operator
\qq
&\Ri{-,-}
[\eele{2}]\,\colon\,
\Vzi\otimes\Vbiip{1}\otimes\Vaii{2} \rightarrow
\Vzi\otimes\Vaii{2}\otimes\Vbiip{1},
\label{2.55}\\
&\Ri{-,-}[\eele{2}] = (\Rl\eel\otimes\Onp[e^{\e\p}t_2,t_2])
\,(I\otimes P).
\label{2.55*}
\nonumber
\qqq
According to \ex{2.34},
\qq
\lefteqn{\dlti{+}
(\Rnm)\motp = \qaoati{-}{\b}{\a}
}\hspace{0.4in}&&
\label{2.56}\\
&&\hspace{-10.5pt}\times\lrbc{ {1\over n!}\,\del_{\el_2}^n\skzi
\Tnkop (\del_{\e\p},\del_{\el_1},\b_1)\,h^{n+k}
}
\lrbar{
(\Ri{-,-}[\eele{2}])\motpe
}{t_2=\qai{2}\atop\el_1,\el_2,\e\p=0}.
\nonumber
\qqq
\end{enumerate}
\item Gaining crossings.
\begin{enumerate}
\item
Positive $\cLpLt$ crossing. Consider an operator
\qq
&\Ri{+,+}
[\eege{2}]\,\colon\,
\Vzi\otimes\Vbiip{1}\otimes\Vaii{2} \rightarrow
\Vzi\otimes\Vaii{2}\otimes\Vbiip{1},
\label{2.57}\\
&\Ri{+,+}[\eege{2}] = (\Rg\eeg\otimes\Onm[e^{\e\p},1])
\,(\Chsn\otimes P).
\label{2.58}
\qqq
According to \ex{2.37},
\qq
\lefteqn{
\dlti{-}\,
(\Rnp)\motp = \qaoat{}{\b}{\a}
}\hspace{0.4in}&&
\label{2.59}
\\
&&\hspace{-62pt}\times\lrbc{ {1\over n!}\,\del_{\eg_2}^n\skzi
\Tnktp (\del_{\e\p},\del_{\eg_1},\b_1,t_2)\,h^{n+k}
}
\lrbar{
(\Ri{+,+}[\eege{2}])\motpe
}{t_2=\qai{2}\atop\eg_1,\eg_2,\e\p=0 }.
\nonumber
\qqq
%
\item
Negative $\cLLpt$ crossing. Consider an operator
\qq
&\Ri{+,-}
[\eege{1}]\,\colon\,
\Vzi\otimes\Vaii{1}\otimes\Vbiip{2} \rightarrow
\Vzi\otimes\Vbiip{2}\otimes\Vaii{1},
\label{2.60}\\
&\Ri{+,-}[\eege{1}] = (I\otimes P)\,
(\Rg\eeg\otimes\Onm[e^{\e\p},1])\,
(\Chsn\otimes \Aot{I}{2}).
\label{2.61}
\qqq
According to \ex{2.38},
\qq
\lefteqn{
\dlti{-}\,
(\Rnp)\motp = \qaoat{-}{\a}{\b}
}\hspace{0.4in}&&
\label{2.63}
\\
&&\hspace{-65pt}\times\lrbc{ {1\over n!}\,\del_{\eg_2}^n\skzi
\Tnktp (\del_{\eg_1},\del_{\e\p},\b_2,t_1)\,h^{n+k}
}
\lrbar{
(\Ri{+,+}[\eege{1}])\motpe
}{t_1=\qai{1}\atop\eg_1,\eg_2,\e\p=0 }.
\nonumber
\qqq
\end{enumerate}
\end{enumerate}

Finally, we turn to $\cLpLpt$ crossings. We introduce two operators
\qq
&\Rppn{\pm}[\e\p_1,\e\p_2]\colon
\Vbiip{1}\otimes\Vbiip{2}\rightarrow \Vbiip{2}\otimes\Vbiip{1},
\label{2.63*}\\
&\Rppn{\pm}[\e\p_1,\e\p_2] = P\,(\Onmp[e^{\e\p_1},1]\otimes
\Onpm[e^{\e\p_2},1]),
\label{2.64}
\qqq
so that according to \eex{2.39},\rx{2.40},
\qq
(\Rnpm)\motp = \qaoati{\pm}{\b}{\b}\,
\lrbc{
\Tnkthpm
(\b_1,\b_2,\del_{\e\p_1},\del_{\e\p_2})\,h^{n+k}
}
\lrbar{
(\Rppn{\pm}[\e\p_1,\e\p_2]) \motp
}{\e\p_1,\e\p_2=0}.
\label{2.65}
\qqq

Now, according to our strategy, we take the expression for
$\BMLLpun$ and
replace $\chR$ and $\chRi$ matrices at $\cLLt$ crossings with the
matrices\rx{2.42},\rx{2.43}, replace $\Rnpm$ matrices at losing
crossings with matrices\rx{2.53},\rx{2.55*}, replace $\Rnpm$
matrices at gaining crossings with matrices\rx{2.58},\rx{2.61} and
replace $\Rnpm$ matrices at $\cLpLpt$ crossings with
matrices\rx{2.64}. The result is the parametrized matrix $\BCMLLpunps$
acting on $\Vzi\otimes\Vtip$,
where $\uue = (\ue;\ue\p)$, $\ue$ and
$\ue\p$ being all the parameters $\e$ and $\e\p$ from the individual
parametrized matrices. We assume that the matrices
\rx{2.53},\rx{2.55*},\rx{2.58} and\rx{2.61} act on the space $\Vzi$
in the following order: losing crossing matrices precede gaining
crossing matrices and among themselves those matrices are ordered
according to the indexing of losing and gaining crossings.

The operator $\BCMLLpunps$ is a tensor product of
operators
\qq
&\BCMLLpunps = \BCMLunv\otimes\bigotimes_{j=1}^{L\p} \BMLjpunv,
\label{2.65*}\\
&\BCMLunv\colon \VBLaiz\rightarrow\VBLaiz,
\qquad
\BMLjpunv\colon \VBLpjubip\rightarrow\VBLpjubip
\nonumber
\qqq
(note that $\BCMLunv$ does not depend on $M$).
To construct the operator $\BCMLunv$, we place the operators
$\Rp$ at positive elementary crossings of $\cL$ and the
operators $\Rm$ at negative elementary crossings of $\cL$. We
also place the operators $\Rl$ at losing crossings and
$\Rg\,(I\otimes\Chsn)$ at gaining crossings. To construct the
operators $\BMLjpunv$ we simply place the appropriate operators
$\Onpm$ at every place where a strand of $\cLj$ crosses some other
strand.

We construct the differential operator $\Dsue$ simultaneously with
$\BMLjpunv$
by taking the products of differential operators which
accompany individual parametrized matrices.
\begin{proposition}
\label{prop.oper}
The operator $\Dsue$ has the following form:
\qq
&\Dsue = \snzi \Dsnued\,h^n,
\quad
\Dsnue\in \IQ[\sqi{\ut},\ube,\uue],
\label{2.66}\\
& D_{\uu{0};0}(\ut;\ube;\uue) = 1.
\label{2.66*}
\qqq
The expansion of the polynomials $\Dsnue$ in powers of $\ut-1$ has
the form
\qq
&\Dsnue = \sumzi \Dsnmue(\ube;\uue)\,(\ut-1)^{\um},
\quad
\Dsnmue(\ube;\uue)\in\IQ[\ube,\uue],
\label{2.67}\\
& \deg_{\ube}\Dsnmue + \deg_{\uep}\Dsnmue \leq |\um| + n.
\nonumber
\qqq
%
%
\end{proposition}
\proof
Equation\rx{2.66*} follows from the similar relations in
\eex{2.30},\rx{2.33*},\rx{2.36} and\rx{2.38*}.
Equations\rx{2.66} and\rx{2.67} follow easily from Theorem\rw{expand}
and Corollary\rw{cor.expand}.\qed

We chose the matrix elements of gaining and losing crossing operators
$\Ri{\pm}$ in such a way that if we take any diagonal element
of the operator
$\Dsue\BMLjpunv$ at
$|\ung|=|\unl|$ in the space $\Vzi$,
then, according to Proposition\rw{prop.subst}, the resulting
operator acting on
$\Vtip$
coincides with the expansion of
$\BMLLpun$,
as described in that proposition. In particular, if we take a
diagonal element corresponding to $f_0$, then
\qq
\BMLLpun = \lrbar{\Dsue \Tr_{\Vzi} (
\PM\lstar^{2\hm}\otimes \Aot{I}{N})
\,\BCMLLpunps}{{\lstar=0\atop \ut=\qua}\atop \uue = 0}
\label{2.68}
\qqq
(here we used an obvious fact that $\PM f_0= f_0$).

Let us remove all operators $\Chsn$ from the expressions for
$\BCMLLpunps$ and $\BCMLunv$. We denote the resulting operators as
$\BMLLpunps$ and $\BMLunv$ (note that $\BMLunv$ no longer depends on
$\uun$).
%
Since the change in the index $m_0$ of the vector
$f_{m_0}\in\Vzi$, as it passes through the
operators\rx{2.54},\rx{2.56},\rx{2.59} and\rx{2.63}, is fixed by the
Kronecker symbols, then the removal of the operators $\Chsn$ from the
\rhs of \ex{2.68} can be compensated by dividing it by the product of
factorials coming from \ex{2.51*1}
\qq
&\BMLLpun = \Fung
\lrbar{\Dsue \Tr_{\Vzi} (\Aot{I}{N}\otimes\PM\lstar^{2\hm})
\,\BCMLLpunps}{{\lstar=0\atop \ut=\qua}\atop \uue = 0},
\label{2.68*}\\
&F(\ung) = \lrbc{ \pioNg {|\ung|-\sjoimo \ng_j\choose \ng_i} }^{-1}.
\qqq

It remains to `parametrize' the expansion\rx{2.41} of the matrix
elements of the operator $\Qulm$. We introduce the numbers
$\Ckn\in\IQ$ according to the formula
\qq
{1\over n!}\,\lrbc{(1+h)^{-1/2} -1}^n =
\skgz \Ckn\,h^k.
\label{2.69}
\qqq
Note that
\qq
C_{0,0} = 1
\label{2.69*}
\qqq
Since
\qq
\l^{2m}\,q^{-m} & = &
\lrbc{\l + \lrbc{ (1+h)^{-1/2} -1}\l }^{2m}
\nonumber\\
& = & \lrbc{\snzi
{ \lrbc{ (1+h)^{-1/2} -1}^n \over n!}\,
(\l\,\del_\l)^n
}
\,\l^{2m}
\nonumber\\
& = &
\lrbc{\snkzi \Ckn\,h^k\,(\l\,\del_\l)^n}\,\l^{2m},
\label{2.70}
\qqq
then the \rhs of the following equation represents the expansion of
the matrix elements of $\Qqom$ described by \ex{2.41}
\qq
\QulomN\,\Qqom  =  \Duulh \,\QulomN,
\label{2.71}
\qqq
where
\qq
\Duulh =
\snkzi \Ckn\,h^k\,\l^n\,\del^n_{\l}.
\label{2.72}
\qqq

Consider a product of operators
\qq
\Dtot = \Duulh \Dsue
\label{2.73}
\qqq
\begin{proposition}
\label{prop.strder}
The operator $\Dtot$ is a power series in $h$
\qq
&\Dtot = \snzi \Dntot\,h^n,\hspace{3pt}
\Duunn(\ut;\l;\ube;\uue;x)
\in \IQ[\sqi{\ut},\l,\ube,\uue,x],
\nonumber\\
&D_{\uu{0},0}(\ut;\l;\ube;\uue;x) = 1,
\label{2.74*}\\
&\deg_{\ue} \Duunn + \deg_{x} \Duunn \leq 2n,
\label{2.74}\\
& \Duunn(\ut;\uul;\ube;\uue;x) =
\sumzi \Dumn(\uul;\ube;\uue;x)\,
(\ut-1)^{\um},
\label{2.75}\\
&\deg_x\Dumn + \deg_{\ube}\Dumn + \deg_{\uep}\Dumn \leq 2n + |\um|.
\label{2.76}
\qqq
and its action on the trace of the operator
$\QulomN\,\BMLLpunps$ yields
the expansion\rx{2.43*5}, so that
\qq
\lefteqn{\AMbun{\ut} = \Fung
\Dntot}&&
\label{2.77}\\
&&\times\lrbar{
\Tr_{\VBLaiz}\PLMs\Quloms\BMLunv \pjoLp
\Tr_{\VBLpjubip}\PLpM\,\Qpulj \BMLjpunv
}{\uue = 0\atop \lstar=0},
\nonumber
\qqq
where
\qq
\PLMs = \PM\otimes\PLM,\qquad
\Quloms = \lstar^{2\hm}\otimes\Qulomi{|\uN|}
\label{2.77*}
\qqq
and
\qq
\Qpulj = \left\{
\begin{array}{cl}
\brot{\l^{2\hm}}{|\uN\p|} & \mbox{if $j>1$ or $\cL$ is non-empty,}\\
\lo^{2\hm}\otimes\brot{\l^{2\hm}}{N\p_j} & \mbox{ if $j=1$ and
$\cL$ is empty}
\end{array}
\right.
\label{2.77**1}
\qqq
\end{proposition}

\proof
Equation\rx{2.74*} follows from \eex{2.66*} and\rx{2.69*}.
The bounds\rx{2.74} and \rx{2.76} follow from
Proposition\rw{prop.oper} and from \ex{2.72}. Equation\rx{2.68}  says
that its \rhs represents the expansion of the matrix elements of its
\lhs in powers of $h$, and according to Proposition\rw{expand},
that expansion leads to the expansion\rx{2.43*6}. Therefore
\qq
\lefteqn{\AMbun{\ut} = \Fung\,\Dntot
}&&
\label{2.77*1}\\
&&\times
\lrbar{
\Tr_{\Vzi\otimes\Vtip} (\PLMs\otimes\PLpM)\Qulmstar \BMLLpunps
}{\uue = 0\atop \lstar=0}
\nonumber
\qqq
According to \ex{2.65*}, the trace of this equation can be presented
as a product of traces, and this leads us to \ex{2.77}.\qed

\subsubsection{Calculating the traces}
\label{s2.2.5}
Our next goal is to calculate the traces of \ex{2.77} in the limit of
$M\rightarrow \infty$. We start with the trace over $\VBLaiz$. We
will calculate it with the help of the following
\begin{lemma}
\label{lem.geom}
Let $\Oh$ be an operator acting on a space $\ICN$ such that the
absolute values of its eigenvalues are less than 1. We extend $\Oh$
to a homorphism $O$ of the algebra $\ICzozN$. Then there exists a
limit
\qq
\limMi \Tr_{\ICzozN} \PM O =
{1\over \det_{\ICN}(1-\Oh)},
\label{2.78}
\qqq
where $\PM$ is a projector onto the polynomials of degree less than
$M$.
\end{lemma}
\proof
Equation\rx{2.78} is a multi-dimensional generalization of the
formula for the sum of a geometric series.\qed

We identify a Verma module $\Vai$ with the algebra $\ICz$
\qq
\fm\leftrightarrow z^m.
\label{2.79}
\qqq
Then it is easy to verify that the
operators\rx{2.42},\rx{2.43},\rx{2.50} and\rx{2.51} are algebra
homomorhisms generated by the following $2\times 2$ matrices
\qq
&\hRp\tteee  =  \pmatrix{
e^{\e_1+\e_3}(1-t_2^{-1}+\e_4) & e^{\e_2}t_1^{-1}\cr
e^{\e_1} & 0},
\label{2.80}\\
&\hRm\tteee = \pmatrix{
0 & e^{\e_2}\cr
e^{\e_1}t_2 & e^{\e_2+\e_3} t_1^{-1} t_2 (1-t_1+\e_4) },
\label{2.81}\\
&\hRg\eeg = \pmatrix{ 1 & 0\cr \e_2 & e^{\e_1} },
\label{2.82}\\
&\hRl\eel = \pmatrix{1 & \e_2 e^{\e_1}\cr 0 & e^{\e_1}}.
\label{2.83}
\qqq
Therefore, we conclude that the whole matrix $\BMLunv$ is a
homomorphism of the algebra $\ICzzzN$ which is generated by a
matrix $\hBMLunv$ acting on $\ICNo$. $\hBMLunv$ is
constructed by putting the operators $\hRpm\tteee$
and $\hRgl\eel$ acting on the appropriate basis vectors of $\ICNo$,
in place of the operators $\Rpm\tteee$ and $\Rgl\eel$ in the
expression for $\BMLunv$.

Since $\Quloms$ is also an algebra homomorphism generated by the
matrix
\qq
\hQulm = \diag{\lstar,\lo,\l,\ldots,\l},
\label{2.84}
\qqq
then we can calculate the first trace of \ex{2.77} in the limit of
$\Mi$ with the help of Lemma\rw{lem.geom}.
\begin{proposition}
There is a vicinity of
$\ue=0$, $\ut=1$, $\lstar,\lzl=0$ in which there
is a uniform convergence
\qq
\lMi \Tr_{\VBLaiz}\PLMs\,\Quloms\,\BMLunv
=
{1\over \det_{\ICNo}(I - \hQulm\hBMLunv) }.
\label{2.85}
\qqq
\end{proposition}
\proof
An inspection the matrix elements of the operators $\Quloms$ and
$\BMLunv$ indicates that there is a vicinity of $\ue=0$, $\ut=1$,
$\lstar,\lo,\lila=0$ in which the absolute values of the eigenvalues
of the operator $\hQulm\hBMLunv$ are less than $1$. Then \ex{2.85}
follows from \ex{2.78}.\qed

We denote
\qq
\Qzutule = \lrbar{\det\nolimits
_{\ICNo}(I - \hQulm\hBMLunv)}{\ue,\lstar=0},
\label{2.86}
\qqq
because we will show that this determinant is indeed the polynomial
mentioned in Proposition\rw{prop.rat}.

\begin{proposition}
The function $\Qzutule$ satisfies the properties\rx{2.87}
and\rx{2.88}.
\end{proposition}
\proof
Relation\rx{2.87} follows from the structure of the matrix elements
of matrices\rx{2.80}--(\ref{2.83}) and\rx{2.84}. To prove \ex{2.88},
note that
\qq
\lrbar{\hRpm\tteee}{\ue=0} = \rhopm\tontw,
\qquad
\lrbar{\hRgl\eeg}{\e_1,\e_2=0} = I.
\label{2.89}
\qqq
Therefore $\hBMLunv$ has a block form at $\ue=0$
\qq
\hBMLunv = \pmatrix{ 1 & 0\cr 0 & \BL[\ut]},
\label{2.90}
\qqq
where $\BL[\ut]$ is the matrix of Burau representation. Similarly
\qq
\lrbar{\hQulm}{\lslo} = \pmatrix{0 & 0\cr 0 & Q},
\label{2.91}
\qqq
where $Q$ is defined by \ex{1.17*1*1}. Then \ex{2.88} follows from
\ex{1.17*2}.\qed

Now it remains to calculate the traces
\qq
\Tr_{\VBLpjubip}\PLpM\,\Qpulj\, \BMLjpunv
\label{2.91*}
\qqq
of \ex{2.77}. Consider the product
of matrix elements of operators which constitute
$\Qpulj\,\BMLjpunv$.
Suppose that we start with a vector $f_{m_j}$ of
the space $\Vbj$, which is the first space in the second tensor
product of\rx{2.21**1}. Since the choice of the numbers $\uun$ fixes
the changes in $m_j$ as this vector evolves along the strands of
$\cLj$, then it is easy to see that, if $M>\b_j + 2|\uun|$, then
the projectors $\PM$ can be removed from the
expression for $\BMLjpunv$ and the projector $\PLpM$ can be removed
from\rx{2.91*}
without affecting the trace. As a result,
the product of the matrix elements has the form $X_j\,Y_j^{m_j}$
%
%
where $X_j(\elt)$ and $Y_j(\elt)$ are the
products of non-negative powers
of $e^{\uep}$, $\lzl$, $\ut$ and $\ut^{-1}$ which depend on $\uun$
(we do not indicate this dependence explicitly in order to simplify
notations).
In particular,
\qq
\lrbar{X_j(\elt)}{\lepcl} = 1,\quad
\lrbar{Y_j(\elt)}{\lepcl} = \pioL t_i^{-l_{ij}}\quad
\mbox{if $\uun=0$.}
\label{2.93*}
\qqq
Calculating the trace
over $\Vbj$ amounts to taking a sum over $m_j$, so
\qq
\lMi
\Tr_{\VBLpjubip}\PLpM\,\Qpulj\, \BMLjpunv
= \smobj X_j\,Y_j^{m_j} = X_j { Y_j^{\b_j} - 1\over Y_j -1 }.
\label{2.93}
\qqq
%
\begin{proposition}
\label{prop.dernum}
For non-negative integers $l,\ull$ the expression
\qq
Z_j(\ube;\ut) = \lrbar{ \del^{l}_{\l} \del^{\ull}_{\uep}
\lrbc{ X_j { Y_j^{\b_j} - 1\over Y_j -1 } }
}{\lepcl}
\label{2.94}
\qqq
has two properties: for fixed $\ube$
\qq
Z_j(\ube;\ut) \in \ZZutti
\label{2.93*1}
\qqq
and its expansion in powers of $\ut-1$ has the form
\qq
&Z_j(\ube;\ut) = \sumzi Z_{j,\um}(\b_j)\,(\ut-1)^{\um},
\nonumber\\
& Z_{j,\um}(\b_j) \in \IQbj, \qquad
\deg Z_{j,\um} \leq |\um| + |\ull| + l +1.
\label{2.96}
\qqq
\end{proposition}
\proof
The property\rx{2.93*1} is easy to see if we reverse \ex{2.93}
\qq
Z_j(\ube;\ut) = \lrbar{ \del^{l}_{\l} \del^{\ull}_{\uep}
\lrbc{ X_j \smobj Y_j^{m_j} }
}{\lepcl}.
\qqq
In order to prove the property\rx{2.96} we expand the
expression\rx{2.93} inside the \rhs of \ex{2.94} in powers of $Y_j-1$
\qq
Z_j(\ube;\ut) = \lrbar{ \del^{l}_{\l} \del^{\ull}_{\uep}
\lrbc{ X_j \snzi {\b_j\choose n+1 }(Y_j-1)^n }
}{\lepcl}.
\label{2.96*}
\qqq
Each of the derivatives $\del_{\l}$, $\del_{\uep}$ reduces the
power of $Y_j-1$ at most by one unit. Therefore the bound\rx{2.96}
follows from \ex{2.96*}, if we observe that since
$\lrbar{Y_j}{\lepcl}$ is a product of the powers of $\ut$,
then the expansion of $\lrbar{(Y_j-1)}{\lepcl}$ in powers of $\ut-1$
starts at degree 1.\qed

Now we have all the ingredients
which are necessary to finish the proof of
Proposition\rw{prop.rat}.
\pr{Proposition}{prop.rat}
Combining \eex{2.93} and\rx{2.85} with \ex{2.77} we conclude that
\qq
\lefteqn{
\lMi\AMbun{\ut} = \Fung
}&&
\label{2.97*1}\\
&&\times
\lrbar{
\Dntot\,
{\pjoLp
X_j { Y_j^{\b_j} - 1\over Y_j -1 }
\over \det_{\ICNo}(I - \hQulm\hBMLunv) }\,
}{\lstar,\uue=0}.
\nonumber
\qqq
Since the sum in the \rhs of \ex{2.41**1} is finite, then the limit
of the \lhs at $\Mi$ is equal to the sum of limits
\qq
\lefteqn{
\lMi\AMb{\ut} = \suneqn\Bigg(\Fung
}&&
\label{2.97}\\
&&\hspace{1in}\times
\lrbar{
\Dntot\,
{\pjoLp
X_j { Y_j^{\b_j} - 1\over Y_j -1 }
\over \det_{\ICNo}(I - \hQulm\hBMLunv) }\,
\Bigg)}{\lstar,\uue=0}.
\nonumber
\qqq

The operators $\Dntot$ are combinations of derivatives,
their structure
is described by Proposition\rw{prop.strder}. The action of the
derivatives $\del_{\l}$ and $\del_{\ue}$ on the denominator in the
\rhs of \ex{2.97} will increase its power and bring its derivatives
to the numerator. Hence the power of $\Qzutule$ in the denominator of
\ex{2.22*3*2} follows from the bound\rx{2.74}.
The action of the derivatives $\del_{\uep}$, $\del_{\l}$ on the
numerator is described by Proposition\rw{prop.dernum}. The order of
the derivatives $\del_{\uep}$ at each monomial of $\ut-1$ in the
expansion\rx{2.75} of $\Dntot$ is restricted by\rx{2.76}, hence the
bound\rx{2.22*3*4}.

Thus we proved that the \rhs of
\ex{2.97} has the same
form as the \rhs of \ex{2.22*3*2} except that we can only claim that
\qq
\Qnubut \in \IQutuul.
\label{2.98}
\qqq
To prove the integrality of the coefficients of $\Qnubut$ we rewrite
\ex{2.22*3*2} as
\qq
\Qnzutule \lMi \AMb{\ut} = \Qnubut.
\label{2.99}
\qqq
Now the integrality of the coefficients of $\Qnubut$ follows from
the integrality of the coefficients of $\Qzutule$ and
from\rx{2.22**5*1}.
Finally, relation\rx{2.88*} follows from \eex{2.74*} and\rx{2.93*}.
\qed
\subsection{Removal of parametrization}
Now it's time to go back to the original (unparametrized and
unregularized) Jones polynomial by setting $\lo=0$, $\l=1$ in
\ex{2.22*3*2}. We define a formal power series
\qq
\lefteqn{\JhrubLLputh =}&&
\label{2.100}\\
&&
\hspace{-1pt}
\left\{
\begin{array}{cl}
\displaystyle
{1\over h}\,q^{
\phlkLpb
+ \ftwo(\ube)}\,
{\PhubLpt\over\AFLut}
\snzi {\lrbar{\Qnubut}{\ulzcl}\over
\lrbc{\PLouti{-1}\,(\thtmh)\,\AFLut}^{2n} }\,h^n, &
\mbox{if $L\geq 1$,}\\
\displaystyle
{1\over h}\,q^{\phlkbas{\empt,\cLp}{\ube}
 + \ftwo(\ube)}\,(1-q^{-\b_1})
\snzi \lrbar{Q_{\ube;n}(\empt;\cLp;\ul)  }{\ulzcl}\,h^n
&\mbox{if $L=0$,}
\end{array}\right.
\nonumber
\qqq
where
\qq
&\ftwo(\ube)  =  \hlf
\lrbc{
\sjoL (l_{jj} - N_j - 1)
-\sjoLp (l\p_{jj} - N\p_j -1)(\b_j-1)
}
+ 1,
\label{2.101}
\qqq
while
$\phlkLpb$ and $\PhubLpt$ are defined by
\eex{1.45**1} and\rx{2.103}.

We used the same notation for the series\rx{2.100} as in \ex{1.35},
because we will ultimately show that both expressions are the same.
\begin{proposition}
Formal power sereis\rx{2.100} satisfies three properties of
$\JhrubLLputh$ (\eex{1.38},\rx{1.39} and\rx{1.42}) listed in Main
Theorem, except that if $L\geq 2$, then \ex{1.42} is satisfied only
if $i\neq 1$.
\end{proposition}
\proof
Equations\rx{1.38},\rx{1.39} follow from \eex{2.22*3*6}
and\rx{2.22*3*7} respectively if we write the latter equations at
$\lo=0$, $\l=1$.
To prove \ex{1.42} we have to consider two cases.

Suppose that $L\geq 2$ and $i\neq 1$.
According to \ex{2.88} and Torres formula\rx{1.15*},
\qq
\lrbar{\Qzutule}{t_i=1\atop{\ulzcl}}
& = & \lrbar{\PLouti{-1}\,(\thtmh)\,\AFLut}{t_i=1}
\label{2.104}\\
& = & \lrbar{\PLouti{-1}\,(\thtmh)}{t_i=1}
\!\lrbc{ \pjoLmo t_j^{l_{ij}/2} -
\pjoLmo t_j^{-l_{ij}/2} }\!
\AFbas{\cL\xrem{i}}{\ut\xrem{i}}
\nonumber
\qqq
According to our assumptions,
\qq
\lrbc{ \pjoLmo t_j^{l_{ij}/2} - \pjoLmo t_j^{-l_{ij}/2} }\not\equiv
0, \qquad
\AFbas{\cL\xrem{i}}{\ut\xrem{i}}\not\equiv 0.
\qqq
Hence, according to \ex{2.104},
$\lrbar{\Qzutule}{t_i=1\atop{\ulzcl}}\not\equiv 0$ and therefore
\ex{1.42} follows from \ex{2.22*3*5}.

Suppose now that $L=1$.
According to \eex{2.88},\rx{1.12} and\rx{1.10},
\qq
\lrbar{\Qili{0}{\cL}{t_1} }{t_1=1\atop{\ulzcl}} =
\lrbar{\Phibasii{-1}{1}{\cL}{t_1} \APbas{\cL}{t_1}}{t_1=1} = 1
\not\equiv 0.
\label{2.105}
\qqq
Therefore we can apply \ex{2.22*3*5}. Equation\rx{2.22**5**4}
indicates that
\qq
\lrbar{\Foi{-\a_1}}{\lo=0} = 0\qquad\mbox{if $\a_1>0$,}
\label{2.106}
\qqq
so that we can drop the term with $\mu=-1$ from the
\lhs of \ex{2.22*3*5} if $\lo=0$. Thus we arrive at \ex{1.42} for
$L=1$.\qed

\begin{proposition}
\label{prop.den}
There exists a set of polynomials $\tPubenLLput\in\ZZutthih$ such that
\qq
&\tPubenLLput \in \ZZutti,
\label{2.107}\\
&\tPubezLLput =
\pjoLp{ 1 - \ppioL{-\b_j} \over 1 -
\ppioL{-1}},
\label{2.106*}
\qqq
an expansion of $\tPubenLLput$ in powers of $\ut-1$ has a form
\qq
&\tPubenLLput = \sumzi \tPubenumLLput\,(\ut-1)^{\um},
\nonumber\\
& \tPubenumLLput \in \IQube,\qquad
\deg \tp_{\um;n} \leq 2(|\um|+n) + 1.
\label{2.108}
\qqq
and
\qq
\lefteqn{
\JhrubLLputh}
\label{2.109}
\\
&&
 =
\left\{
\begin{array}{cl}   \displaystyle
q^{\phlkLpb}\snzi \tPubenLp\,h^n
&\mbox{if $L=0$,}
\vspace{5pt}
\\
\displaystyle
{1\over h}\,
q^{
\phlkLpb}
{\PhubLpto\over \AFbas{\cL}{t_1}}
\snzi {\tPbas{\ube}{n}{\cL}{\cLp}{t_1}\over \APpbas{2n}{\cL}{t_1} }
\,h^n
&\mbox{if $L=1$,}
\vspace{5pt}
\\
\displaystyle
{1\over h}\, q^{
\phlkLpb}
{\PhubLpt\over \AFLut}
\snzi
{\tPubenLLput\over
\bigg(
(t_1-1)\,\AFLut
\bigg)^{2n} }
\,h^n
&\mbox{if $L\geq 2$,}
\end{array}
\right.
\nonumber
\qqq
\end{proposition}

\proof
The expressions\rx{2.109} can be derived directly from \ex{2.100}. If
$L=0$, then we have to expand the factor
$q^{\fth(\ube)+ 1}\,(1-q^{-\b_1})$ in powers of $h$. If
$L\geq 1$, then we expand the factor
$q^{\ftwo(\ube)}$
in powers of $h$ and combine the factors
$\PLouti{-2n}$ with the polynomials $\lrbar{\Qnubut}{\ulzcl}$.
Equation\rx{2.106*} follows from \eex{2.100},\rx{2.103}
and\rx{2.88*}.
The property\rx{2.107} follows from\rx{2.87} and from
%
the fact that $\ftwo(\ube)\in\ZZ$, because it is not hard to prove
that
\qq
l_{jj} - N_j - 1,\; l\p_{jj} - N\p_j - 1 \in \ZZ.
\label{2.109*}
\qqq
%
The bound\rx{2.108} on the degree of the
polynomials $\tP_{\um;n}$ follows easily from the similar
bound\rx{2.22*3*4}.\qed

\nsection{Uniqueness arguments}

\subsection{Proof of Main Theorem}

The claims of Proposition\rw{prop.den} are very similar to those of
Main Theorem. In fact we will prove all claims of Main Theorem
if we show two things: the formal power series $\JhrubLLputh$ is an
invariant of $\ccLLp$ (that is, its value does not depend on a
choice of a braid $\cBLLp$ used to calculate it) and if $L\geq 2$,
then polynomials $\tPubenLLput$ are divisible by $(t_1-1)^{2n}$.

\begin{theorem}[Uniqueness Theorem]
Let $\bcLs$ be a subset of the set of all links, such that the empty
link and all 1-component links belong to it, and if
$\cL\in\bcLs$ and a knot $\cL_0$ is algebraically connected to $\cL$,
then $\cL_0\cup\cL\in\bcLs$.
%
Suppose that for any pair of
oriented links $\cLp$, $\cL\in\bcLs$ and for
any number $l$, $1\leq l\leq L$, there exists a non-empty set $\CLl$
such that for any $c\in\CLl$ there exist formal power series
%
\qq
\GcLua \in \IQbua,\qquad
\FcnLab\in\IQbuab,\quad n\geq 0,
\label{3.1}
\qqq
which satisfy the following properties:
\qq
\GcLua\not\equiv 0,
\label{3.2}
\qqq
if
$\cL_j$ is algebraically connected to $\cL\xrem{j}$ then
\qq
\lrbar{\GcLua}{a_j=0}\not\equiv 0
\label{3.3}
\qqq
and a formal power series
\qq
\IcLLabe =
\left\{
\begin{array}{cl}
\displaystyle
\snzi \Fbas{c}{n}{\cLp}{\ube}\,\e^n &
\mbox{if $L=0$,}\\
\displaystyle
\snzi {\Fbas{c}{n}{\ccLLp}{a_1;\ube} \over
\lrbc{ \Gbas{c}{\cL}{a_1} }^{2n+1} }\,\e^n &
\mbox{if $L=1$,}
\\
\displaystyle
\snzi \,
{\FcnLab\over \lrbc{a_l\,\GcLua}^{2n+1} }\,
\e^n &
\mbox{if $L\geq 2$.}
\end{array}
\right.
\label{3.4}
\qqq
satisfies the following properties
\begin{itemize}
\item[(1)]
If $\cL$ is empty, then $\Ibas{c}{\empt,\cLp}{\ube}$ is an invariant
of an oriented link $\cLp$, that is, it is the same for all
$c\in\Ci{\cLp}{}$.

\item[(2)]
For any $j$, $1\leq j\leq L\p$ and for any $c\in\CLl$ there exists
$c\p\in \Ci{\cL,\cLpro}{l}$ such that
\qq
\lrbar{\IcLLabe}{\b_j=1} =
\Ibas{c}{\cL,\cLp\xrem{j}}{\ua;\ube\xrem{j};\e}.
\label{3.5}
\qqq

\item[(3)]
Suppose that
for a number $i$, $1\leq i\leq L$, $\cL\xrem{i}\in\bcLs$ and
$\cL_i$ is algebraically connected to $\cL\xrem{i}$.
Then for any $c\in\CLl$ there exists
$c\p\in\Ci{\cLro,\cL_L\cup\cLp}{l}$  and for any
$c\p\in\Ci{\cLro,\cL_L\cup\cLp}{l}$ there exists $c\in\CLl$ such
that
%
%
\qq
\lefteqn{
\Ibas{c\p}{\cL\xrem{i},\cL_j\cup\cLp}{\ua\xrem{i};\b_0,\ube;\e} }
&&\label{3.7}\\
& = &
\left\{
\begin{array}{cl}
\displaystyle
\Ibas{c}{\ccLLp}{\e\b_0;\ube;\e} &
\mbox{if $L=1$,}\\
\displaystyle
\smumu\,
e^{i\pi\mu\b_0\sjoLni l_{ij}a_j}
\lrbar{\IcLLabe}{a_i=\e\b_0} &
\mbox{if $L\geq 2$ and $i\neq l$.}
\end{array}
\nonumber
\right.
\qqq
\end{itemize}
Then
\begin{itemize}
\item[(1)]
The formal power series $\IcLLabe$ is an invariant
of a pair consisting
of an oriented link $\cL\in\bcLs$ and a link $\cLp$,
that is, for any two numbers $l_1$, $l_2$,
$1\leq l_1,l_2 \leq L$ and for two elements
$c_{1,2}\in\Ci{\ccLLp}{l_{1,2}}$,
\qq
\Ibasi{1} = \Ibasi{2}.
\label{3.7*}
\qqq

\item[(2)]
The invariant $\Ibasws$ is an odd function of $\ube$, that is,
for any $j$, $1\leq j\leq L\p$
\qq
\Ibasw{\ccLLp}{\ua;\b_1,\ldots,-\b_j,\ldots,\b_{L\p}} =
- \Ibasw{\cLLp}{\ua;\ube}
\label{3.7*3}
\qqq
(we dropped the index $(c)$, because in view of \ex{3.7*}, $\IcLLabe$
does not depend on it).

\item[(3)]
If $\cL$ is algebraically connected, then
the formal power series $\FcnLab$ is divisible by $a_l^{2n+1}$
in $\IQbuab$.

\item[(4)]
If the series $\GcLua$ is an invariant of $\cL$, then for any link
$\cL\in\bcLs$ the series $\FcnLab$ is divisible by $a_l^{2n+1}$ in
$\IQbuab$.


\end{itemize}
\end{theorem}

\begin{corollary}
\label{cor.uniq}
Suppose that in addition to the conditions of Uniqueness Theorem the
invariant $\Ibasw{\empt,\cLp}{\ube}$ does not depend on the
orientation of $\cLp$. Then
\begin{itemize}

\item[(1)]
The invariant $\Ibasws$ does not depend on the orientation of $\cLp$.

\item[(2)]
Let $\cLb$ be a link which differs from $\cL\in\bcLs$ only by the
orientation of its component $\cL_j$. Then
\qq
\Ibasw{\cLb,\cLp}{\ua;\ube}
= -
\Ibasw{\cLLp}{a_1,\ldots,-a_j,\ldots,a_L;\ube}
\label{3.7*1}
\qqq

\item[(3)]
If we take an oriented link $\cL\in\bcLs$ and reverse the
orientations of all of its components, then for the resulting link
$\cLb$
\qq
\Ibasw{\cLb,\cLp}{\ua;\ube} = \Ibasw{\cL,\cLp}{\ua;\ube}
\label{3.7*2}
\qqq

\end{itemize}
\end{corollary}

We will prove Uniqueness Theorem
and Corollary\rw{cor.uniq} in the next subsection. Here we will
show how they help to derive Main Theorem from
Proposition\rw{prop.den}.

\noindent
\emph{Proof of Main Theorem}
Let $\bcLs$ be a set of all links whose Alexander-Conway function is
not identially equal to zero. Let $\CLl$ be the set of braids
$\cBLLp$ whose closure is $\cLLp$. Finally, let $\IcLLabe$ be a formal
power series $h\JhrubLLputh$ in which we perform the following
substitution
\qq
\ut = e^{\ua}, \qquad
h = e^{\e} - 1.
\label{3.8}
\qqq
A simple manipulation with \ex{2.109} brings it to the form\rx{3.4}
with
\qq
\GcLua =
\left\{
\begin{array}{cl}
\displaystyle
\APbas{\cL}{e^{a_1}} &
\mbox{if $L=1$,}\\
\displaystyle
\AFbas{\cL}{e^{\ua}} &
\mbox{if $L\geq 2$.}
\end{array}
\right.
\label{3.9}
\qqq
Now a simple manipulation brings $\JhrubLLputh$ of \ex{2.109} to the
form\rx{3.4}.
Therefore Proposition\rw{prop.den} implies that all conditions of
Uniqueness Theorem hold. In its turn, Uniqueness Theorem says that
for $L\geq 2$, $\tPubenLLp{e^{\ua}}$ has to be divisible by
$a_1^{2n}$. In view of \rx{2.107} this means that
\qq
{\tPubenLLput\over (t_1-1)^{2n}}\in
\ZZutti.
\label{3.9*}
\qqq
Thus we can set
\qq
\PpubenLLput =
\left\{
\begin{array}{cl}
\displaystyle
\tPubenLLput &
\mbox{if $L=0,1$,}\\
\displaystyle
{\tPubenLLput \over (t_1-1)^{2n}} &
\mbox{if $L\geq 2$,}
\end{array}
\right.
\label{3.10}
\qqq
and \ex{2.109} will turn into \ex{1.45**0} with the relation\rx{3.9}
guaranteeing\rx{1.45**2}.

Equation\rx{1.35} follows from \ex{1.45**0} if we combine $\PhubLpt$
with the polynomials $\PpubenLLput$ and expand the factor
$q^{\phlkLpb}$ in powers of $h$. Equation\rx{1.43} follows from
\ex{2.106*}. Relation\rx{1.33}
follows from\rx{1.45**2}
because
\qq
\PhubLpt\in\ZZutthi,\qquad \phlkLpb\in \hlf\,\ZZ.
\label{3.10*}
\qqq

The property\rx{1.38} of $\Jhremp$ together with the independence of
the colored Jones polynomial $\JubLp$ of the orientation of the link
$\cLp$ implies that the condition of Corollary\rw{cor.uniq} is
satisfied. Therefore Uniqueness Theorem and
Corollary\rx{cor.uniq} together say that the formal power series
$\JhrubLLputh$ (and hence the polynomials $\PubenLLput$) is an
invariant of an oriented link $\cL$ and a link $\cLp$.

Relations\rx{1.44},\rx{1.45} follow from \eex{3.7*1}, \rx{3.7*2} and
from the properties\rx{1.14},\rx{1.15} of the Alexander-Conway
function. Relation\rx{1.45*2} follows from \ex{3.7*3}.\qed

\subsection{Proof of Uniqueness Theorem}

\subsubsection{Algebraically connected links}
First of all, we will prove an obvious lemma.
\begin{lemma}
The set $\bcLs$ of Uniqueness Theorem contains all algebraically
connected links.
\end{lemma}
\proof
We use induction in the number $L$ of link components. Since $\bcLs$
contains the empty link, then lemma holds for $L=0$. Suppose that all
the algebraically connected links with $L-1$ components belong to
$\bcLs$. Let $\cL$ be an $L$-component algebraically connected link.
Then $\cLro\in\bcLs$. Since, by definition, $\cL_L$ is algebraically
connected to $\cLro$, then, according to the property of $\bcLs$,
$\cL\in\bcLs$.\qed

\begin{proposition}
\label{prop.algcon}
Claims~1 and~3 of Uniqueness Theorem hold for algebraically connected
links.
\end{proposition}
We need two lemmas in order to prove this proposition.
\begin{lemma}
\label{lem.den}
If the series $\IcLLabe$ is unique for an $L$-component algebraically
connected link $\cL$
with $L\geq 2$, then $\FcnLab$ is divisible by $a_l^{2n+1}$ in
$\IQbuab$.
\end{lemma}
\proof
For a number $l$, $1\leq l\leq L$ and an element $c\in\CLl$
choose another number $l\p\neq l$, $1\leq l\leq L$ and an element
$c\p\in\Ci{\ccLLp}{l\p}$. Uniqueness means that for any $n\geq 0$
\qq
{\Fbasi{c}\over \Big(a_{l}\GcLua\Big)^{2n+1} } =
{\Fbasi{c\p}\over \Big(a_{l\p}\Gbas{c\p}{\cL}{\ua}\Big)^{2n+1} },
\label{3.11*}
\qqq
or, equivalently,
\qq
\Fbasi{c} \Big(a_{l\p}\Gbas{c\p}{\cL}{\ua}\Big)^{2n+1}
=
\Fbasi{c\p} \Big( a_{l}\GcLua\Big)^{2n+1}.
\label{3.12*}
\qqq
The \rhs of this equation is divisible by $a_l^{2n+1}$, while
$\Gbas{c\p}{\cL}{\ua}$ is not divisible by $a_l$, because, since
$\cL$ is algebraically connected, then $\cL_l$ is algebraically
connected to $\cL\setminus \cL_l$ and as a result, according to
\ex{3.3},
\qq
\lrbar{\Gbas{c\p}{\cL}{\ua}}{a_l}\not\equiv 0.
\label{3.13*1}
\qqq
Therefore $\Fbasi{c}$ has to be divisible by $a_l^{2n+1}$.\qed

\begin{lemma}
\label{lem.uniq}
Let $L\geq 2$. Assume that the series $\IcLLabe$ is unique for all
algebraically connected $L-1$-component links. Suppose that the
components of an $L$-component link $\cL$ are indexed in such a way
that the link $\cLro$ is also algebraically connected. Then
for any two numbers $l_1,l_2\neq L$ and for
any two elements $c_1\in\Ci{\cLLp}{l_1}$, $c_2\in\Ci{\cLLp}{l_2}$
\qq
\Ibasi{1} = \Ibasi{2}.
\label{3.11}
\qqq
\end{lemma}
\proof
Since $\cL$ is algebraically connected, then $\cL_L$ is algebraically
connected to $\cLro$ and then, according to \ex{3.3},
\qq
\lrbar{\Gbasi{i}}{a_L=0}\not\equiv 0,\qquad
i=1,2.
\label{3.12}
\qqq
Therefore we can expand $\Ibasi{i}$ in powers of
$a_L$
\qq
&\Ibasi{i} = \smnzi
\Ibasii{i}
\,a_L^{m}\,\e^{n}, \qquad i=1,2,
\label{3.13}\\
&\Ibasii{i} = {\Hcmni{i} \over a_{l_i}^{2n+1}\,
\Big( \Gbasi{i} \Big)^{2n+m+1} },
\label{3.13*}
\qqq
where $\Hcmni{i}\in\IQ[\ube][\uaro]$. Equation\rx{3.11} would follow
from the relations
\qq
&\Iot = 0,\qquad m,n\geq 0,
\label{3.14}
\qqq
where
\qq
&\mbox{ $\Iot = \Ibasii{1}-\Ibasii{2}$.}
\label{3.14*}
\qqq

We are going to prove relations\rx{3.14}.
Since, according to our assumptions, $\cLro$ is algebraically
connected and
the series $\Ibas{c}{\cLro,\cLp}{\uaro;\ube}$ is unique,
then, in view of \ex{3.7},
\qq
\hspace{-30pt}\smuLmu\,
e^{i\pi\muaL\sjoLmo l_{jL}a_j}
\Big(\Ibas{c_1}{\ccLLp}{\uaro;\muaL;\ube;\e} -
\Ibas{c_1}{\ccLLp}{\uaro;\muaL;\ube;\e} \Big) = 0.
\label{3.15}
\qqq
Rewriting the series $\Ibas{c_i}{\ccLLp}{\uaro;\muaL;\ube;\e}$ with
the help of \ex{3.13}, we find that
\qq
\snzi \e^n
\lrbs{
\smuLmu\,e^{i\pi\muaL\sjoLmo l_{jL}a_j}
\lrbc{
\smzn
\Iotm\,\mu_L^m \a_L^m
}} = 0.
\label{3.16}
\qqq
As a function of $\a_L$, each of the two terms in the sum over
$\mu_L$ is a product of an exponential and a polynomial.
Since, according to our assumptions,
\qq
\muaL\sjoLmo l_{jL}a_j \not\equiv 0,
\label{3.17}
\qqq
then the exponents in both terms are different and therefore each
term in the sum over $\mu_L$ has to be equal to zero. Hence we come
to \ex{3.14} and to \ex{3.11}.\qed

\pr{Proposition}{prop.algcon}
The proof is by induction in the number of link components $L$.

If $L=0$, then the claim of Uniqueness Theorem is equivalent to
Property~1 of the series $\Ibas{c}{\cLp}{\ube}$.

If $L=1$, then the uniqueness of $\Ibas{c}{\cL_1,\cLp}{a_1;\ube}$
follows from \ex{3.7} and from the uniqueness of
$\Ibas{c\p}{\cLp}{\e a_1,\ube}$.

If $L=2$, then we have to consider two cases. For $l=1,2$ and for
$c_{1,2}\in\Ci{\ccLLp}{l}$ the uniqueness relation\rx{3.7}
follows from Lemma\rw{lem.uniq}. The case of $l_1=1$, $l_2=2$
requires a special consideration. Let us expand
\qq
{\Fbas{c_j}{n}{\ccLLp}{a_1,a_2;\ube}\over
\lrbc{\Gbas{c_j}{\cL}{a_1,a_2}}^{2n+1} }
= \smotzi \Hbasi{j}\,a_1^{m_1}a_2^{m_2}.
\label{3.18}
\qqq
Then
the uniqueness of the series $\Ibas{c}{\ccLLp}{a_1,a_2;\ube}$ would
follow from the relation
\qq
\Hbasii{1}{2n+m_1+1,m_2,n} =
\Hbasii{2}{m_1,2n + m_2 + 1, n}\qquad
\mbox{for all $m_1$, $m_2$ and $n$.}
\label{3.19}
\qqq

Let us apply the property\rx{3.7} to
$\Ibas{c_j}{\ccLLp}{a_1,a_2;\ube}$ twice
first reducing it to\\
$\Ibas{c\p_j}{\cL_j, \cL_k\cup\cLp}{a_j;\e \a_k,\ube}$,
where $1\leq k\neq j \leq 2$, and then reducing it further to
$\Ibas{c^{\prime\prime}_j}{\empt;\cLLp}{\e \a_1,\e \a_2,\ube}$.
We introduce a notation
\qq
\lefteqn{
\snmot
\XXbas{m_1,m_2,n}(\ccLLp;\ube)\, \a_1^{m_1} \a_2^{m_2} \e^n =
\smumu\,e^{i\pi\e\, l_{12} \a_1\a_2}
}
\label{3.20}\\
&&
\times
\snmotzi \Big( \Hbasi{1}\,\mu^{m_2} \a_1^{m_1-2n-1}
\a_2^{m_2}
\nonumber\\
&&\hspace{2in} +
\Hbasi{2}\, \mu^{m_1}\a_1^{m_1}\a_2^{m_2-2n-1}
\Big)
\e^{m_1+m_2-n-1}
\nonumber
\qqq
(note that only a finite number of coefficients $\Hbasi{j}$
contributes to any given coefficient $X$). The uniqueness of
$\Ibas{c^{\prime\prime}}{\empt;\cLLp}{\e \a_1,\e \a_2,\ube}$ implies
the following relation
\qq
\XXbas{m_1,m_2,n}(\ccLLp;\ube) = 0\qquad \mbox{for all $m_1,m_2,n$.}
\label{3.21}
\qqq
It is easy to see from \ex{3.20} that for any two integers
$m_1$, $m_2$
\qq
\lefteqn{
\sumn \XXbas{m_1-n,m_2-n,m_1+m_2-n}\,\a_1^{m_1-n}\a_2^{m_2-n}
\e^{m_1+ m_2 -n}=
\smumu
\,e^{i\pi\e\, l_{12} \a_1\a_2}
}
\label{3.22}\\
&&\times\sumn
\Big( \Hbasii{1}{m_1+n-1,m_2-n,n} - \Hbasii{2}{m_1-n,m_2+n-1,n}
\Big)
\a_1^{m_1-n}\a_2^{m_2-n}
\e^{m_1+ m_2 -n}.
\nonumber
\qqq
Therefore, the relation\rx{3.21} implies that the \rhs of \ex{3.22}
has to be identically equal to zero. Since the indices of the
coefficients $H$ have to be non-negative, then the sum over $n$ in
\ex{3.22} contains only finitely many terms. Therefore since
$l_{1,2}\neq 0$, then the \rhs of \ex{3.22} can be identically equal
to zero only if
\qq
\Hbasii{1}{m_1+n+1,m_2-n,n} - \Hbasii{2}{m_1-n,m_2+n+1,n}=0
\quad
\mbox{for any $m_1$, $m_2$ and $n$.}
\label{3.23}
\qqq
This condition is obviously equivalent to\rx{3.19}, hence we proved
the uniqueness of the series $\IcLLabe$ for 2-component links.

Suppose now that $L>2$ and that the series $\IcLLabe$ is unique for
all algebraically connected links with $L-1$ components. Consider an
$L$-component algebraically connected link $\cL$. Let us construct a
special graph $\G(\cL)$ with $L$ vertices identified with link
components. Two vertices corresponding to link components $\cL_{j_1}$
and $\cL_{j_2}$ are connected by an edge if there exists a third
component $\cL_{k}$ such that $k\neq j_1,j_2$ and
the sublink $\cL\setminus\cL_k$ is algebraically connected.
\begin{lemma}
\label{lem.con}
For an algebraically connected link $\cL$, the graph $\G(\cL)$ is
connected.
\end{lemma}
\proof We leave the proof to the reader.\qed

According to Lemma\rw{lem.uniq}, if two vertices $j_1$ and $j_2$ of
the graph are connected, then for any two elements
$c_1\in\Ci{\cLLp}{j_1}$, $c_2\in\Ci{\cLLp}{j_2}$ the corresponding
series $\IcLLabe$ are equal. Therefore, the uniqueness of this series
for $\cL$ follows from Lemma\rw{lem.con}. This completes the
induction argument for the proof of uniqueness of $\IcLLabe$ for
algebraically connected links. The divisibility of $\FcnLab$ by
$a_l^{2n+1}$ follows from Lemma\rw{lem.den}.\qed

\subsubsection{General links}

Now we can prove the claims of Uniqueness Theorem for any link
$\cL\in\bcLs$.

\noindent
\emph{Proof Claim~1 of Uniqueness Theorem.}
Let $\cL\in\bcLs$ be an $L$-component link. Let us add an extra
0-th component $\cLz$ such that the link $\cLLz$ is algebraically
connected. According to \ex{3.7}, for any integer $l$,
$1\leq l\leq L$ and for any $c\in\CLl$ there exists
$c\p\in\Ci{\cLz,\cLp}{l}$ such that
\qq
\Ibas{c}{\cL,\cLzLp}{\ua;\a_0,\ube;\e} =
\smuzmu\,
e^{i\pi\mu_0\a_0\sjoLmo l_{jL}a_j}
\Ibas{c\p}{\cLLz,\cLp}{\mu_0\a_0,\ua;\ube;\e}.
\label{3.24}
\qqq
Equation\rx{3.5} allows us to remove the extra component $\cLz$ from
the \lhs of \ex{3.24}
\qq
\IcLLabe =
\smuzmu\,
\lrbar{
e^{i\pi\mu_0\a_0\sjoLmo l_{jL}a_j}
\Ibas{c\p}{\cLLz,\cLp}{\mu_0\a_0,\ua;\ube;\e} }{\a_0=1}.
\label{3.25}
\qqq
The link $\cLLz$ is algebraically connected. Therefore, according to
Proposition\rw{prop.algcon}, the \rhs of \ex{3.25} does not depend on
$c\p$, hence the \lhs does not depend on $c$. This proves the
uniqueness of $\IcLLabe$.\qed

\noindent
\emph{Proof of Claim~2 of Uniqueness Theorem.}
Consider an oriented $L$-component link $\cL\in\bcLs$ and an
oriented
$L\p$-component link $\cLp$.
Choose a number $k$, $1\leq k\leq L\p$. Let us add an extra
component $\cL_0$ to $\cL$ in such a way that
$\lpp_{0k}=1$ and, in addition, $l_{0j}=1$ for all $j$,
$1\leq j\leq L$. It is easy to see that both $\cLLz$
and $\cLkLLz$ are algebraically connected and therefore
$\cLLz,\cLkLLz\in\bcLs$.

According to \ex{3.7},
\qq
\lefteqn{
\Ibasw{\cLLz,\cLp}{a_0,\ua;\ube}
}
\label{3.27}
\\
&&=
\smumu e^{i\pi\mu \b_k\lrbc{a_0 + \sjoL \lpp_{jk}a_j} }
\Ibasw{\cLkLLz,\cLp\xrem{k}}{\mu\e\b_k,a_0,\ua;\ube\xrem{k}}.
\nonumber
\qqq
On the other hand, a combination of \eex{3.5}
and\rx{3.7} implies that
\qq
\Ibasw{\cL,\cLp}{\ua;\ube} =
\lrbar{
\smumu e^{i\pi\mu \a_0\sjoL l_{0j}a_j}
\Ibasw{\cLLz,\cLp}{\mu\e\a_0,\ua;\ube} }{\a_0=1}.
\label{3.28}
\qqq
Equation\rx{3.27} indicates that
$\Ibasw{\cLLz,\cLp}{a_0,\ua;\ube}$ is an odd function of $\b_k$. Then
it follows from \ex{3.28}, that $\Ibasw{\cL,\cLp}{\ua;\ube}$ is also
an odd function of $\b_k$. \qed

\noindent
\emph{Proof of Claim~4 of Uniqueness Theorem.}
For a number $l$, $1\leq l\leq L$ and an element $c\in\CLl$
choose another number $l\p\neq l$, $1\leq l\leq L$ and an element
$c\p\in\Ci{\ccLLp}{l\p}$. Uniqueness means that for any $n\geq 0$
\ex{3.12*} holds. If $\GcLua$ is unique (that is, if it does not
depend on $l$ and $c$), then we can cancel it from both sides of that
equation, so that
\qq
\Fbasi{c}\, a_{l\p}^{2n+1}
=
\Fbasi{c\p}\, a_l^{2n+1}.
\label{3.12**}
\qqq
Divisibility follows from this equation.\qed

\subsection{Proof of Corollary\rw{cor.uniq}}

\pr{Claims~1 and 2 of Corollary}{cor.uniq}
Let $\cL\in\bcLs$ be an $L$-component oriented link and let $\cLp$
be an oriented $L\p$-component link. Consider a group $\GcLLp$ which
is generated by $L$ elements $g_j$, $1\leq j\leq L$ each of which
acts on $\cL$ by changing the orientation of the corresponding link
component $\cL_j$ and on $\ua$ by changing the sign of $a_j$, and by
$L\p$ elements $g\p_j$, $1\leq j\leq L\p$ each of which acts on
$\cLp$ by changing the orientation of the corresponding link
component $\cLp_j$.
For $g\in\GcL$ let $|g|$ denote the number of components of $\cL$
whose orientation it changes.

Consider a new invariant
\qq
\IGbas{\ccLLp}{\ua;\ube}
= 2^{-(L+L\p)}
\sgp\,
\Ibasw{g(\ccLLp)}{g(\ua),\ube}.
\label{3.30}
\qqq
It is not hard to see that $\IGbas{\cLLp}{\ua;\ube}$ satisfies the
conditions of Uniqueness Theorem (note that condition\rx{3.7} at
$L=1$ follows from Claim~2 of Uniqueness Theorem applied to
$\Ibasws$). Since we assume that $\Ibasw{\empt,\cLp}{\ube}$ does not
depend on the orientation of $\cLp$, then
\qq
\IGbas{\empt;\cLp}{\ube} = \Ibasw{\empt,\cLp}{\ube}.
\label{3.31}
\qqq
Therefore Uniqueness Theorem says that
\qq
\IGbas{\ccLLp}{\ua;\ube} = \Ibasws.
\label{3.32}
\qqq
Since the \rhs of \ex{3.30}
satisfies Claims~1 and~2 of Corollary\rw{cor.uniq}, then so does
the \rhs of \ex{3.32}.\qed

\pr{Claim~3 of Corollary}{cor.uniq}
The proof is similar to the previous one.
For an oriented link $\cL$ let $\cLb$ denote the same link with
reversed orientation of all its compomenents. This time we define
the new invariant $\IGbas{\ccLLp}{\ua;\ube}$ as
\qq
\IGbas{\ccLLp}{\ua;\ube} = \hlf \Big( \Ibasws +
\Ibasw{\cLb,\cLp}{\ua;\ube} \Big).
\label{3.33}
\qqq
It follows easily from Claim~1 of Corollary\rw{cor.uniq} that
$\IGbas{\ccLLp}{\ua;\ube}$ satisfies conditions of Uniqueness
Theorem and \ex{3.31}. Therefore, according to Uniqueness Theorem,
\ex{3.32} holds, and since the \rhs of \ex{3.33} satisfies Claim~3 of
Corollary\rw{cor.uniq}, then so does the \rhs of \ex{3.32}.\qed

\section*{Appendix}
\appendix
\nsection{Non-abelian connections in link complements}
\label{a1}

\def\AFLua{ \AFL{e^{2\pi i \ua}} }

In this appendix we will prove Theorem\rw{red.conn}. First, we have
to establish two simple lemmas.
\begin{lemma}
\label{la1}
Let $\cL$ be an algebraically connected
$L$-component link with admissible
indexing of its components.
Then there exist positive numbers
$\uc=\mtibas{c}{L}$ and positive integers $\un=\mtibas{n}{L-1}$ such
that $\AFLua\neq 0$ for any $\ua\in\Aspucna$.
\end{lemma}

\proof
We will prove the lemma by induction in $L$. If $L=1$, then the claim
of the lemma follows from \eex{1.10} and\rx{1.12}. Now suppose
that $L\geq 2$ and that the lemma holds for $(L-1)$-component links.
Let us expand
$\AFLut$ in powers of $t_L-1$. The first term is this expansion is
described by the Torres
formula\rx{1.15*}, so
\qq
&
\AFLut = \lrbc{ \pjoLmo t_j^{l_{jL}/2} - \pjoLmo t_j^{-l_{jL}/2} }
\AFbas{\cL\xrem{L}}{\ut\xrem{L}}
+ (t_L-1) F(\cL;\ut),
\label{a1.1}\\
&
\qquad F(\cL;\ut) \in \ZZutthi.
\nonumber
\qqq
Now we substitute $\ut = e^{2\pi i\ua} $ in \ex{a1.1}.
It is easy to see that there exists an area of the form\rx{1.69} in
which
\qq
\left| \lrbc{ \pjoLmo
e^{i \pi l_{jL}a_j} - \pjoLmo e^{-i \pi l_{jL}a_j} }
\AFbas{\cL\xrem{L}}{e^{2\pi i \ua\xrem{L}}} \right| > \left|
(e^{2\pi i a_L}-1) F(\cL;e^{2\pi i\ua}) \right|.
\label{a1.2}
\qqq
Since $\cL_L$ is algebraically connected to $\cL\xrem{L}$, then there
exists another area of the form\rx{1.69} in which
\qq
\pjoLmo
e^{i \pi l_{jL}a_j} - \pjoLmo e^{-i \pi l_{jL}a_j}
\neq 0.
\label{a1.3}
\qqq
Finally, according to the assumption of induction, there exists an
area of the form\rx{1.69} in which
\qq
\AFbas{\cL\xrem{L}}{e^{2\pi i \ua\xrem{L}}}
\neq 0.
\label{a1.4}
\qqq
According to \ex{a1.1}, the intersection of these three areas is the
area of the form\rx{1.69} in which $\AFLua\neq 0$.\qed

\def\lg{ \mathfrak{g} }
\def\lh{ \mathfrak{h} }
\def\Dp{ \Delta_+ }
\def\po{ \pi_1 }
\def\compL{ S^3\setminus \cL }
\def\Nr{ {\uN} }
\def\cLa#1{ \cL_{l(#1)} }
\def\xa#1{ x_{l(#1)} }
\def\fhom#1{ f_{#1} }
\def\fA{ \fhom{A} }
\def\fAl{ \fhom{A,\l} }

The second lemma says that the Alexander-Conway function measures an
obstruction to an irreducible deformation of a $U(1)$-reducible
connection. The proof is based on the observation that the Alexander
matrix gives a system of linear equations governing this deformation
in the linear approximation.

Let $G$ be a simple Lie group, $\lg$ be its Lie algebra. Let $H$ be a
maximal torus of $G$ and let
let $\lh$ be the corresponding Cartan subalgebra. We denote by
$\Dp\subset\lh^*$
a set of positive roots of $\lg$.

Let $\cL$ be an $L$-component link
in $S^3$. We denote by $\po$ the group of $\cL$ which is the
fundamental group of its complement: $\po = \po(\compL)$. A flat
connection $A$
of the trivial $G$-bundle over $\compL$ is determined by
the corresponding homomorphism
\qq
\fA:\;\po\rightarrow G
\label{a1.5}
\qqq
up to an overall conjugation. The homomorphism\rx{a1.5} is, of
course, generated by the holonomies of the flat connection along the
elements of $\po$. We call a flat connection $H$-reducible, if there
exists a maximal torus $H\subset G$ such that the
homomorphism\rx{a2.1} maps $\po$ into $H$:
\qq
\fA:\;\po\rightarrow H\subset G
\label{a1.5*}
\qqq
Since $H$ is abelian, then
such a homomorphism is determined by the images of the meridians of
the link components
\qq
\fA:\;\mer(\cL_j) \mapsto \exp(x_j),\qquad \ux\in\lh
\label{a1.6}
\qqq
up to the action of the Weyl group.

\begin{lemma}
\label{al1.2}
Consider a flat $H$-reducible connection in the $G$-bundle over the
link complement $\compL$, determined by \eex{a1.5*} and\rx{a1.6}.
Suppose that for
any positive root $\l\in \Dp$:
\begin{itemize}
\item[(1)]
if $L=1$, then
\qq
\APbas{\cL}{e^{(\l,\ux)}} \neq 0;
\label{a1.7*1}
\qqq
\item[(2)]
if $L\geq 2$, then
\qq
\AFbas{\cL}{e^{(\l,\ux)}} \neq 0,
\label{a1.7}
\qqq
and for at least one value of $j$
\qq
(\l,x_j)\neq 0.
\label{a1.7*}
\qqq
\end{itemize}
Here $(\cdot,\cdot)$ denotes a natural pairing between the elements
of $\lh$ and its dual $\lh^*$. Then any infinitesimal deformation of
this connection should also be $H$-reducible.
\end{lemma}

\proof
First of all, we recall the Wirtinger presentation of $\po$ (see,
\egt Chapter 3 of\cx{BZ}). The link $\cL$ is projected onto a plane,
so that each of its components is split by overcrossings into arcs.
We denote the number of arcs for $\cL_j$ as $N_j$. Then $\po$ is
generated by $\Nr=\sjoL N_j$ generators
$s_i$ which are meridians of the arcs. We introduce a function $l(i)$
such that $s_i$ is the meridian of an arc which is a part of
$\cLa{i}$. We can always choose a projection of $\cL$
in such a way that $\Nr$
equals the number of crossings in that projection.
Then the meridians $s_i$, as elements of $\po$,
are related by $\Nr$ equations
\qq
r_i = 1, \qquad 1\leq i\leq \Nr,
\label{a1.8}
\qqq
where the relators $r_i$ come from the crossings in the projection of
$\cL$ and have a form
\qq
r_i = s_i s_k^{-1} s_m s_k\quad\mbox{or}\quad
r_i = s_i s_k s_m^{-1} s_k^{-1},\qquad l(i)=l(m),
\label{a1.9}
\qqq
depending on the type of the $i$-th crossing.

A homomorphism $\fA$ describing a flat $H$-reducible
connection maps the elements $s_i$ as
\qq
\fA(s_i) = \exp( \xa{i} ).
\label{a1.10}
\qqq
Consider an infinitesimal deformation of this map
\qq
\fA(s_i) = \exp( \e_i ) \exp( \xa{i}),
\label{a1.11}
\qqq
where $\e_i$ are infinitesimal elements of $\lg$. The images of the
elements $s_i$ should satisfy the relations\rx{a1.9}. It is easy to
see that to the first order in $\e$ this means that for any relator
$r_i$,
\qq
\fA\lrbc{\del r_i\over \del s_i} (\e_i) +
\fA\lrbc{\del r_i\over \del s_k} (\e_k) +
\fA\lrbc{\del r_i\over \del s_m} (\e_m) = 0.
\label{a1.12}
\qqq
Here $(\del r / \del s)$ denote Fox derivatives (see
\egt\cx{BZ}, Chapter 9,B) and for a group element $g\in G$ and Lie
algebra element $\e\in\lg$, $g(\e)$ denotes the adjoint action of $g$
on $\e$.

Let us decompose $\e_i$ into sums of Cartan algebra elements and root
space elements. Since the flat connection $A$ is $H$-reducible, then
equations\rx{a1.12} do not mix differents spaces. Therefore, since we
investigate the existence of irreducible deformations, then we may
assume that all $\e_i$ belong to the same root space
$V_\l\subset \lg$ corresponding to a root $\l$:
\qq
\e_i = b_i \e,\qquad
b_i\in\IR,\quad\e\in V_\l.
\label{a1.13}
\qqq
Let $\fAl:\;\po\rightarrow \IC$ be a homomorphism such that
\qq
\fAl(s_i) = e^{(\l,\xa{i})},\qquad 1\leq i\leq \Nr
\qqq
($\fAl$ is well-defined, because $\fAl(r_i)=1$). Then
the system of
equations\rx{a1.12} for elements\rx{a1.13} can be rewritten as
\qq
\fAl\lrbc{\del r_i\over \del s_i} b_i +
\fAl\lrbc{\del r_i\over \del s_k} b_k +
\fAl\lrbc{\del r_i\over \del s_m} b_m = 0,
\qquad 1\leq i\leq \Nr.
\label{a1.14}
\qqq

If $L=1$, then according to Lemma~9.12 of\cx{BZ}, the matrix of
equations\rx{a1.14} is degenerate and has a solution of the form
\qq
b_i = b,\qquad 1\leq i\leq \Nr, \qquad b\in\IR.
\label{1.14*}
\qqq
This is the only solution, because, according to Theorem~9.10
of\cx{BZ}, the determinant of any $(\Nr-1)\times(\Nr-1)$ minor of the
matrix of\rx{a1.14} is proportional to $\APbas{\cL}{e^{(\l,\ux)}}$,
so it is not equal to zero in view of condition\rx{1.17*1}.

If $L\geq 2$, then
according to the proof of Proposition~9.16 of\cx{BZ}, the matrix
of equations\rx{a1.14} is degenerate and has a solution of the form
\qq
b_i = (e^{(\l,\xa{i})} - 1)\,b,\qquad 1\leq i\leq\Nr,\qquad b\in\IR.
\label{a1.15}
\qqq
This is the only solution, because the determinant of
an $(\Nr-1)\times(\Nr-1)$ minor of the matrix of\rx{a1.14} obtained
by removing the $i$-th column and the $j$-th row, is proportional to
$(e^{(\l,\xa{i})}-1)\AFbas{\cL}{e^{(\l,\ux)}}$, the coefficients
being the powers of $e^{(\l,\ux)}$, so at least some
minors are not equal to zero in view of conditions\rx{a1.7}
and\rx{a1.7*}.

The deformation\rx{a1.11} with $\e_i$ given by
either \ex{1.14*} or by \eex{a1.13}
and\rx{a1.15} is in fact $H$-reducible, because it is easy to see
that $\fA(s_i)$ commute with each other in the first order in $\e$
and therefore can be conjugated back into $H$ at this order.\qed

Let us introduce the following notation: if for two Lie group
elements $g_1,g_2\in G$ there exists $h\in G$ such that
$g_2= hg_1 h^{-1}$, then we denote this as $g_1 \sim g_2$.

Now we can prove the following
\begin{theorem}
\label{at1.1}
Let $\cL$ be an $L$-component algebraically connected link with
admissible indexing of its components. Then for a simple Lie group
$G$ there exist
positive numbers
$\uc=\mtibas{c}{L}$ and positive integers $\un=\mtibas{n}{L-1}$ such
that if $\ux=\mtibas{x}{L}$ is a list of $L$ elements of $\lh$ such
that for any positive root $\l\in\Dp$ we have $|(\l,\ux)|\in\Aspucna$,
then any flat connection A in the link complement $\compL$ which
satisfies the conditions
\qq
\fA(\mer(\cL_j)) \sim \exp(x_j)
\label{a1.16}
\qqq
is $H$-reducible.
\end{theorem}

\proof
We prove the theorem by induction in the number of components of
$\cL$.

Let $L=1$. The
moduli space of flat $G$-connections is an algebraic variety (which
comes, \egt from solving the Wirtinger equations\rx{a1.8} for the
images of $s_i$ under the homomorphism $\fA$. Therefore if the claim
of the theorem does not hold, then this means that the trivial
connection can be infinitesimally deformed into a non-$H$-reducible
connection. However, this contradicts Lemma\rw{al1.2} coupled with
\ex{1.10}.

Suppose that the theorem holds for $(L-1)$-component link and does
not hold for an $L$-component link $\cL$. Since the moduli space
of flat $G$-connections in $\compL$ is an algebraic variety, then
this means that for any $L-1$ positive numbers $\uc\xrem{L}$ and for
any $L-2$ positive integers $\un\xrem{L}$ there exist $L$ Cartan
subalgebra elements $\ux$ such that $x_L=0$ and for any positive root
$\l$ we have
$|(\l,\ux\xrem{l})|\in\Asp{\uc\xrem{L},\un\xrem{L}}{}$ and there
exists a flat connection $A$ in $\compL$ which satisfies\rx{a1.16}
and is either not $H$-reducible, or it is $H$-reducible, but it can
be infinitesimally deformed into a connection which is not
$H$-reducible. Since the theorem holds for $\cL\xrem{L}$, then there
exist $\uc\xrem{L}$ and $\un\xrem{L}$ for which only the first
possibility is impossible. According to Lemma\rw{la1}, there also
exist such $\uc\xrem{L}$ and $\un\xrem{L}$ for which\rx{a1.7} holds,
and so according to Lemma\rw{al1.2}, the second possibility is also
impossible. The intersection of these two areas
$\Asp{\uc\xrem{L},\un\xrem{L}}{}$ contains another area of the same
form in which the assumption of induction leads to contradiction.\qed

\pr{Theorem}{red.conn}
This theorem
is a particular case of Theorem\rw{at1.1} for
$G=SU(2)$. \qed

\nsection{The first polynomials}
\label{a2}

\def\qpt{ (q^{-m_2-l} - q^{-\gr_2})(q^{m_1-n+l} - 1) }
\def\qlmo{ \qemo{l} }
\def\oqamt{ \oqemo{-\gr_2} }
\def\slon{ \sum_{l=1}^n }
\def\Tiiu#1#2#3{ T_{#1,#2}^{(#3)} }
\def\Tiip#1#2{ \Tiiu{#1}{#2}{+} }
\def\qemo#1{ q^{#1} - 1 }
\def\oqemo#1{ (1-q^{#1}) }

The calculation of the polynomials $\Tjkpm$ of \eex{2.31},\rx{2.32} is
a relatively straightforward excercise. For example, one could first
present the interesting part of \rx{2.6} as an exponential
\qq
\lefteqn{
{ \plon \qpt \over \plon (\qlmo) }\,
q^{-m_2(\gr_1-1) + m_1 m_2 + \hlf n(n+1) } =
}
\label{a2.1}\\
&= &  {m_1 \choose n} \oqamt^n \,q^{-m_2\gr_1}
\exp \left[
\lrbc{ m_2 + m_1 m_2 + \hlf n(n+1)} \log(1+h)
\right.
\nonumber\\
&& \left.
+ \slon\log \lrbc{ (1+h)^{-m_2-l} - q^{-\gr_2} \over 1 -
q^{-\gr_2} } + \log \lrbc{ (1+h)^{m_1-n+l} - 1 \over m_1-n+l} +
\log\lrbc{ (1+h)^l -1 \over l} \right].
\nonumber
\qqq
The exponent should be expanded in powers of $h$. Since the
coefficients in this expansion have a polynomial dependence on $l$,
the sums over $l$ produce Bernoulli polynomials of $n$. After
calculating the exponential of the resulting power series in $h$, we
arrive at \ex{2.31}.

Here is the list of the first polynomials
$\Tjkp(m_1,m_2,n)$ which were obtained
in this way:
\qq
\Tiip{0}{0} & = & 1,
\label{a2.2}\\
\Tiip{0}{1} & = & \hlf (m_1 n + n + 2m_1m_2 + 2m_2),
\label{a2.3}\\
\Tiip{1}{0} & = & -\hlf (n +2  m_2 + 1),
\label{a2.4}\\
\Tiip{0}{2} & = & {1\over 24} \left(
(3m_1^2 + 5m_1+2)n^2 + (12 m_1^2 m_2 + m_1^2 + 24m_1 m_2 - 5m_1
+12m_2 - 6) n \right.
\label{a2.5}\\
&&\qquad \left.+ 12(m_1^2m_2^2 + 2m_1m_2^2 + m_2^2 - m_1m_2 -
m_2)\right),
\nonumber\\
\Tiip{1}{2} & = & - {1\over 12} \left(
(3m_1+1) n^2 + 3(4m_1m_2 + m_1+2m_2-1)n
\right).
\label{a2.6}\\
&&\qquad\left.
+ 2(6m_1m_2^2+3m_1m_2
+3m_2^2-3m_2-2) \right),
\nonumber\\
\Tiip{2}{2}&=&{1\over 24} \left( 3n^2 + (12m_2+5)n +
2( 6m_2^2+6m_2+1) \right).
\label{a2.7}
\qqq
The corresponding polynomials $\Tjkm$ can be either derived in a
similar way from \ex{2.32} or from the list of $\Tjkp$ through the
relation\rx{2.33}.

\nsection{The Symmetry Principle}
\label{c2}

We are going to show that the Symmetry Principle of R.~Kirby
and P.~Melvin\cx{KiMe} follows rather easily from the Main
Theorem.

\begin{theorem}[R.~Kirby, P.~Melvin]
Let $\cL$ be an $L$-component link in $S^3$. For a positive integer
$K$ and for an integer $j$, $1\leq j\leq L$ consider two sets of
positive integer numbers $\ual$, $\ualp$ such that
$1\leq \ual\leq K-1$ and
\qq
\a\p_j = \a_k\;\;\mbox{if $k\neq j$},\quad
\a\p_j = K - \a_j
\label{c.1}
\qqq
Then
\qq
\JuapL = (-1)^{\soknjL l_{jk} (\a_k - 1)} \JuaL.
\label{c.2}
\qqq
\end{theorem}

\proof
First, we will prove this theorem for algebraically connected links
with $L\geq 2$. Let $\cL$ be such a link. It is easy to see that we
can assume without loss of generality that its indexing is admissible
and $j\neq 1$.

According to \eex{1.45**0} and\rx{1.45**} written for the case of the
empty link $\cLp$
\qq
\JruapL & = &
{1\over h}\,q^{\lkLuap}\, q^{
\phlkL}\,
{1\over \AFLuqap}
\,\snzi \frpPubLLputAFxeap\,h^n
\nonumber
\\
& = &
{1\over h}\,q^{\lkLuap}\, q^{\phlkL}\,
\PLuqap\snzi F_n(q^{\a\p})\,h^n,
\label{c.3}
\qqq
where
\qq
F_n(\ut) = {1\over \PLut\AFLut}\,\frpPubLLputAFxe.
\label{c.4}
\qqq
Since, according to relations\rx{1.17},\rx{1.16} and\rx{1.45**2},
$F_n(\ut)$ are rational functions of $\ut$ and since
$q^{\a\p_j} = q^{-\a_j}$, then obviously
\qq
F_n(q^{\ua\p}) = F_n(q^{\uapp} ),
\label{c.5}
\qqq
where
\qq
\app_j = \a_k\;\;\mbox{if $k\neq j$},\quad
\app_j =  - \a_j.
\label{c.6}
\qqq
At the same time, \eex{1.9*4} and\rx{1.16} indicate that
\qq
q^{\lkLuap}\,\PLuqap =
-(-1)^{\soknjL l_{jk} (\a_k - 1)} q^{\lkLuapp}\,\PLuqapp.
\label{c.7}
\qqq
A combination of \eex{c.3},\rx{c.5} and\rx{c.7} indicates that
\qq
\JruapL = -(-1)^{\soknjL l_{jk} (\a_k - 1)}\JruappL.
\label{c.8}
\qqq

Let us write \ex{1.45*1} for the case of an empty link $\cLp$ and for
the colors $\ualp$
\qq
\serc{\sumuLum\JruapL}{\ualp} =
\JuapL.
\label{c.9}
\qqq
Assume that in this formula $\ualp$ are fixed positive numbers as
described by \ex{c.1} (so that expansion in the \lhs goes in fact
over the
powers of $h$ as described by \eex{1.40},\rx{1.41}). Then we can
substitute the expression\rx{c.8} for $\JuapL$ in the \lhs of
\ex{c.9}. If we also substitute there $-\mu_j$ instead of $\mu_j$ so
that $\umu\uapp$ becomes $\umu\ua$ and compare the result with
\ex{c.9} written for the colors $\ual$, then we arrive at the
relation\rx{c.2}.

Thus we proved \ex{c.2} for algebraically connected links with
$L\geq 2$. For an arbitrary non-empty link $\cL$ one can always find
a knot $\cLz$ such that the link $\cLzL$ is algebraically connected
and has at least two components. Then \ex{c.2} holds for $\cLz$. If
we set there $\a_0=1$, then we see that it also holds for $\cL$,
because
\qq
\Jbas{1,\ual}{\cLzL} = \JuaL.
\label{c.10}
\qqq
This completes the proof of the Symmetry Principle.\qed

\end{document}